\documentclass[11pt, reqno]{amsart}

\usepackage[usenames,dvipsnames]{color}
	
\usepackage[OT2, T1]{fontenc}
\usepackage{url}
\usepackage{amsmath}
\usepackage{array}
\usepackage{graphicx}
\usepackage{amssymb}
\usepackage{amsthm}
\definecolor{link}{RGB}{11,0,128}
\usepackage[colorlinks=true,citecolor=NavyBlue,linkcolor=Bittersweet,urlcolor=link]{hyperref}
\usepackage{hyperref}
\usepackage{paralist}
\usepackage{colonequals}			
\usepackage{color}
\usepackage{mathabx}			
\usepackage{amscd}
\usepackage[all,cmtip]{xy}			
\usepackage{sseq}				
\usepackage{verbatim}
\usepackage{parskip}
\usepackage{microtype}
\usepackage[in]{fullpage}
\usepackage{mathrsfs} 			
\usepackage{cleveref}
\usepackage{enumitem}
\usepackage{moreenum}			
\usepackage[alphabetic,lite]{amsrefs} 	
\usepackage{stmaryrd}			
\usepackage[foot]{amsaddr}
\usepackage{subfiles}
\usepackage{ragged2e}
\usepackage[bbgreekl]{mathbbol}     
\DeclareSymbolFontAlphabet{\mathbb}{AMSb}   

\newcommand{\gA}{\alpha}
\newcommand{\al}{\alpha}
\newcommand{\ssup}{>}

\newcommand{\bA}{\mathbb{A}}

\newcommand{\bC}{\mathbb{C}}

\newcommand{\bF}{\mathbb{F}}
\newcommand{\bG}{\mathbb{G}}

\newcommand{\bN}{\mathbb{N}}
\newcommand{\bP}{\mathbb{P}}
\newcommand{\bQ}{\mathbb{Q}}

\newcommand{\bZ}{\mathbb{Z}}

\newcommand{\bbB}{\mathbf{B}}

\newcommand{\cA}{\mathcal{A}}
\newcommand{\cB}{\mathcal{B}}

\newcommand{\cH}{\mathcal{H}}

\newcommand{\cM}{\mathcal{M}}

\newcommand{\cO}{\mathcal{O}}
\newcommand{\cP}{\mathcal{P}}

\newcommand{\cS}{\mathcal{S}}
\newcommand{\cT}{\mathcal{T}}

\newcommand{\cX}{\mathcal{X}}
\newcommand{\cY}{\mathcal{Y}}

\newcommand{\fb}{\mathfrak{b}}
\newcommand{\fc}{\mathfrak{c}}

\newcommand{\fg}{\mathfrak{g}}

\newcommand{\fm}{\mathfrak{m}}

\newcommand{\fp}{\mathfrak{p}}

\newcommand{\fs}{\mathfrak{s}}
\newcommand{\ft}{\mathfrak{t}}
\newcommand{\fu}{\mathfrak{u}}

\newcommand{\sB}{\mathscr{B}}

\newcommand{\sE}{\mathscr{E}}
\newcommand{\sF}{\mathscr{F}}
\newcommand{\sG}{\mathscr{G}}

\newcommand{\sI}{\mathscr{I}}

\newcommand{\sL}{\mathscr{L}}

\newcommand{\sO}{\mathscr{O}}

\newcommand{\sR}{\mathscr{R}}

\newcommand{\sV}{\mathscr{V}}

\newcommand{\sX}{\mathscr{X}}

\newcommand{\kg}{\mathfrak{g}}
\newcommand{\kt}{\mathfrak{t}}

\newcommand{\g}{\gamma}

\DeclareMathOperator{\Ad}{Ad}		
\DeclareMathOperator{\Br}{Br}			
\DeclareMathOperator{\Cauchy}{Cauchy}	
\DeclareMathOperator{\Char}{char}		
\DeclareMathOperator{\Coker}{Coker}	
\DeclareMathOperator{\Eq}{Eq}			
\DeclareMathOperator{\Gal}{Gal}		
\DeclareMathOperator{\GL}{GL}		
\DeclareMathOperator{\Hom}{Hom}		
\DeclareMathOperator{\im}{Im}			
\DeclareMathOperator{\Ker}{Ker}		
\DeclareMathOperator{\Lie}{Lie}		
\DeclareMathOperator{\Mat}{Mat}		
\DeclareMathOperator{\Null}{Null}		
\DeclareMathOperator{\ord}{ord}		
\DeclareMathOperator{\PGL}{PGL}		
\DeclareMathOperator{\Pic}{Pic}		
\DeclareMathOperator{\Res}{Res}		
\DeclareMathOperator{\rk}{rk}			
\DeclareSymbolFont{cyrletters}{OT2}{wncyr}{m}{n}
\DeclareMathSymbol{\Sha}{\mathalpha}{cyrletters}{"58}	
\DeclareMathOperator{\SL}{SL}			
\DeclareMathOperator{\Sp}{Sp}		
\DeclareMathOperator{\Spec}{Spec}		
\DeclareMathOperator{\Sym}{Sym}		
\DeclareMathOperator{\Tor}{Tor}		

\newcommand{\ad}{\mathrm{ad}}			
\newcommand{\ani}{{\mathrm{ani}}}		
\newcommand{\ce}{\colonequals}
\newcommand{\const}{\mathrm{const}}	
\newcommand{\cont}{{\mathrm{cont}}}		

\newcommand{\eps}{\epsilon}
\newcommand{\et}{\mathrm{\acute{e}t}}	
\newcommand{\fet}{\mathrm{f\acute{e}t}}	
\newcommand{\fin}{\mathrm{fin}}		
\newcommand{\heart}{{\heartsuit}}		
\newcommand{\hra}{\hookrightarrow}
\renewcommand{\i}{^{-1}}
\newcommand{\id}{\mathrm{id}}			

\newcommand{\isomto}{\overset{\sim}{\longrightarrow}}

\newcommand{\llb}{\llbracket}			
\newcommand{\llp}{(\!(}			
\newcommand{\ov}{\overline}
\providecommand{\p}[1]{\left(#1\right)}
\providecommand{\qxq}[1]{\quad\text{#1}\quad}
\providecommand{\qx}[1]{\quad\text{#1}}
\providecommand{\xq}[1]{\text{#1}\quad}

\DeclareSymbolFontAlphabet{\mathbbl}{bbold}

\newcommand{\ra}{\rightarrow}
\newcommand{\Ra}{\Rightarrow}
\newcommand{\red}{{\mathrm{red}}}		
\newcommand{\reg}{\mathrm{reg}}		
\newcommand{\rrp}{)\!)}			
\newcommand{\rrb}{\rrbracket}			
\newcommand{\rs}{\mathrm{rs}}		
\newcommand{\sh}{\mathrm{sh}}			

\newcommand{\surjects}{\twoheadrightarrow}
\newcommand{\tensor}{\otimes} 			
\newcommand{\tors}{\mathrm{tors}}		

\newcommand{\wh}{\widehat}
\newcommand{\wt}{\widetilde}
\newcommand{\xra}{\xrightarrow}

\newcommand{\csub[1]}{\centering\subsection{\text{#1} \nopunct} \hfill \justify}

\providecommand{\up}[1]{{\upshape(}#1{\upshape)}}
\providecommand{\uref}[1]{{\upshape\ref{#1}}}
\providecommand{\uS}{{\upshape\S}}
\providecommand{\f}[2]{\frac{#1}{#2}}

\renewcommand{\b}{\textbf}
\providecommand{\ucolon}{{\upshape:} }
\providecommand{\uscolon}{{\upshape;} }

\newcommand{\brems}{\begin{rems} \hfill \begin{enumerate}[label=\b{\thenumberingbase.},ref=\thenumberingbase]}

\newcommand{\erems}{\end{enumerate} \end{rems}}
\newcommand{\bexs}{\begin{exas} \hfill \begin{enumerate}[label=\b{\thenumberingbase.},ref=\thenumberingbase]}

\newcommand{\eexs}{\end{enumerate} \end{exas}}
\newcommand{\m}{\item}
\newcommand{\bsm}{\begin{smallmatrix}}
\newcommand{\esm}{\end{smallmatrix}}
\newcommand{\blem}{\begin{lemma}}
\newcommand{\elem}{\end{lemma}}

\newcommand{\bconj}{\begin{conj}}
\newcommand{\econj}{\end{conj}}
\newcommand{\bprob}{\begin{Problem}}
\newcommand{\eprob}{\end{Problem}}

\newcommand{\bq}{\begin{Q}}
\newcommand{\eq}{\end{Q}}
\newcommand{\benum}{\begin{enumerate}[label={{\upshape(\alph*)}}]}
\newcommand{\benuma}{\begin{enumerate}[label={{\upshape(\arabic*)}}]}
\newcommand{\benumr}{\begin{enumerate}[label={{\upshape(\roman*)}}]}
\newcommand{\eenum}{\end{enumerate}}
\newcommand{\bitem}{\begin{itemize}}
\newcommand{\eitem}{\end{itemize}}
\newcommand{\bc}{}
\newcommand{\bd}{\begin{defn}}
\newcommand{\ed}{\end{defn}}

\newcommand{\beg}{\begin{eg}}
\newcommand{\eeg}{\end{eg}}
\newcommand{\bex}{\begin{exa}}
\newcommand{\eex}{\end{exa}}

\newcommand{\bcl}{\begin{claim}}
\newcommand{\ecl}{\end{claim}}
\newcommand{\lab}{\label}

\newcommand{\x}{\text}

\newcommand{\q}{\quad}
\newcommand{\qq}{\quad\quad}
\newcommand{\qqq}{\quad\quad\quad}

\newcommand{\tst}{\textstyle}
\newcommand{\ba}{\begin{aligned}}
\newcommand{\ea}{\end{aligned}}
\newcommand{\be}{\begin{equation}}
\newcommand{\ee}{\end{equation}}
\newcommand{\bpf}{\begin{proof}}
\newcommand{\epf}{\end{proof}}
\newcommand{\bthm}{\begin{thm}}
\newcommand{\ethm}{\end{thm}}
\newcommand{\bprop}{\begin{prop}}
\newcommand{\eprop}{\end{prop}}
\newcommand{\bcor}{\begin{cor}}
\newcommand{\ecor}{\end{cor}}
\newcommand{\brem}{\begin{rem}}
\newcommand{\erem}{\end{rem}}
\newcommand*{\QED}{\hfill\ensuremath{\qed}}

\newcommand{\der}{\mathrm{der}}		

\usepackage{aliascnt}
\newaliascnt{numberingbase}{subsubsection}
\numberwithin{equation}{numberingbase}

\newtheoremstyle{thms}{0pt}{0pt}{\itshape}{}{\bfseries}{.}{ }{}
\theoremstyle{thms}
\newtheorem{conj}[numberingbase]{Conjecture}
\newtheorem{cor}[numberingbase]{Corollary}
\newtheorem{lemma}[numberingbase]{Lemma}
\newtheorem{lem}[numberingbase]{Lemma}
\newtheorem{prop}[numberingbase]{Proposition}
\newtheorem{Q}[numberingbase]{Question}
\newtheorem{thm}[numberingbase]{Theorem}
\newtheorem{theorem}[numberingbase]{Theorem}
\newtheorem{variant}[numberingbase]{Variant}
\newtheorem{bclaim}[equation]{Claim}

\newtheoremstyle{claims}{0pt}{0pt}{}{}{\itshape}{.}{ }{}
\theoremstyle{claims}
\newtheorem{claim}[equation]{Claim}

\newtheoremstyle{defs}{0pt}{0pt}{}{}{\bfseries}{.}{ }{}
\theoremstyle{defs}
\newtheorem{defn}[numberingbase]{Definition}

\newtheorem{eg}[numberingbase]{Example}

\newtheorem{rem}[numberingbase]{Remark}
\newtheorem*{rems}{Remarks}

\Crefname{claim}{Claim}{Claims}
\Crefname{bclaim}{Claim}{Claims}
\Crefname{conj}{Conjecture}{Conjectures}
\Crefname{cor}{Corollary}{Corollaries}
\Crefname{defn}{Definition}{Definitions}
\Crefname{eg}{Example}{Examples}
\Crefname{prop}{Proposition}{Propositions} 
\Crefname{Q}{Question}{Questions}
\Crefname{rem}{Remark}{Remarks}
\Crefname{thm}{Theorem}{Theorems}
\Crefname{variant}{Variant}{Variants}

\theoremstyle{thms}
\newtheorem{thm-tweak}[subsection]{Theorem}
\Crefname{thm-tweak}{Theorem}{Theorems}
\newtheorem{lemma-tweak}[subsection]{Lemma}
\Crefname{lemma-tweak}{Lemma}{Lemmas}
\newtheorem{cor-tweak}[subsection]{Corollary}
\Crefname{cor-tweak}{Corollary}{Corollaries}
\newtheorem{prop-tweak}[subsection]{Proposition}
\Crefname{prop-tweak}{Proposition}{Propositions} 
\newtheorem{conj-tweak}[subsection]{Conjecture}
\Crefname{conj-tweak}{Conjecture}{Conjectures} 

\theoremstyle{defs}
\newtheorem{defn-tweak}[subsection]{Definition}
\Crefname{defn-tweak}{Definition}{Definitions}
\newtheorem{eg-tweak}[subsection]{Example}
\Crefname{eg-tweak}{Example}{Examples}
\newtheorem*{rems-tweak}{Remarks}
\newtheorem{rem-tweak}[subsection]{Remark}
\Crefname{rem-tweak}{Remark}{Remarks}

\newtheoremstyle{subsection-tweak}
   {0pt}
   {0pt}%
   {}
   {}%
   {\bfseries}
   {}%
   {.5em}
   {\thmnumber{\@{#1}{}\@{#2}.}%
    \thmnote{~{\bfseries#3.}}}    
    
\theoremstyle{subsection-tweak}
\newtheorem{pp}[numberingbase]{}
\newcommand{\bpp}{\begin{pp}}
\newcommand{\epp}{\end{pp}}

\theoremstyle{subsection-tweak}
\newtheorem{pp-tweak}[subsection]{}

\makeatletter
\def\@tocline#1#2#3#4#5#6#7{
    \begingroup 
    \@ifempty{#4}{}{}

    \parindent\z@ \leftskip#3\relax \advance\leftskip\@tempdima\relax
    #5\hskip-\@tempdima
      \ifcase #1
       \or\or \hskip 2em \or \hskip 1em \else \hskip 3em \fi%
      #6\nobreak\relax
    \dotfill\hbox to\@pnumwidth{\@tocpagenum{#7}}\par
    \nobreak
    \endgroup
 }
 \def\l@section{\@tocline{1}{0pt}{1pc}{}{}}

\renewcommand{\tocsection}[3]{%
  \indentlabel{\@ifnotempty{#2}{\makebox[1.3em][l]{%
    \ignorespaces#1 \bfseries{#2}.\hfill}}}\bfseries{#3}
    \vspace{-5pt}}

\renewcommand{\tocsubsection}[3]{%
  \indentlabel{\@ifnotempty{#2}{\hspace*{-0.5em}\makebox[2.1em][l]{%
    \ignorespaces#1#2.\hfill}}}#3
    \vspace{-5pt}}

\makeatother 

\makeatletter 
\newcommand\appendix@section[1]{%
  \refstepcounter{section}%
  \orig@section*{Appendix \@Alph\c@section. #1}%
}
\let\orig@section\section
\g@addto@macro\appendix{\let\section\appendix@section}
\makeatother

\author{Alexis Bouthier}

\author{K\k{e}stutis \v{C}esnavi\v{c}ius}
\address{Institut de Math\'{e}matiques de Jussieu-PRG, Universit\'{e} Pierre et Marie Curie, 4 place Jussieu, 75005 Paris, France}
\address{CNRS, UMR 8628, Laboratoire de Math\'{e}matiques d'Orsay, Universit\'{e} Paris-Saclay, 91405 Orsay, France}
\email{alexis.bouthier@imj-prg.fr, kestutis@math.u-psud.fr}

\date{\today}

\def\UTFviii@defined#1{%
  \ifx#1\relax
      \PackageError{inputenc}{Unicode\space char\space\expandafter
                              \UTFviii@splitcsname\string#1\relax
                              \MessageBreak
                              not\space set\space up\space
                              for\space use\space with\space LaTeX}\@eha
  \else\expandafter
    #1%
  \fi
}

\makeatletter
\def\UTFviii@defined#1{%
  \ifx#1\relax
      ?%
  \else\expandafter
    #1%
  \fi
}

\makeatother


\begin{document}

\title{Torsors on loop groups and the Hitchin fibration}

\maketitle

\begin{abstract} In his proof of the fundamental lemma, Ng\^o established the product formula for the Hitchin fibration over the anisotropic locus. One expects this formula over the larger generically regular semisimple locus, and we confirm this by deducing the relevant vanishing statement for torsors over $R\llp t\rrp$ from a general formula for $\Pic(R\llp t \rrp)$. In the build up to the product formula, we present general algebraization, approximation, and invariance under Henselian pairs results for torsors, give short new proofs for the Elkik approximation theorem and the Chevalley isomorphism $\fg\!\sslash\! G \cong \ft/W$, and improve results on the geometry of the Chevalley morphism $\fg \ra \fg\!\sslash\!G$.  \end{abstract}

\renewcommand{\abstractname}{R\'{e}sum\'{e}}

\begin{abstract}
Dans sa preuve du lemme fondamental, Ng\^o \'{e}tablit une formule du 
produit au-dessus du lieu anisotrope. On s'attend \`{a} ce qu'une telle 
formule s'\'{e}tende au-dessus de l'ouvert g\'{e}n\'{e}riquement  r\'{e}gulier 
semisimple. Nous \'{e}tablissons cette formule en la d\'{e}duisant d'un 
r\'{e}sultat d'annulation de torseurs sous des groupes de lacets \`{a} partir 
d'une formule g\'{e}n\'{e}rale pour $\Pic(R\llp t\rrp)$.
Dans le cours de la preuve, nous montrons des r\'{e}sultats g\'{e}n\'{e}raux 
d'alg\'{e}brisation, d'approximation et d'invariance hens\'{e}lienne pour des 
torseurs ; nous donnons de nouvelles preuves plus concises du th\'{e}or\`{e}me 
d'alg\'{e}brisation d'Elkik et de l'isomorphisme de Chevalley $\fg\!\sslash\! G \cong \ft/W$ et am\'{e}liorons les \'{e}nonc\'{e}s sur la g\'{e}om\'{e}trie du morphisme de Chevalley $\fg \ra \fg\!\sslash\!G$.
  \end{abstract}

\hypersetup{
    linktoc=page,     
}
\renewcommand*\contentsname{}
\q\\
\vspace{-50pt}

\tableofcontents

\vspace{-20pt}

\section{Introduction}

\begin{pp-tweak}[The product formula for the Hitchin fibration]
A key insight in Ng\^{o}'s proof of the fundamental lemma in \cite{Ngo10} is to relate the affine Springer fibration, which over an equicharacteristic local field geometrically encodes the properties of orbital integrals, to the Hitchin fibration, which is global and whose geometric properties are easier to access. The mechanism that supplies the relation between the two is the product formula that Ng\^{o} established over the anisotropic locus $\cA^\ani$ of the Hitchin base $\cA$ in \cite{Ngo06}*{Th\'{e}or\`{e}me 4.6} and \cite{Ngo10}*{Proposition 4.15.1} and expected to also hold over the larger generically regular semisimple locus $\cA^\heart \subset \cA$ in \cite{Ngo10}*{before Corollaire 4.15.2}. The product formula over $\cA^\heart$ has already been used, for instance, in \cite{Yun11}*{Proposition~2.4.1},  \cite{Yun14}*{Section 5.5, equation (34)}, or \cite{OY16}*{proof of Proposition 6.6.3 (1)}, and one of the main goals of this article is to establish it in \Cref{thm:product-formula}. 

Roughly speaking, the product formula is a geometric incarnation of the Beauville--Laszlo glueing for torsors: it translates this glueing into geometric properties of the morphism of algebraic stacks that relates affine Springer fibers, which parametrize torsors over formal discs $R\llb t \rrb$, to Hitchin fibers, which parametrize torsors over a fixed proper smooth curve $X_R$ (for a variable base ring $R$). Under this dictionary, the product formula eventually reduces to a statement that torsors on $X_R$ are obtained from the ``Kostant--Hitchin torsor'' over a fixed open $U_R \subset X_R$ by glueing along the punctured formal discs $R\llp t \rrp$ at the $R$-points in $X \setminus U$. One is thus led to studying torsors over $R\llp t \rrp$. 

Over $\cA^\ani$ the Hitchin fibration is separated 
and the intervening stacks are Deligne--Mumford. For the product formula, these additional properties allowed Ng\^{o} to reduce to only considering those $R$ that are algebraically closed fields $k$, a case in which $k\llp t \rrp$ is a field with relatively simple arithmetic. Over $\cA^\heart$, however, such a reduction does not seem available, and we need to study more general~$R\llp t \rrp$.
\end{pp-tweak}

\numberwithin{equation}{subsection}
\begin{pp-tweak}[Torsors under tori over $R\llp t \rrp$]
The product formula says that the comparison morphism is a universal homeomorphism, so, due to the valuative criteria for stacks, the $R$ that are most relevant are fields and discrete valuation rings. Nevertheless, the valuative criterion for universal closedness assumes that the map is quasi-compact, so, to avoid verifying this assumption directly, it is convenient to allow more general $R$ (see \Cref{norm}). Our $R$ will in fact be seminormal, strictly Henselian, and local, and the key torsor-theoretic input to the product formula is then \Cref{twisted}: for such an $R$ and an $R\llp t \rrp$-torus $T$ that splits over a finite \'{e}tale Galois cover whose degree is invertible in $R$,
\be \label{eqn:key-torsor-input}
H^{1}(R\llp t\rrp,T) \cong 0.
\ee
Relative purity results from \cite{SGA4III}*{Expos\'{e} XVI}, whose essential input is the relative Abhyankar's lemma, reduce this vanishing to $T = \bG_m$. In this case, there is in fact a general formula 
\be \label{Gabber-thm-ann}
\Pic(R\llp t \rrp) \cong \Pic(R[t\i]) \oplus H^1_\et(R, \bZ)
\ee
due to Gabber \cite{Gab19} that is valid for any ring $R$ and is presented in the slightly more general setting of an arbitrary $R$-torus in \Cref{Gabber}. For seminormal $R$, we have $\Pic(R[t\i]) \cong \Pic(R)$, so if $R$ is also strictly Henselian local, then all the terms in \eqref{Gabber-thm-ann} vanish and \eqref{eqn:key-torsor-input} follows. 

In addition, the vanishing \eqref{eqn:key-torsor-input} implies that for a seminormal, strictly Henselian, local ring $R$ and any $n > 0$ less than any positive residue characteristic of $R$, every regular semisimple $n \times n$ matrix with entries in $R\llp t \rrp$ is conjugate to its companion matrix---see \Cref{conj}, which gives a general conjugacy to a Kostant section result of this type.
\end{pp-tweak}

Overall the argument for the product formula is fairly short---it suffices to read \S\S\ref{section:H1-description}--\ref{section:tame-tori}, \S\ref{section:Hitchin-fibration}, and review \S\ref{section:Chevalley-Kostant}---but we decided to complement it with the following improvements and generalizations to various broadly useful results that enter into its proof.

\begin{pp-tweak}[Algebraization of torsors and approximation]
A practical deficiency of the Laurent power series ring $R\llp t \rrp$ is that its formation does not commute with filtered direct limits and quotients in $R$, so one often prefers its Henselian counterpart $R\{t\}[\f 1t]$ reviewed in \S\ref{hens-series}. We show that such ``algebraization'' does not affect torsors: by \Cref{Hi-algebraize}, for any ring $R$ and any smooth, quasi-affine, $R\{t\}[\f 1t]$-group~$G$,
\be \label{torsor-inv}
\tst H^1(R\{t\}[\f 1t], G) \isomto H^1(R\llp t\rrp, G),
\ee
which generalizes a result of Gabber--Ramero \cite{GR03}*{Theorem 5.8.14} valid in the presence of a suitable embedding $G \hra \GL_n$. To prove \eqref{torsor-inv}, we exhibit a general procedure for showing that 
\[
\tst F(R\{t\}[\f 1t]) \isomto F(R\llp t\rrp)
\]
for functors $F$ that are invariant under Henselian pairs: the idea, which appears to be due to Gabber, is to consider the ring of $t$-adic Cauchy sequences (and double sequences) valued in $R\{t\}[\f 1t]$ and to show that this ring is Henselian along the ideal of nil sequences, see \Cref{filter-Hens} and \Cref{Gabber-thm}. To verify that our functor $F(-) = H^1(-, G)$ is invariant under Henselian pairs, we use recent results on Tannaka duality for algebraic stacks, see \Cref{cool-prop} and \Cref{Hens-pair-inv}.

The idea of considering Cauchy sequences and, more generally, Cauchy nets also leads to a new proof and a generalization of the the Elkik approximation theorem, for which we exhibit new non-Noetherian versions in \Cref{thm:basic-Elkik,thm:beyond-affine,Elk-Gab}. We then use them to extend the algebraization results to non-affine settings in \S\ref{section:beyond-affine}: for instance, we show that for a Noetherian ring $R$ that is Henselian along an ideal $J$ and the $J$-adic completion $\wh{R}$, 
\[
\Br(U) \isomto \Br(U_{\wh{R}}) \qxq{for every open} \Spec(R) \setminus V(J) \subset U \subset \Spec(R),
\]
a result that was announced in \cite{Gab93}*{Theorem 2.8~(i)}; see \Cref{Br-comp} for further statements of this sort and the results preceding it in \S\ref{section:beyond-affine} for sharper non-Noetherian versions. For a concrete situation in which such passage to completion is useful, see \cite{brauer-purity}*{Proposition 3.3 and the proof of Theorem 5.3}. 
\end{pp-tweak}

\begin{pp-tweak}[The Chevalley isomorphism and small characteristics]
The construction of the Hitchin fibration for a reductive group $G$ with Lie algebra $\fg$ rests on the Chevalley isomorphism 
\be \label{eqn:chev-iso-ann}
\fg\!\sslash\! G \cong \ft  /W,
\ee
 where $G$ acts on $\fg$ by the adjoint action, $\ft$ is the Lie algebra of a maximal torus $T \subset G$, and $W \ce N_G(T)/T$ is the Weyl group. In \Cref{adj}, we give a short proof for \eqref{eqn:chev-iso-ann} that is new even over $\bC$ but works over any base scheme $S$ as long as $G$ is \emph{root-smooth} (see \S\ref{pp:root-smoothness}; this condition holds if $2$ is invertible on $S$ or if the geometric fibers of $G$ avoid types $C_n$). The main idea is to consider the Grothendieck alteration 
 \[
 \wt{\fg} \surjects \fg
 \]
where $\wt{\fg}$ is the Lie algebra of the universal Borel subgroup of $G$, and to extend the $W$-action from the regular semisimple locus $\wt{\fg}^\rs$ to the maximal locus $\wt{\fg}^\fin = \wt{\fg}^\reg$ over which the alteration is finite. The result generalizes work of Chaput--Romagny \cite{CR10}, who adapted a classical proof to the case of a general base $S$ under more restrictive assumptions.
 
In \S\ref{section:Chevalley-Kostant}, we use the Chevalley isomorphism to review the constructions that go into setting up the product formula, such as building the group $J$ that descends to $\ft/W$ the centralizer of the universal regular section of $\fg$, and we take the opportunity to improve their assumptions: roughly, it suffices to assume that each residue characteristic of the base $S$ is not a torsion prime for the root datum of the respective fiber of $G$. This condition, described precisely in \S\ref{pp:torsion-primes}, is less restrictive than the order of the Weyl group $W$ being invertible on $S$, as is often assumed in \cite{Ngo10}, and is slightly weaker than the conditions that appear in \cite{Ric17}, so we improve several results in these references. A key advance that permits this is the construction of the Kostant section under less restrictive assumptions than before that was recently carried out  in \cite{AFV18}*{Section 2}. 
\end{pp-tweak}

\begin{pp-tweak}[Notation and conventions] \label{conv}
All our rings are commutative, with unit, except that sometimes we also use \emph{nonunital} rings, whose only distinctive feature is that they need not have a multiplicative identity (but are still commutative). Every ring is a nonunital ring and, more generally, so is every ideal in any nonunital ring. For a module $M$ over ring $A$ and an $a \in A$, we let $M\langle a\rangle$ denote the elements of $M$ killed by $a$ and set 
\[
\tst M\langle a^\infty\rangle \ce \bigcup_{n > 0} M\langle a^n\rangle.
\]
For brevity, we call the spectrum of an algebraically closed field a \emph{geometric point}. For a scheme $S$, we denote a choice of a geometric point above an $s \in S$ by $\ov{s}$. We say that $S$ is \emph{seminormal} if every universal homeomorphism $S' \ra S$ that induces isomorphisms on residue fields has a unique section (compare with \cite{SP}*{Lemma \href{https://stacks.math.columbia.edu/tag/0EUS}{0EUS}, Definition \href{https://stacks.math.columbia.edu/tag/0EUT}{0EUT}}, where the seminormalization of $S$ is defined as the initial object of the category of such universal homeomorphisms); by \cite{SP}*{Lemma~\href{https://stacks.math.columbia.edu/tag/0EUQ}{0EUQ}}, every seminormal $S$ is reduced. For a vector bundle $\sV$ on a scheme $S$, we often identify $\sV$ with the $S$-scheme $\underline{\Spec}_{\sO_S}(\Sym(\sV^\vee))$ whose functor of points is $S' \mapsto \Gamma(S', \sV \tensor_{\sO_S} \sO_{S'})$ (see \cite{SGA3Inew}*{Expos\'{e}~I, Corollaire 4.6.5.1}). Unless indicated otherwise, we form cohomology in the fppf topology, but whenever the coefficient sheaf is a smooth group scheme we implicitly make the identification \cite{Gro68c}*{Th\'{e}or\`{e}me 11.7} with \'{e}tale cohomology.

 We follow \cite{SGA3IIInew} for the basic theory of reductive group schemes, which, in particular, are required to have connected fibers (see the definition \cite{SGA3IIInew}*{Expos\'{e} XIX, D\'{e}finition 2.7}). For instance, we freely use the \'{e}tale local existence of splittings and pinnings (see \cite{SGA3IIInew}*{Expos\'{e}~XXII, D\'{e}finition 1.13, Corollaire 2.3; Expos\'{e} XXIII, D\'{e}finition 1.1}) or the classification of split pinned reductive groups by root data (see \cite{SGA3IIInew}*{Expos\'{e} XXV, Th\'{e}or\`{e}me 1.1}). We let $\sR(G)$ denote the root datum associated to a splitting of $G$ by \cite{SGA3IIInew}*{Expos\'{e} XXII, Proposition 1.14}; the choice of a splitting will not matter when we use $\sR(G)$. For a Lie algebra $\fg$ and an $a\in\kg$, we let 
 \[
 \ad(a) \colon \fg \ra \fg \qxq{be the map} x \mapsto [a,x].
 \]
  Similarly, we denote the adjoint action of a group $G$ on $\Lie(G)$ by $\Ad(-)$. We let $C_G(-)$ denote a centralizer subgroup of $G$, and we let $\mathrm{Cent}(-)$ denote the center.

For a scheme $S$ and an affine $S$-group $G$ acting on an affine $S$-scheme $X$, we let $X\!\sslash\! G$ denote the affine $S$-scheme given by the $\underline{\Spec}_{\sO_S}$ of the equalizer between the action and the inclusion of a factor maps (so that the coordinate rings of $X\!\sslash\! G$ are the rings of invariants):
\[
X\!\sslash\! G \ce \underline{\Spec}_{\sO_S}(\Eq( \sO_X \rightrightarrows  \sO_X \tensor_{\sO_S} \sO_G)).
\]
This and other quotient notation is abusive when used for left actions, for instance, in \eqref{eqn:chev-iso-ann}, where the quotients are of left actions. 
The construction of $X\!\sslash\! G$ commutes with flat base change in $S$, see \cite{Ses77}*{page 243, Lemma 2}. If $G$ is finite locally free over $S$, then we abbreviate $X\!\sslash\! G$ to $X/G$ because, by, for instance, \cite{Ryd13}*{Theorem 4.1, Definition 3.17}, it agrees with the coarse moduli space of the algebraic stack quotient $[X/G]$, which we always form in the fppf topology. Often $G$ will be a Weyl group scheme $W$ of a reductive group scheme; we recall from \cite{SGA3II}*{Expos\'{e} XII, Th\'{e}or\`{e}me 2.1} that such a $W$ is always finite \'{e}tale. We let $\times^G $ denote a contracted product, that is, the quotient of the product of an object with a right $G$-action and an object with a left $G$-action by the combined action given by $g\cdot (x, y) = (xg\i, gy)$. 
\end{pp-tweak}

\subsection*{Acknowledgements} 
We thank Ofer Gabber for many helpful interactions; as the reader will notice, this article owes a significant intellectual debt to his various unpublished results. We thank Laurent Moret-Bailly for numerous suggestions of improvements and for exhaustive comments, especially, for his ideas presented in \S\ref{pp:Cauchy} that streamlined and generalized much of \S\S\ref{section-invariance}--\ref{section-approximation}. We thank the referee for a very thorough report that significantly improved the manuscript. We thank Jingren Chi, Christophe Cornut, Philippe Gille, Ning Guo, Aise Johan de Jong, G\'{e}rard Laumon, Kei Nakazato, Ng\^{o} Bao 
Ch\^au, Sam Raskin, Timo Richarz, Simon Riche, and Peter Scholze for helpful conversations and correspondence. Some results in this article complement the work of the first named author with David Kazhdan and Yakov Varshavsky, and the first named author thanks them both for constant support and encouragement. The second named author thanks MSRI for support during a part of the preparation period of this article while in residence for the Spring semester of 2019 supported by the National Science Foundation under Grant No.~1440140. This project received funding from the European Research Council (ERC) under the European Union's Horizon 2020 research and innovation programme (grant agreement No.~851146).

\numberwithin{equation}{subsubsection}



\section{Approximation and algebraization} \lab{chapter-invariance}

Our first goal is the algebraization of torsors over loop groups $R\llp t \rrp$, see \Cref{Hi-algebraize}. For this, we consider rings of Cauchy sequences, which we discuss in \S\ref{section-invariance} in the general setting of Gabber--Ramero triples. The latter, in addition to supplying more general topologies, also simultaneously capture the recurrent framework of Henselian pairs---consequently, the aforementioned \Cref{Hi-algebraize} contains the more basic \Cref{Hens-pair-inv} as a special case. We combine the Cauchy sequence technique with Beauville--Laszlo glueing in \S\ref{section-approximation} to give a new proof of the Elkik approximation theorem, see \Cref{thm:Elkik-reproof}. The latter plays a role in extending the algebraization statements to nonaffine settings in \S\ref{section:beyond-affine}, see \Cref{Br-comp} for a concrete consequence.

\csub[Invariance of torsors under Henselian pairs and algebraization] \label{section-invariance}

The principal goal of this section is to show in \Cref{Hi-algebraize} that for any ring $R$ the Laurent power series and the Henselian Laurent power series rings $R\llp t \rrp$ and $R\{t \}[\f{1}{t}]$ (see \S\ref{hens-series}) possess the same collection of torsors under a given smooth, affine group scheme. The first step towards this is the invariance of such collections under Henselian pairs, which we obtain in \Cref{Hens-pair-inv}.

\bpp[Zariski and Henselian pairs] \label{Zar-Hens}
We recall that a pair $(A, I)$ consisting of a ring $A$ and an ideal $I \subset A$ is \emph{Zariski} if $I$ lies in every maximal ideal, that is, if $1 + I \subset A^\times$. A Zariski pair $(A, I)$ is \emph{Henselian} if it satisfies the Gabber criterion\footnote{For an earlier Henselianity criterion of this sort, a ``Newton's lemma,'' see \cite{Gru72}*{Corollaire I.3} or \cite{Gre69}*{Theorem~5.11}.} in the sense that  every polynomial
\be \label{eqn:Gabber-poly}
\tst \qq f(T) = T^{N}(T - 1) + a_{N}T^{N} + \cdots + a_1T + a_0 \qq \text{with} \q a_i \in I \q \text{and} \q N \ge 1
\ee
has a (necessarily unique) root in $1 + I$ (by \cite{Gab92}*{Proposition 1} or \cite{SP}*{Lemma \href{http://stacks.math.columbia.edu/tag/09XI}{09XI}}, this agrees with other definitions). As observed by Gabber, these are properties of $I$, so we say that the nonunital ring $I$ is \emph{Zariski} or \emph{Henselian}, respectively, and extend this terminology to pairs with nonunital $A$. This terminology is not entirely abusive because every nonunital ring $I$ appears as an ideal in a commutative ring, for instance, in the ring $\bZ \oplus I$ with multiplication $(z, i)(z', i') = (zz', zi' + z'i + ii')$. 

If $(A, I)$ is Zariski or Henselian, then, by \cite{SP}*{Lemmas \href{http://stacks.math.columbia.edu/tag/0DYD}{0DYD} and \href{http://stacks.math.columbia.edu/tag/09XK}{09XK}}, so is $(A', I'A')$ for any ideal $I' \subset I$ and any integral morphism $A \ra A'$ (such as a surjection). By \cite{SP}*{Lemmas \href{http://stacks.math.columbia.edu/tag/0EM6}{0EM6} and~\href{http://stacks.math.columbia.edu/tag/0CT7}{0CT7}}, the category of Zariski (respectively, Henselian) pairs is closed under filtered direct limits, inverse limits, and contains nilpotent thickenings, so $(A, IJ)$ is Henselian whenever $A \isomto \varprojlim_{m > 0} (A/IJ^m)$. The Zariskization and, by \cite{SP}*{Lemma \href{https://stacks.math.columbia.edu/tag/0AGV}{0AGV}}, also the Henselization of a Noetherian ring along any ideal is Noetherian. 

We will say that a functor $F$ defined on some subcategory of commutative rings is \emph{invariant under Zariski pairs} (respectively, is \emph{invariant under Henselian pairs}) if for every Zariski (respectively, Henselian) pair $(A, I)$ with both $A$ and $A/I$ in the domain of definition of $F$, we have
\[
F(A) \isomto F(A/I).
\]
\epp

The following instance of the Henselian pair formalism is particularly relevant for this article. 

\bpp[Henselian power series] \label{hens-series}
For a ring $R$, we let $R\{t \}$ denote the \emph{Henselian power series} ring over $R$, that is, the Henselization of $R[t]$ with respect to $(t)$. In comparison to the power series ring $R\llb t \rrb$, to which it admits the $t$-adic completion map $R\{ t \} \ra R\llb t \rrb$, this ring is better behaved in non-Noetherian settings: for instance, the functor $R \mapsto R\{ t \}$ commutes with filtered direct limits and also with quotients in the following sense: for any ideal $I \subset R$, we have 
\[
R\{t\}/I R\{t\} \cong (R/I)\{t \}
\]
(see \S\ref{Zar-Hens}). 
These properties persist to the \emph{Henselian Laurent power series} ring $R\{t\}[\f{1}{t}]$ that comes equipped with the map $R\{t\}[\f{1}{t}] \ra R\llp t \rrp$ to the usual Laurent power series ring $R\llp t \rrp \cong (R\llb t \rrb)[\f{1}{t}]$. For any ring $R$, the map 
\[
\tst R\{ t \} \ra R\llb t \rrb \qxq{is injective, and hence so is} R\{t\}[\f{1}{t}] \ra R\llp t \rrp;
\]
indeed, for Noetherian $R$ this follows from the Krull intersection theorem \cite{SP}*{Lemma \href{https://stacks.math.columbia.edu/tag/00IQ}{00IQ}} and, in general, $R$ is a filtered direct union of its finite type $\bZ$-subalgebras $R_i$, so
\be \lab{eqn:hens-series}
\tst R\{t \} \cong \bigcup_i R_i\{t\} \hra \bigcup_i R_i\llb t \rrb \subset R\llb t \rrb.
\ee
\epp

The argument of \Cref{Hens-pair-inv} will use the following general lemma.

\blem \lab{xexc-comp}
For a Henselian pair $(A, I)$ with $A$ Noetherian, if the geometric fibers of $A \ra \wh{A}$ with $\wh{A} \ce \varprojlim_{m > 0} (A/I^m)$ are regular,\footnote{By \cite{EGAIV2}*{Proposition 7.4.6}, this holds if the formal fibers of $A$ are geometrically regular, for instance, if $A$ is quasi-excellent.} then for any functor
\[
F \colon \text{$A$-algebras} \ra \text{Sets}
\]
that commutes with filtered direct limits,
\[
F(A) \hra F(\wh{A}) \qxq{and an element of} F(A/I) \qxq{lifts to} F(A) \qxq{if and only if it lifts to} F(\wh{A}).
\]
In particular, for a Noetherian ring $B$ whose formal fibers are geometrically regular, 
a functor 
\[
G \colon \text{$B$-algebras} \ra \text{Sets}
\]
that commutes with filtered direct limits is invariant under Henselian pairs if and only if it is invariant under those Henselian pairs that are obtained by completing a finite type $B$-algebra along  some ideal.
\elem

\bpf
By Popescu's smoothing theorem \cite{SP}*{Theorem \href{http://stacks.math.columbia.edu/tag/07GC}{07GC}}, there is a filtered direct system $\{A_j\}_{j \in J}$ of smooth $A$-algebras such that 
\[
\tst \wh{A} \cong \varinjlim_{j \in J} A_j, \qq \text{so that also} \qq F(\wh{A}) \cong \varinjlim_{j \in J} F(A_j).
\]
By \cite{Gru72}*{Th\'{e}or\`{e}me I.8} (which is the affine scheme case of \Cref{cool-prop} below), the smooth map $A \ra A_j$ has a retraction: the $A/I$-point of $A_j$ inherited from $\wh{A}$ lifts to an $A$-point. Thus, the map $F(A) \ra F(A_j)$ also has a retraction, so that $F(A) \hra F(\wh{A})$ and any lift of an element of $F(A/I)$ to $F(\wh{A})$ first descends to a lift in some $F(A_j)$ and then maps via the retraction to a desired lift~in~$F(A)$.

For the assertion about $G$, by a limit argument, we may assume that the Henselian pair $(A, I)$ in question is the Henselization of a finite type $B$-algebra along some ideal. By \S\ref{Zar-Hens}, such an $A$ is Noetherian and, by \cite{SP}*{Proposition \href{https://stacks.math.columbia.edu/tag/07PV}{07PV}, Lemma \href{https://stacks.math.columbia.edu/tag/0AH3}{0AH3}}, it inherits geometric regularity of formal fibers from $B$, so the claim about $F$ applied to $G$ allows us to replace $(A, I)$ by its $I$-adic completion and to conclude. 
\epf

We combine \Cref{xexc-comp} with recent results on Tannaka duality for algebraic stacks to extend the lifting result for smooth affine schemes \cite{Gru72}*{Th\'{e}or\`{e}me I.8} used above to algebraic stacks as follows.

\bprop \label{cool-prop}
For a Henselian pair $(A, I)$ and a smooth algebraic $A$-stack $\sX$ with quasi-affine diagonal, the pullback map $\sX(A) \ra \sX(A/I)$ is essentially surjective.\footnote{Added in proof: The surjectivity also holds when $\sX$ is a smooth scheme (with no assumption on its diagonal), see \cite{torsors-regular}*{Proposition 6.1.1~(a)}.}
\eprop

\bpf
By combining \cite{SP}*{Lemmas \href{https://stacks.math.columbia.edu/tag/04Y9}{04Y9}, \href{https://stacks.math.columbia.edu/tag/04XL}{04XL}, and \href{https://stacks.math.columbia.edu/tag/06FJ}{06FJ}, Definition \href{https://stacks.math.columbia.edu/tag/04YB}{04YB}}, we see that any $A/I$-point of $\sX$ factors through a quasi-compact open substack, so we lose no generality by assuming that $\sX$ is quasi-compact. The assumptions then ensure that $\sX$ is of finite presentation, so limit formalism for algebraic stacks \cite{LMB00}*{Proposition 4.18}
 reduces us to the case when $(A, I)$ is the Henselization of a finitely generated $\bZ$-algebra along some ideal. 
 By \S\ref{Zar-Hens}, this $A$ is Noetherian and, by \cite{SP}*{Lemmas \href{https://stacks.math.columbia.edu/tag/0AH2}{0AH2} and \href{https://stacks.math.columbia.edu/tag/0AH3}{0AH3}}, the geometric fibers of the map $A \ra \wh{A}$ to the $I$-adic completion are regular. Thus, due to limit formalism again, \Cref{xexc-comp} applies and reduces us to the case when $A$ is Noetherian and $I$-adically complete.
 
In this case, by a general continuity result \cite{HR19}*{Corollary 1.5 (ii)} for values of a Noetherian algebraic stack with quasi-affine diagonal (alternatively, see \cite{BHL17}*{Corollary 1.5}), the pullback morphism 
\[
\tst \sX(A) \ra \varprojlim_{n \ge 1} \sX(A/I^n) \q \x{is an equivalence of categories.}
\]
By the infinitesimal lifting criterion for smoothness \cite{SP}*{Lemma \href{https://stacks.math.columbia.edu/tag/0DP0}{0DP0}}, every $A/I$-point of $\sX$ extends to a compatible sequence of $A/I^n$-points, and the desired conclusion follows.
\epf

\beg \label{Elk-Gab-b}
For instance, $\sX$ could be a smooth, quasi-separated $A$-algebraic space (or even a scheme): by \cite{SP}*{Lemma \href{https://stacks.math.columbia.edu/tag/03HK}{03HK}}, the diagonal of a morphism of algebraic spaces is representable, separated, locally quasi-finite, so, by \cite{SP}*{Lemma \href{https://stacks.math.columbia.edu/tag/02LR}{02LR}}, its quasi-compactness implies quasi-affineness.
\eeg

The $H^1$ aspect of the following consequence of the work above was announced in \cite{Str83}*{Theorem~1} for smooth affine $G$ and was then proved in \cite{GR03}*{Theorem 5.8.14} for smooth $G$ that admit suitable embeddings into $\GL_n$. In turn, the $H^2$ aspect is a generalization of an unpublished result of Gabber, who established it in the case $G = \bG_m$.

\bthm \label{Hens-pair-inv}
Let $(A, I)$ be a Henselian pair and $G$ an $A$-group scheme.
\benum
\m \label{HPI-a}
If $G$ is 
smooth and quasi-separated, then 
\[
\q H^1(A, G) \hra H^1(A/I, G).
\]

\m \label{HPI-b}
If $G$ is quasi-affine, of finite presentation, and flat over $A$, then\footnote{Added in proof: For a similar result due to To\"{e}n in the case when $A$ is local Henselian and $G$ is no longer quasi-affine, see \cite{topology-torsors}*{Proposition B.13}.}
\[
\qqq H^1(A, G) \surjects H^1(A/I, G) \q \x{{\upshape(}respectively,} \q 
\Ker(H^2(A, G) \ra H^2(A/I, G)) = \{*\} \x{\upshape{)}.}
\]
\eenum
In particular, $\Br(A) \cong \Br(A/I)$ and if $G$ is quasi-affine and $A$-smooth, then 
\[
H^1(A, G) \isomto H^1(A/I, G).
\]
\ethm

\bpf
For \ref{HPI-a}, by \cite{SP}*{Lemma \href{https://stacks.math.columbia.edu/tag/04SK}{04SK}}, 
the functor of isomorphisms between two fixed $G$-torsors is representable by a smooth, quasi-separated $A$-algebraic space $X$. 
Thus, \Cref{cool-prop} ensures that 
$X(A) \surjects X(A/I)$, as desired.

In \ref{HPI-b}, we interpret the noncommutative $H^2$ in terms of gerbes $\sG$ bound by $G$ (for agreement with the derived functor cohomology when $G$ is commutative, see \cite{Gir71}*{Chapitre IV, Th\'{e}or\`{e}me~3.4.2~(i)}). By \cite{SP}*{Theorem \href{https://stacks.math.columbia.edu/tag/06FI}{06FI}}, the classifying stack $\bbB G$ that parametrizes $G$-torsors is algebraic and, by, for instance, \cite{topology-torsors}*{Lemma A.2 (b), Proposition A.3}, it is smooth, quasi-compact, with quasi-affine diagonal. General descent results for algebraic stacks \cite{SP}*{Theorem \href{https://stacks.math.columbia.edu/tag/06DC}{06DC}, Lemmas \href{https://stacks.math.columbia.edu/tag/0429}{0429} and \href{https://stacks.math.columbia.edu/tag/0423}{0423}} then imply that each $\sG$ is also algebraic, smooth, quasi-compact, with quasi-affine diagonal. Thus, both assertions of \ref{HPI-b} are immediate from \Cref{cool-prop}.

The assertion about $\Br(-)$ follows from the rest for $G = \GL_N$ and $G = \PGL_N$: indeed, by definition, 
\[
\tst \Br(A) = \bigcup_{N \ge 1} \im(H^1(A, \PGL_N) \ra H^2(A, \bG_m)_\tors), \qx{and likewise over $A/I$.} \qedhere
\]
\epf

\brem
For smooth, quasi-separated $G$, the injection 
\[
H^1(A, G) \hra H^1(A/I, G)
\]
of \Cref{Hens-pair-inv}~\ref{HPI-a} need not be surjective. Indeed, this may fail already for $G = \bZ$: if $A$ is a normal domain but $A/I$ is not, then $H^1(A, \bZ) = 0$, but $H^1(A/I, \bZ) \neq 0$ is possible (see \cite{Wei91}*{Remark~5.5.2}). 
\erem

\brem
Contrary to the assertion of the main theorem of \cite{Str84}, the injection 
\[
H^2(A, \bG_m) \hra H^2(A/I, \bG_m)
\]
of \Cref{Hens-pair-inv}~\ref{HPI-b} need not be surjective. Indeed, if it were, then for any regular ring $A$ and any ideal $I \subset A$, the group $H^2(A/I, \bG_m)$ would be torsion: one could lift cohomology classes to the Henselization of $A$ along $I$ and apply \cite{Gro68b}*{Corollaire 1.8}, which says that the values of $H^2(-, \bG_m)$ on regular rings are torsion abelian groups. However, this fails for 
\[
A \ce \bC[x, y, z] \q \x{and} \q I \ce (z(y^2 - (x^3 - x - z))) \subset\bC[x, y, z],
\]
for which $\Spec(A/I) \subset \bA^3_\bC$ is the union of two copies of $\bA^2_\bC$ whose intersection $C$ is the punctured elliptic curve $y^2 = x^3 - x$ in the $\{z = 0\}$ plane: indeed, the Mayer--Vietoris sequence \cite{Bou78}*{Chapitre IV, Corollaire 5.2} shows that $\Pic(C) \subset H^2(A/I, \bG_m)$, and $\Pic(C)$ has many nontorsion elements  given by restricting nontorsion elements of $\Pic^0(\ov{C}) \simeq \bC/\bZ$ to $C$, where $\ov{C}$ is the smooth compactification of $C$.
\erem

Our next goal is \Cref{Hi-algebraize}---the analogue of \Cref{Hens-pair-inv} that for any ring $R$ compares $G$-torsors over $R\{t\}[\f{1}{t}]$ and over $R\llp t\rrp$. The following convenient formalism of Gabber--Ramero introduced in \cite{GR03}*{Section 5.4} unifies this situation with the Henselian pair setting (see \Cref{eg:unify}).

\bpp[Gabber--Ramero triples] \label{GR-triples}
A \emph{Gabber--Ramero triple} is a datum $(A, t, I)$ of a commutative ring $A$, an element $t \in A$, and an ideal $I \subset A$. Such a triple is \emph{bounded} if, in the notation of \S\ref{conv},
\[
I\langle t^\infty \rangle = I\langle t^N \rangle \qxq{for some} N > 0,
\]
that is, if $I\langle t^\infty \rangle \cap t^N I = 0$ for some $N > 0$, equivalently, if the $t$-adic topology of $I$ induces the discrete topology on $I\langle t^\infty \rangle$. For instance, this happens when $t$ is a nonzerodivisor. In general, we set 
\[
\ov{A} \ce A/A\langle t^\infty\rangle \qxq{and} \ov{I} \ce I/I\langle t^\infty\rangle \cong I \ov{A} \subset \ov{A},
\]
so that the Gabber--Ramero triple $(\ov{A}, t, \ov{I})$ is bounded, in fact, $t$ is even a nonzerodivisor in $\ov{A}$. It is often useful to consider the intermediate bounded Gabber--Ramero triple 
\[
(\wt{A}, t, \ov{I}) \qxq{with} \wt{A} \ce A/I\langle t^\infty\rangle.
\]
The restriction of the surjection $\wt{A}  \surjects \ov{A}$ to each coset of $\ov{I}$ is injective, to the effect that this surjection is a local homeomorphism for the topologies that we are about to define.

We first endow the ideal $I$ with its $t$-adic topology, then $A$ with the unique ring topology for which the ideal $I \subset A$ is open, and, finally, $A[\f 1t]$ with the unique ring topology for which the map $A \ra A[\f 1t]$ is continuous and open \up{the existence of such ring topologies is ensured by the axioms \cite{BouTG}*{Chapitre~III, section 6, num\'{e}ro 3, axiomes (AV$_{\text{I}}$)--(AV$_{\text{II}}$)}}. Concretely, the resulting \emph{$(t, I)$-adic} ring topologies on $A$ and $A[\f{1}{t}]$ are determined by the respective open neighborhood bases of zero
\[
\tst \{t^m I\}_{m \ge 0} \q \x{and} \q \{\im(t^m I \ra A[\f{1}{t}])\}_{m \ge 0}.
\]
In particular, we have the identification of topological rings
\be \label{eqn:put-in-bar}
\tst A[\f 1t] \isomto \ov{A}[\f 1t].
\ee
We caution that the topology on $A[\f 1t]$ is not defined by ideals, in other words, it is not $A[\f 1t]$-linear, only $A$-linear. We let $\wh{A}$ and  $\wh{A[\f{1}{t}]}$ be the 
completions, so that, explicitly,
\be \label{comp-eq}
\tst \wh{A} \ce \varprojlim_{m > 0} (A/t^mI) \q \x{and} \q \wh{A[\f{1}{t}]} \cong \varprojlim_{m > 0} (A[\f{1}{t}]/\im(t^m I \ra A[\f{1}{t}])),
\ee
where the individual terms that appear in the last limit need not be rings, only their limit $\wh{A[\f{1}{t}]}$ is. By \eqref{eqn:put-in-bar}, the topological ring $\wh{A[\f{1}{t}]}$ depends only on the Gabber--Ramero triple $(\ov{A}, t, \ov{I})$. 

Both $\wh{A}$ and $\wh{A[\f{1}{t}]}$ are complete topological rings and $\wh{A}$ comes equipped with the ideals 
\[
\tst \wh{t^n I} \cong \varprojlim_{m \ge n} (t^nI/t^mI) \subset \wh{A} \q \x{that form a neighborhood base of zero.}
\]
By \cite{SP}*{Lemma \href{https://stacks.math.columbia.edu/tag/05GG}{05GG}}, we have $\wh{t^n I} = t^n \wh{I}$, so the topology of $\wh{A}$ is $(t, \wh{I})$-adic and 
\be \label{R-Rhat-comp}
A/t^nI \cong \wh{A}/t^n\wh{I}, \qxq{and hence also} A/(t^n) \cong \wh{A}/(t^n), \q \x{for every} \q n > 0.
\ee
The \emph{completion} of a Gabber--Ramero triple $(A, t, I)$ is the Gabber--Ramero triple $(\wh{A}, t, \wh{I})$. 

Maps $(A, t, I) \ra (A', t', I')$ of Gabber--Ramero triples are ring homomorphisms $f\colon A \ra A'$ that satisfy $f(t) = t'$ and $f(I) \subset I'$. Any such is continuous and induces continuous homomorphisms 
\[
\tst A[\f{1}{t}] \ra A'[\f{1}{t'}] \qxq{and} \wh{A[\f{1}{t}]} \ra \wh{A'[\f{1}{t'}]}.
\]
A common example is the map $(A, t, I) \ra (\wh{A}, t, \wh{I})$ to the completion. Other useful cases are when $(A', t', I')$ is a localization of $(A, t, I)$ or when $A = A'$ with $I \subset I'$.

A Gabber--Ramero triple $(A, t, I)$ is \emph{Zariski} (respectively, \emph{Henselian}) if so is the pair $(A, tI)$ (a weaker assumption than the same for $(A, t)$), that is, if so is the nonunital ring $tI$. For instance, by \S\ref{Zar-Hens}, the Gabber--Ramero triple $(\wh{A}, t, \wh{I})$ is always Henselian. 
\epp

\beg \label{eg:unify}
The unification alluded to before \S\ref{GR-triples} manifests itself as follows.
\benuma
\item \label{GRt-eg-1}
The case $t = 1$ amounts to that of a pair $(A, I)$ that is often assumed to be Henselian. In addition, $A[\f{1}{t}] \cong A$ and $\wh{A} \cong \wh{A[\f{1}{t}]} \cong A/I$, which is precisely the setting of \Cref{Hens-pair-inv}.

\item \label{GRt-eg-2}
The case $I = A$ amounts to that of a $t \in A$; a typical example: $A = R\{ t \}$ as above. The topology is then $t$-adic, one often assumes that $A$ has bounded $t^\infty$-torsion, so that $\wh{A[\f{1}{t}]} \cong \wh{A}[\f{1}{t}]$ by the following lemma, and the goal is to compare algebraic structures over $A[\f{1}{t}]$ and $\wh{A}[\f{1}{t}]$.
\eenum
\eeg

\blem \label{nzd-example}
For a bounded Gabber--Ramero triple $(A, t, I)$, 
\be \label{comp-GR}
\tst \wh{A}[\f{1}{t}] \isomto \wh{A[\f{1}{t}]} \q \x{as topological rings,} \q A\langle t^\infty\rangle \isomto \wh{A}\langle t^\infty\rangle, \q I\langle t^\infty\rangle \isomto \wh{I}\langle t^\infty\rangle
\ee
and $\ov{\wh{A}} \cong \wh{\ov{A}}$ in such a way that $\ov{t^n\wh{I}} 
= \wh{t^n\ov{I}}$ for $n > 0$ \up{in particular, compatibly with topologies}.
\elem

\bpf
In the case when $t$ is a nonzerodivisor in $A$, by \eqref{comp-eq}, it is also a nonzerodivisor in $\wh{A}$, and we only need to argue the first part of \eqref{comp-GR}. However, then we have the open continuous~injection  
\[
\tst \wh{A} \hra \wh{A[\f 1t]}
\]
of topological~rings, which becomes a topological isomorphism after inverting $t$, as desired.

In general, the boundedness assumption implies that $A\langle t^\infty \rangle \cdot t^mI$, the \emph{a priori} larger ideal $A\langle t^\infty \rangle \cap t^mI$, and the even larger ideal $t^m\cdot I\langle t^\infty\rangle$ all vanish for $m \ge N$. In particular, since 
\[
\Tor_1^A(A/I_1, A/I_2) \cong (I_1 \cap I_2)/(I_1 \cdot I_2)
\]
for any ideals $I_1, I_2 \subset A$, the system 
\[
\{ \Tor_1^A(\ov{A}, A/t^mI ) \}_{m \ge 0}
\]
is essentially zero. Thus, 
by forming $\varprojlim_{m \ge 0}$ over the exact sequences
\[
\Tor_1^A(\ov{A}, A/t^mI ) \ra A\langle t^\infty \rangle \tensor_A A/t^mI \ra A/t^mI \ra \ov{A}/t^mI\ov{A} \ra 0
\]
we obtain a short exact sequence: indeed, the transition maps of the system $\{A\langle t^\infty \rangle \tensor_A A/t^mI\}_{m \ge 0}$, and hence also the maps between the images of the second arrow, are surjective, whereas, by the previous sentence, the images of the first arrow form an essentially zero system, whose $\varprojlim_{m \ge 0}$ vanishes; the resulting vanishing of the $\varprojlim^1_{m \ge 0}$ terms then gives the exactness. We conclude that 
\[
A\langle t^\infty\rangle \cong \wh{A}\langle t^\infty \rangle, \qxq{that} \ov{\wh{A}} \cong \wh{\ov{A}}, \qxq{and, by replacing $A$ by $I$, that also} I\langle t^\infty\rangle \cong \wh{I}\langle t^\infty \rangle, \q \ov{\wh{I}} = \wh{\ov{I}}.
\]
By \eqref{R-Rhat-comp}, we have $\wh{t^n\ov{I}} \cong t^n \wh{\ov{I}}$ and, by \eqref{eqn:put-in-bar} and the settled case when $t$ is a nonzerodivisor, we have $\wh{A[\f{1}{t}]} \cong \wh{\ov{A}}[\f 1t]$, so we obtain the desired topological identification $\wh{A[\f{1}{t}]} \cong \wh{A}[\f{1}{t}]$.
\epf

The extension of \Cref{Hens-pair-inv} to the setting of Gabber--Ramero triples (so also to $R\{t\}[\f{1}{t}]$ and $R\llp t \rrp$) is a special case of the axiomatic algebraization theorem \ref{Gabber-thm} that generalizes an unpublished result of Gabber. Its proof given below rests on the following material on the Henselianity of rings of Cauchy sequences 
along their ideals of null sequences. 
To the best of our knowledge, the idea of considering rings of Cauchy sequences is due to Gabber. The idea of conceptualizing and generalizing our initial arguments by working with general filters and nonunital rings was suggested by Moret-Bailly.

\bpp[Rings of Cauchy nets] \label{pp:Cauchy}
We fix a nonunital topological ring $B$ and, for concreteness, assume that $B$ has an open neighborhood base of zero consisting of abelian subgroups. For a poset $S$ in which every two elements have a common upper bound,
we say that a function $f \colon S \ra B$
is \emph{null} 
(respectively, \emph{Cauchy}) if for every neighborhood $U \subset B$ of zero there is an $s \in S$ such that (respectively, the pairwise differences of) the values of $f$ on $S_{\ge s}$ lie in $U$. 
We consider the nonunital rings of functions
\[\ba
\Null(S, B) &\ce \{ \x{functions } f\colon S \ra B \x{ that are null}\}, \\
\Cauchy(S, B) &\ce \{ \x{functions } f\colon S \ra B \x{ that are Cauchy}\},
\ea\]
as well as the resulting nonunital rings of germs of functions
\[
\tst \Null_S(B) \ce \varinjlim_{s \in S} \Null(S_{\ge s}, B) \qxq{and} \Cauchy_S(B) \ce \varinjlim_{s \in S} \Cauchy(S_{\ge s}, B).
\]
Of course, $\Null(S, B)$ is an ideal in $\Cauchy(S, B)$ and $\Null_S(B)$ is an ideal in $\Cauchy_S(B)$. 
\epp

\blem \label{filter-Hens}
Let $S$ be a poset in which every two elements have a common upper bound, and let $B$ be a nonunital topological ring that has an open nonunital subring $B' \subset B$ 
whose induced topology has an open neighborhood base of zero consisting of ideals of $B'$. 
If $B'$ is Zariski \up{respectively, Henselian},  
then so is the nonunital ring $\Null_S(B)$. 
\elem

\bpf
Every null function eventually takes values in $B'$, so $\Null_S(B) \cong \Null_S(B')$ and we may assume that $B = B'$. 
The assumption on the topology then ensures that $\Null_S(B)$ is an ideal in 
\[
\tst \varinjlim_{s \in S}(\prod_{S_{\ge s}} B).
\]
Thus, since, by \S\ref{Zar-Hens}, being Zariski (respectively, Henselian) is stable under products, filtered direct limits, and passing to ideals, the claim 
follows. 
\epf

\brem
\Cref{filter-Hens} holds with the same proof for the `$N$-Henselian' property that is defined for Zariski pairs by only considering polynomials of a fixed degree $N$ in \eqref{eqn:Gabber-poly}. Similarly, instead of assuming that $B'$ be Henselian in \Cref{filter-Hens}, it is enough to assume that $B'$ has an $N$-Henselian open ideal for every $N \ge 0$, equivalently, an open neighborhood base of $0$ consisting of $N$-Henselian ideals (possibly with no single ideal being $N$-Henselian for all $N$). We will not pursue it, but this last improvement to \Cref{filter-Hens} leads to analogous improvements to results below. 
\erem

\bthm \label{Gabber-thm}
Let $B$ be a topological ring that has a Zariski \up{respectively, Henselian} open nonunital subring $B' \subset B$ whose induced topology has an open neighborhood base of zero consisting of ideals of $B'$, let $S$ be the poset given by some neighborhood base of zero of $B$ with the order relation given by $U \le U'$ if and only if $U' \subset U$, 
and consider a functor
\[
\tst F \colon \text{$B$-algebras} \ra \text{Sets}.
\]
\benum
\item \label{GT-a}
If $F$ satisfies 
\[
\tst \qq \varinjlim_{U \in S} F(B_U) \hra F(\varinjlim_{U \in S} B_U) \qxq{and} F(C) \hra F(C/J)
\]
for $S$-indexed direct systems of $B$-algebras $\{B_U\}$ and Zariski \up{respectively, Henselian} pairs $(C, J)$,~then
\[
\tst \qq F(B) \hra F(\wh{B}).
\]

\item \label{GT-b}
If $F$ satisfies 
\[
\tst \qq \varinjlim_{U \in S} F(B_U) \isomto F(\varinjlim_{U \in S} B_U) \qxq{and} F(C) \isomto F(C/J)
\]
for $S$-indexed direct systems of $B$-algebras $\{B_U\}$ and Zariski \up{respectively, Henselian} pairs $(C, J)$,~then
\[
\tst \qq F(B) \isomto F(\wh{B}).
\]
\eenum
\ethm


\bpf
By \cite{BouTG}*{Chapitre III, section 6, num\'{e}ro 5 and Chapitre III, section 7, num\'{e}ro 2, Corollaire~1} that identifies the completion $\wh{B}$ with the inverse limit of quotients of $B$ by open additive subgroups, the assumption on the topology of $B$ ensures that our poset $S$ has the property that
\be \label{eqn:Bhat-id}
\wh{B} \cong \Cauchy_S(B)/\Null_S(B).
\ee
Likewise, since the diagonal is cofinal in the product poset $S \times S$, we have 
\be \label{eqn:Bhat-id-2}
\tst \wh{B} \cong \Cauchy_{S \times S}(B)/\Null_{S \times S}(B).
\ee

\benum
\item
For every $U \in S$, the map given by the constant functions 
\[
\tst\qqq \const\colon B \ra \Cauchy(S_{\ge U}, B) \q \text{has a retraction} \q \mathrm{ev}_U\colon \Cauchy(S_{\ge U}, B) \ra B
\]
given by the evaluation at $U$. Thus, it induces an injection 
\[
\tst\qq F(B) \hra F(\Cauchy(S_{\ge U}, B)).
\]
It remains to form the direct limit over $U$ and combine \Cref{filter-Hens} with \eqref{eqn:Bhat-id}.

\item
To begin with, we claim that in any category, for a commutative diagram
\[
\qq \xymatrix{
X \ar[d]^-{i} \ar@{^(->}[r] & Y \ar[r] \ar@{^(->}[d] & X \ar[d]^-{i} \\
\wt{X} \ar@{^(->}[r] & \wt{Y} \ar[r] & \wt{X}
}
\]
 with monomorphisms as indicated and the bottom horizontal composition being the identity, the top composition is also the identity, $i$ is a monomorphism, and the left square is Cartesian; in addition, if the middle vertical monomorphism is split, then these assumptions, and hence also the conclusions, are preserved by any functor. Indeed, the functoriality is clear and, for the rest, the Yoneda embedding reduces us to the case of sets, when $X = \wt{X} \cap Y$ in $\wt{Y}$. 

For every $U \in S$, we apply the above to the~diagram
\[
\qqq\xymatrix@C=48pt{
B \ar@{^(->}[r]^-{\const} \ar[d]^-{\const} & \Cauchy(S_{\ge U}, B) \ar@{^(->}[d]^-{\mathrm{pr}_1^*} \ar[r]^-{\mathrm{ev}_U} & B \ar[d]^-{\const} \\
\Cauchy(S_{\ge U}, B) 
\ar@{^(->}[r]^-{\mathrm{pr}_2^*} & \Cauchy((S\times S)_{\ge (U,\, U)}, B) \ar[r]^-{\x{$\mathrm{ev}_{U \times \Delta}$}} 
& \Cauchy(S_{\ge U}, B) 
}
\]
in which the maps $\mathrm{pr}_i^*$ are obtained from the projections $\mathrm{pr}_i \colon S\times S \ra S$ and the map $\mathrm{ev}_{U\times \Delta}$ is given by the restriction of a Cauchy net indexed by $(S \times S)_{\ge (U,\,  U)}$ to the subnet indexed by $\{U \times U'\}_{U' \in S_{\ge U}}$.
 The monomorphism $\mathrm{pr}_1^*$ is split by the analogous map $\mathrm{ev}_{\Delta\times U}$. Thus, by first applying the functor $F$ and then the aforementioned claim, we obtain a Cartesian square
\[
\qqq \xymatrix@C=54pt{
F(B) \ar@{^(->}[r]^-{F(\const)} \ar@{^(->}[d]^-{F(\const)} & F(\Cauchy(S_{\ge U}, B)) \ar@{^(->}[d]^-{F(\mathrm{pr}_1^*)} \\
F(\Cauchy(S_{\ge U}, B)) 
\ar@{^(->}[r]^-{F(\mathrm{pr}_2^*)} & F(\Cauchy((S \times S)_{\ge (U,\, U)}, B)). 
}
\]
By forming the direct limit over $U$, we then obtain the Cartesian square 
\[
\qqq \xymatrix@C=54pt{
F(B) \ar@{^(->}[r]^-{F(\const)} \ar@{^(->}[d]^-{F(\const)} & F(\Cauchy_{S}(B)) \ar@{^(->}[d]^-{F(\mathrm{pr}_1^*)} \\
F(\Cauchy_{S}(B)) 
\ar@{^(->}[r]^-{F(\mathrm{pr}_2^*)} & F(\Cauchy_{S\times S}( B)). 
}
\]
Here the right vertical map is bijective: by \Cref{filter-Hens} and \eqref{eqn:Bhat-id}--\eqref{eqn:Bhat-id-2}, its source and target are compatibly identified with $F(\wh{B})$. Thus, the left vertical map is also bijective. Since \Cref{filter-Hens} and \eqref{eqn:Bhat-id} identify its target with $F(\wh{B})$, the conclusion follows. \qedhere
\eenum
\epf

The proof of \Cref{Gabber-thm}~\ref{GT-a} also shows the following variant.

\begin{variant} \label{variant}
Let $B$ and $S$ be as in Theorem \uref{Gabber-thm} and consider a functor
\[
\tst F \colon B\text{-algebras} \ra \text{Pointed sets}
\]
If $F$ satisfies 
\[
\tst  \Ker(\varinjlim_{U \in S} F(B_U) \ra F(\varinjlim_{U \in S} B_U)) =\{*\} \qxq{and} \Ker(F(C) \ra F(C/J)) = \{*\}
\]
for $S$-indexed direct systems of $B$-algebras $\{B_U\}$ and Zariski \up{respectively, Henselian} pairs $(C, J)$,~then
\[
\tst  \Ker(F(B) \ra F(\wh{B})) = \{*\}.
\]
\end{variant}

\brem
\Cref{Gabber-thm} and Variant \ref{variant} continue to hold with the same proof if $B$ is only a nonunital ring, where we interpret a ``$B$-algebra'' to mean a nonunital ring $C$ equipped with a morphism $B \ra C$ of nonunital rings (an abelian group homomorphism compatible with multiplication). However, this extension to the nonunital case does not formally specialize to the versions above when $B$ is unital because, due to the requirement that morphisms of commutative rings preserve the multiplicative unit, the categories of $B$-algebras differ in the two cases. 
\erem

\beg \label{filter-eg}
\Cref{Gabber-thm} and Variant \ref{variant} apply in the case when $B \ce A[\f 1t]$ for a Zariski (respectively, Henselian) Gabber--Ramero triple $(A, t, I)$ with $B'$ being the image of $tI$ and $\wh{B} = \wh{A[\f 1t]}$ (see \S\ref{Zar-Hens} and \S\ref{GR-triples}). In this case, the images of the $t^nI$ for $n > 0$ form a countable open neighborhood base of zero, so one may choose $S \ce \bN$ (if $t$ is a unit in $A$, then one may even choose $S$ to be a singleton) and \Cref{Gabber-thm} gives criteria for the following pullback to be injective or bijective:
\[
\tst F(A[\f 1t]) \ra F(\wh{A[\f 1t]}).
\] 
\eeg

Thanks to this example, the following three corollaries apply to Henselian Gabber--Ramero triples $(A, t, I)$ and show that the corresponding functors have the same values on $A[\f 1t]$ and $\wh{A[\f 1t]}$, for instance, that they have the same values on $R\{t\}[\f 1t]$ and $R\llp t \rrp$ for any ring $R$. They also apply in the context of rigid geometry, namely, to Henselian Huber rings that were defined in \cite{Hub96}*{Definition~3.1.2}.


\bcor \label{cor:stays-connected}
For a topological ring $B$ that has an open nonunital subring $B' \subset B$ that is Henselian and whose induced topology has an open neighborhood base of zero consisting of ideals of $B'$, the map $B \ra \wh{B}$ induces a bijection on idempotents. In particular, for any ring $R$, all the~maps~in
\[
\tst R \ra R\{t\} \ra R\llb t \rrb \qxq{and} R \ra R\{t\}[\f{1}{t}] \ra R\llp t \rrp \qx{induce bijections on idempotents,}
\]
so that $\Spec(R\llp t \rrp)$ is connected if and only if so is $\Spec(R)$.
\ecor

\bpf
The functor $F$ that sends a ring $C$ to the set of idempotents in $C$ commutes with filtered direct limits and, by \cite{SP}*{Lemma \href{https://stacks.math.columbia.edu/tag/09XI}{09XI}}, is invariant under Henselian pairs. Thus, \Cref{Gabber-thm}~\ref{GT-b} implies all the claims except for the assertion about $R \ra R\{t\}[\f{1}{t}]$. For the latter, we may first replace $R$ by $R^\red$ and then consider $R \ra R\llp t \rrp$ instead. It remains to note that for reduced $R$, by considering the term of lowest degree, the inclusion $R\llb t \rrb \hra R\llp t \rrp$ induces a bijection on idempotents.
\epf

The case of the map $R\{t\}[\f{1}{t}] \ra R\llp t \rrp$ is also of practical interest in the following special case.

\bcor \label{cor:finite-etale}
For a topological ring $B$ that has an open nonunital subring $B' \subset B$ that is Henselian and whose induced topology has an open neighborhood base of zero consisting of ideals of $B'$, pullback gives an equivalence between the categories of finite \'{e}tale algebras over $B$ and $\wh{B}$ and
\[
\tst R\Gamma_\et(B, \sF) \isomto R\Gamma_\et(\wh{B}, \sF) \qxq{for every torsion abelian sheaf} \sF \qxq{on} B_\et. 
\]
\ecor 

\bpf
By \cite{SP}*{Lemma \href{https://stacks.math.columbia.edu/tag/09ZL}{09ZL}}, the functor that associates to a ring the set of isomorphism classes of finite \'{e}tale algebras (respectively, the set of morphisms between fixed finite \'{e}tale algebras) is invariant under Henselian pairs. It also commutes with filtered direct limits, so \Cref{Gabber-thm}~\ref{GT-b} applies to give the claim about finite \'{e}tale algebras. For the rest, it suffices to similarly observe that, by \cite{SGA4II}*{Expos\'{e} VII, Corollaire 5.8}, for each $i \in \bZ$ the functor $B \mapsto H^i_\et(B, \sF)$ commutes with filtered direct limits and, by the affine analogue of proper base change \cite{Gab94}*{Theorem 1}, is invariant under Henselian pairs. 
\epf

The following example presents some functors that are well-behaved with respect to Zariski pairs.

\beg \label{sta-free-eg}
We recall that a module $M$ over a ring $C$ is \emph{stably free} if $M \oplus C^{\oplus n} \simeq C^{\oplus n'}$ for some $n, n' \ge 0$. The functors
\[ 
F \colon C \mapsto \{ \x{finite projective $C$-modules} \}\, /\simeq \qxq{and} F' \colon C \mapsto \{ \x{stably free $C$-modules} \}\, /\simeq
\]
satisfy
\[
F(C) \hra F(C/J) \qxq{and} F'(C) \isomto F'(C/J)  \qxq{for Zariski pairs}  (C, J).
\]
Indeed, if $M$ is a projective $C$-module and $M'$ is a $C$-module, then any $C/J$-morphism 
\[
f \colon M/JM \ra M'/JM'  \qxq{lifts to a $C$-morphism} \wt{f} \colon M \ra M',
\]
and Nakayama's lemma \cite{SP}*{Lemma \href{http://stacks.math.columbia.edu/tag/00DV}{00DV}} ensures that $\wt{f}$ is surjective whenever $M'$ is finitely generated and $f$ is surjective. Thus, since, by \cite{Mat89}*{Theorem 2.4}, a surjective endomorphism of a finite module is an isomorphism, we have $F(C) \hra F(C/J)$. Finally, to lift a stably free $C/J$-module to $C$, we note that any surjection 
\[
(C/J)^{\oplus n'} \surjects (C/J)^{\oplus n} \qxq{lifts to a necessarily split surjection} C^{\oplus n'} \surjects C^{\oplus n}.
\] 
\eeg

We are ready for the promised extension of \Cref{Hens-pair-inv} that includes \cite{GR03}*{Theorem 5.8.14} a special case.

\bcor \label{Hi-algebraize}
Let $B$ be a topological ring that has an open nonunital subring $B' \subset B$ that is Zariski and whose induced topology has an open neighborhood base of zero consisting of ideals~of~$B'$.
\benum
\item \lab{Hi-0}
Base change induces an injection
\[
\qq \{ \x{finite projective $B$-modules} \} /\!\simeq\,\, \hra \{ \x{finite projective $\wh{B}$-modules} \} /\!\simeq,
\]
as well as a bijection
\[
\qq \{ \x{stably free $B$-modules} \} /\!\simeq\,\, \isomto \{ \x{stably free $\wh{B}$-modules} \} /\!\simeq.
\]
 If $G$ is an inner form of $\GL_{n,\, B}$, then 
\[
\qq H^1(B, G) \hra H^1(\wh{B}, G).
\]

\item \label{Hi-00}
If $B'$ is Henselian and $G$ is a smooth, quasi-separated $B$-group scheme, then
\[
\tst \qq H^1(B, G) \hra H^1(\wh{B}, G).
\]

\item \lab{Hi-a}
If $B'$ is Henselian and $G$ is a quasi-affine, smooth $B$-group scheme, then
\[
\tst \qq H^1(B, G) \isomto H^1(\wh{B}, G).
\]

\item \lab{Hi-b}
If $B'$ is Henselian and $G$ is a quasi-affine, finitely presented, and flat $B$-group scheme, then
\[
\tst \qq \Ker(H^2(B, G) \ra H^2(\wh{B}, G)) =\{*\}.
\]
\eenum
In particular, for a ring $R$ and  a quasi-affine, smooth \up{respectively, quasi-affine, finitely presented, and flat} $R\{t\}[\f{1}{t}]$-group scheme $G$,
\[
\tst H^1(R\{t\}[\f{1}{t}], G) \isomto H^1(R\llp t \rrp, G) \q \x{{\upshape(}respectively,} \q \Ker(H^2(R\{t\}[\f{1}{t}], G) \ra H^2(R\llp t \rrp, G)) = \{*\}) .
\]
\ecor

\bpf
In \ref{Hi-0}, by \cite{Gir71}*{Chapitre III, Proposition 2.6.1~(i)}, up to translation inner forms have identical collections of torsors, so the claim about $G$ follows from the rest. 
For the rest of \ref{Hi-0}--\ref{Hi-b}, by \Cref{sta-free-eg}, \Cref{Hens-pair-inv}, and limit formalism (that we already discussed in the proofs of \Cref{cool-prop} and \Cref{Hens-pair-inv}), \Cref{Gabber-thm} and \Cref{variant} apply to the functors in question and give the claims. The `in particular' then follows from Examples \ref{filter-eg} and \ref{eg:unify}~\ref{GRt-eg-2}.
\epf

To further illustrate the Zariski aspects of \Cref{Gabber-thm}, in \Cref{Lam-cor}~\ref{LC-d} we reformulate the following conjecture of Lam that appeared in \cite{Lam78}*{equation (H) on page XI} and \cite{Lam06}*{equation~(H$'$) on page 180}.


\bconj[Lam]
For a local ring $R$, every stably free $R[t]$-module is free.
\econj


\bthm \label{Lam-cor}
Let $R$ be a ring and set $R(t) \ce (R[t]_{1 + tR[t]})[\f{1}{t}]$, so that $R(t)$ can be identified with the localization of $R[t]$ with respect to the monic polynomials \up{see the proof for this identification}.
\benum
\item \label{LC-b}
Nonisomorphic finite projective $R$-modules cannot become isomorphic over $R\llp t \rrp$.

\item \label{LC-a}
Nonisomorphic finite projective $R(t)$-modules cannot become isomorphic over $R\llp t \rrp$.

\item \label{LC-c}
Base change identifies the set of isomorphism classes of stably free $R(t)$-modules with the set of isomorphism classes of stably free $R\llp t \rrp$-modules.

\item \label{LC-d}
If $R$ is local, then every stably free $R[t]$-module is free iff every stably free $R\llp t \rrp$-module~is~free. 
\eenum
In particular, the following pullback maps are injective for $n \ge 0$\ucolon 
\[
H^1(R, \GL_n) \hra H^1(R(t), \GL_n) \hra H^1(R\llp t \rrp, \GL_n).
\]
\ethm

\bpf
Firstly, to exhibit the claimed identification of $R(t)$, we let $R(t)'$ be the localization  of $R[t]$ with respect to the multiplicative set of the monic polynomials, view $R(t)'$ as a localization of $R[t, t\i]$, and then note that the involution $t \mapsto t\i$ exchanges $R(t)$ and $R(t)'$.

Parts \ref{LC-a} and \ref{LC-c} are immediate from \Cref{Hi-algebraize}~\ref{Hi-0}. Part \ref{LC-b} follows by combining \ref{LC-a} with \cite{Lam06}*{Chapter V, Proposition 2.4}, which is the version of \ref{LC-b} in which $R\llp t \rrp$ is replaced with $R(t)$.  Part \ref{LC-d} follows from \ref{LC-c} and \cite{BR83}*{Theorem A}, which is the version of \ref{LC-d} in which $R\llp t \rrp$ is replaced with $R(t)$.\footnote{The blanket Noetherianity assumption of \emph{op.~cit.} is not needed for \cite{BR83}*{Theorem A}: indeed, both \cite{BR83}*{Theorem 2.2} and \cite{Lam06}*{Chapter IV, Horrocks' Theorem 2.1}, which give the two implications, are susceptible to limit arguments.  }
\epf

\brem
\Cref{Lam-cor} \ref{LC-c} and \ref{LC-d} continue to hold even if one drops the finite generation implicit in the definition of stable freeness: more precisely, by Gabel's trick \cite{Lam06}*{Chapter I, Proposition 4.2}, for a ring $A$, every $A$-module $M$ that is not finitely generated and such that $M \oplus A^{\oplus n}$ is free is itself free.
\erem


\csub[The Elkik--Gabber--Ramero approximation via Cauchy nets] \lab{section-approximation}

The Cauchy net technique used in \S\ref{section-invariance} leads to non-Noetherian versions of the Elkik approximation theorem via a short argument that is also new in the Noetherian case, see \Cref{thm:basic-Elkik,thm:beyond-affine}. It also strengthens the non-Noetherian version of this theorem presented in \cite{GR03}*{Proposition~5.4.21}: in \Cref{Elk-Gab} below, $t$ need not be a nonzerodivisor and the open $U$ need not be $\Spec(A[\f{1}{t}])$. In some sense, we invert the argument: Gabber and Ramero deduced their version from the Noetherian case that was settled by Elkik in \cite{Elk73} and restated in \cite{GR03}*{Lemma 5.4.12}, whereas we first settle the general non-Noetherian version and then deduce the Noetherian statement from it in \Cref{thm:Elkik-reproof}.

 
In spite of its slightly different flavor, Elkik approximation is spiritually close to \S\ref{section-invariance}: for instance, the approximation statements could be considered as nonabelian incarnations of the phenomenon that $H^1_{\{t = 0\}}$ tends to depend only on the formal $t$-adic neighborhood of the ring in question, and hence should not change upon passage to $t$-adic completions. For a concrete such statement, see \Cref{eg:cute-affine-Gras}, which ensures that the affine Grassmannian may be formed with Henselian loops. 

Our core approximation idea is captured by \Cref{thm:basic-Elkik}, formulated using the following topologies.

\bpp[Topology on points of affine schemes] \label{pp:affine-top}
We recall from \cite{Con12b}*{Proposition 2.1} that for any topological ring $B$ there is a unique way to topologize the sets $X(B)$ for affine $B$-schemes $X$ of finite type in such a way that
\benuma
\item \label{top-1}
any $B$-morphism $X \ra X'$ induces a continuous map $X(B) \ra X'(B)$;

\item \label{top-2}
for each $n \ge 0$, the identification $\bA^n(B) \cong B^n$ is a homeomorphism; 

\item \label{top-3}
a closed immersion $X \hra X'$ induces an embedding $X(B) \hra X'(B)$.
\eenum

Indeed, one chooses a closed immersion $X \hra \bA^n_{B}$ and checks (\emph{loc.~cit.})~that the resulting subspace topology on $X(B) \subset B^n$ does not depend on the choice. In terms of any such embedding, elements of $X(B)$ are topologically close if and only if the resulting values in $B$ of the corresponding standard coordinates of $\bA^n_{B}$ are close. Moreover, by \cite{Con12b}*{Proposition 2.1, Example 2.2} and the definitions,
\benuma \addtocounter{enumi}{3}
\item \label{top-7}
the identifications $(X \times_{X''} X')(B) \isomto X(B) \times_{X''(B)} X'(B)$ are homeomorphisms;

\item \label{top-4}
for a continuous homomorphism $B \ra B'$ of topological rings, the maps $X(B) \ra X(B')$ are continuous, for instance, for an open ideal $I \subset B$, the fibers of $X(B) \ra X(B/I)$ are open;

\m \label{top-new}
if a continuous homomorphism $B \ra B'$ of topological rings is a (respectively, open; respectively, closed; respectively, discrete) embedding, then so are the maps $X(B) \ra X(B')$ of topological spaces because   
\be\label{top-emb}
\qqq X(B) = X(B') \cap \bA^n(B) \qxq{in} \bA^n(B') \qxq{for a closed immersion} X \hra \bA^n;
\ee

\m \label{top-new-2}
if the restriction of a continuous homomorphism  $B \ra B'$ of topological rings to an open ideal $I \subset B$ is a (respectively, open) topological embedding, then the maps $X(B) \ra X(B')$ are locally on their sources (respectively, open) topological embeddings because \eqref{top-emb} holds locally on $\bA^n(B)$.
\eenum
\epp

\bthm \label{thm:basic-Elkik}
Let $B$ be a topological ring that has an open nonunital subring $B' \subset B$ that is Zariski and whose induced topology has an open neighborhood base of zero consisting of ideals of $B'$, and let $X$ be a smooth, affine $B$-scheme. 
If either 
\benumr
\m \label{BE-i}
$B'$ is Henselian\uscolon or

\m \label{BE-ii}
$X(C) \surjects X(C/J)$ for every Zariski $B$-pair $(C, J)$ with $C/J \cong \wh{B}$\uscolon\!\!\footnote{\label{foot:surj}The surjectivity assumption holds for every Zariski $B$-pair $(C, J)$ if $X$ is an fpqc inner form of $\GL_{n, B}$. Indeed, such an $X$ is the open subscheme cut out in the corresponding inner form $\wt{X}$ of the matrix algebra $\Mat_{n \times n, B}$ by the nonvanishing of the inner form of the determinant; thus, since $\wt{X}(C) \ra \wt{X}(C/J)$ is identified with the reduction modulo $J$ surjection for a projective $C$-module and $C^\times$ is the preimage of $(C/J)^\times$, we obtain the desired $X(C) \surjects X(C/J)$.} 
\eenum
then the pullback map
\[
X(B) \ra X(\wh{B}) \qx{has dense image.} 
\]
\ethm

\beg \label{eg:basic-eg}
For instance, $B' \subset B$ could be a Henselian pair $J \subset C$ with $J$ endowed with its coarse topology. In this case, the completion $\wh{C} \cong C/J$ is discrete and \Cref{thm:basic-Elkik} recovers \cite{Gru72}*{Th\'{e}or\`{e}me I.8}: for a smooth, affine $C$-scheme $Y$, we have 
\[
Y(C) \surjects Y(C/J).
\]
This case is an input to the proof of \Cref{thm:basic-Elkik}, see also \Cref{cool-prop} for a generalization.
\eeg

\bpf
We use the notation of rings of Cauchy nets introduced in \S\ref{pp:Cauchy} and let $S$ be the poset consisting of the neighborhoods of zero $U \subset B$ where $U \le U'$ if and only if $U' \subset U$, so that, as in \eqref{eqn:Bhat-id}, we have the identification
\[
\wh{B} \cong \Cauchy_S(B)/\Null_S(B).
\]
\Cref{filter-Hens} and \Cref{eg:basic-eg} 
then give us the surjection
\[
\tst \varinjlim_{U \in S} \p{X\p{\Cauchy(S_{\ge U}, B)}} \cong X\p{\Cauchy_S(B)} \surjects X(\wh{B}), 
\] 
so that every element of $X(\wh{B})$ comes from some $X(\Cauchy(S_{\ge U}, B))$. The elements of this last set may be visualized as $S_{\ge U}$-indexed posets of elements of $X(B)$ such that the values in $B$ of affine coordinates of $X$ converge to elements of $\wh{B}$. Thus, by considering a finite set that generates the coordinate ring of $X$ as an $B$-algebra, we conclude  that each neighborhood of every element of $X(\wh{B})$ meets the image of $X(B)$ (see \S\ref{pp:affine-top}), and the claimed density follows. 
\epf

We are going to extend \Cref{thm:basic-Elkik} in two directions: beyond affine $X$ in \Cref{thm:beyond-affine} and beyond affine $B$ in the context of Gabber--Ramero triples in \Cref{Elk-Gab}. To remove the affineness assumption on $X$, it is important to be able to extend the definition of the topology on $X(B)$. The standard approach to this is to assume that $B$ is local, so that every $B$-point of any $X$ factors through an affine open. We do not wish to restrict to local $B$---for instance, $R\llp t \rrp$ is rarely local---so in \S\ref{xI-top} we pursue the approach based on \Cref{GR-lemma}. For the sake of brevity, we use the following notion that avoids a quasi-compactness assumption inherent in quasi-projective morphisms.

\bd
A scheme morphism $X \ra S$ is \emph{subprojective} if $X$ is an open subscheme of a projective $S$-scheme, so that, in particular, $X$ is locally of finite type and separated over $S$. 
\ed

We also use the following generalization of quasi-affineness that avoids quasi-compactness assumptions. This definition seems to be due to Gabber and is given in  \cite{SP}*{Definition \href{https://stacks.math.columbia.edu/tag/0AP6}{0AP6}}.

\bd \label{def:ind-qa}
A scheme $X$ is \emph{ind-quasi-affine} if every quasi-compact open of $X$ is quasi-affine. A scheme map $f \colon X \ra S$ is \emph{ind-quasi-affine} if $f\i(U)$ is ind-quasi-affine for every affine open $U \subset S$.
\ed

An ind-quasi-affine morphism is separated because the intersection of any two affine opens in an ind-quasi-affine scheme is affine.
By \cite{SP}*{Lemma \href{https://stacks.math.columbia.edu/tag/0AP9}{0AP9}}, separated, locally quasi-finite morphisms, such as immersions or separated \'{e}tale maps, are ind-quasi-affine. By \cite{SP}*{Lemmas \href{https://stacks.math.columbia.edu/tag/0F1V}{0F1V}, \href{https://stacks.math.columbia.edu/tag/0AP7}{0AP7}, and~\href{https://stacks.math.columbia.edu/tag/0AP8}{0AP8}}, ind-quasi-affineness is stable under composition, base change, is~fpqc~local on the base. By \cite{SP}*{Lemma \href{https://stacks.math.columbia.edu/tag/0APK}{0APK}} that is due to Gabber, fpqc descent is effective for ind-quasi-affine morphisms.

\blem \label{GR-lemma}
For a ring $B$ and a $B$-scheme $X$ that is either ind-quasi-affine or subprojective, every $B$-point of $X$ factors through some affine open of $X$. 
\elem

\bpf
Every $B$-point factors through a quasi-compact open subscheme of $B$, so in the ind-quasi-affine case we may assume that $X$ is quasi-affine, in particular, that $X = \Spec(C) \setminus V(J)$ for a $B$-algebra $C$ and an ideal $J \subset C$. The $B$-point $C \surjects B$ in question maps $J$ to the unit ideal, so there is a linear combination $\sum_i c_ij_i$ with $c_i \in C$ and $j_i \in J$ such that the $B$-point factors through the affine open 
\[
\tst \Spec(C[\f{1}{\sum_i c_ij_i}]) \subset X. 
\]
The subprojective case is \cite{GR03}*{Lemma 5.4.17}: similarly to the proof of \Cref{find-affine} below, one constructs a hyperplane not meeting the $B$-point in question, and hence reduces to the ind-quasi-affine case. 
\epf


\bpp[Topology on points beyond affine schemes] \label{xI-top}\label{extend-top}
Let $B$ be a topological ring such that the unit group $B^\times \subset B$ is open and the inversion map of $B^\times$ is continuous for the subspace topology. 
By the proofs of \cite{Con12b}*{Proposition 3.1} and of \Cref{GR-lemma} above, these assumptions ensure that an open immersion of affine $B$-schemes of finite type induces an open embedding on $B$-points for the topology defined in \S\ref{pp:affine-top}. \Cref{GR-lemma} then allows us to extend the definition of this topology beyond affine $X$, namely, to topologize the set $X(B)$ for every locally of finite type $B$-scheme $X$ that is either ind-quasi-affine or subprojective by declaring a subset to be open if and only if its intersection with $U(B)$ is open for every affine open $U \subset X$. With this definition, thanks to \Cref{GR-lemma},
\benuma \addtocounter{enumi}{7}
\item \label{top-6}
an open immersion $X \hra X'$ induces an open embedding $X(B) \hra X'(B)$
\eenum
and \S\ref{pp:affine-top} \ref{top-1}--\ref{top-new-2} hold, granted that we omit the closed embedding aspect of \ref{top-new} and require $B'$ in \ref{top-4}--\ref{top-new-2} to be such that $B'^\times \subset B'$ is open and the inversion map of $B'^\times$ is continuous. If $B$ is local, then \Cref{GR-lemma} holds for all $B$-schemes and in this paragraph we may remove the assumption that the locally of finite type $B$-scheme $X$ be either ind-quasi-affine or subprojective.

\epp

The promised extension of \Cref{thm:basic-Elkik} to nonaffine $X$ rests on \Cref{find-affine}, which is an extension of \Cref{GR-lemma} to a case when $X$ is no longer defined over $B$. The present version of this lemma was suggested by Laurent Moret-Bailly and uses the following material on topologizing modules. 

\bpp[Canonical topology on finite modules] \label{pp:can-topo}
Let $B$ be a topological ring. We endow every finite $B$-module $M$ with the quotient of the product topology with respect to a surjection $\pi \colon B^{\oplus n} \surjects M$ of $B$-modules. By \cite{GR18}*{Lemma 8.3.34, Corollary 8.3.37}, this \emph{canonical topology} of $M$ does not depend on the choice of $\pi$, makes $\pi$ an open map, and makes $M$ a topological $B$-module. By \cite{GR18}*{Proposition 8.3.36}, its formation is compatible with finite products, 
in particular, for projective $M$ it agrees with the subspace topology of $M$ as a direct summand of a finite free $B$-module. By \cite{GR18}*{Proposition 8.3.36, Corollary 8.3.37}, every morphism (respectively, surjection) of finite $B$-modules is continuous (respectively, open) in the canonical topology.
\epp

\begin{bclaim} \label{lem:LMB}
If $B^\times \subset B$ is open and $M$ and $M'$ are finite projective $B$-modules, then the surjective homomorphisms form an open subset of the finite $B$-module $\Hom_B(M, M')$. 
\end{bclaim}

\bpf
There are finite $B$-modules $\wt{M}$ and $\wt{M}'$ such that $M \oplus \wt{M} \simeq B^{\oplus n}$ and $M' \oplus \wt{M}' \simeq B^{\oplus n'}$, and by adding a free summand we arrange that there exist a $B$-module surjection $f \colon \wt{M} \surjects \wt{M}'$. In particular, $\Hom_B(M, M')$ is a direct summand of the free $B$-module $\Hom_B(B^{\oplus n}, B^{\oplus n'})$. By considering its coset given by $f$ and vanishing crossed terms $M \ra \wt{M}'$ and $\wt{M} \ra M'$, we therefore reduce to considering $\Hom_B(B^{\oplus n}, B^{\oplus n'})$. The latter is the $B$-module of $n' \times n$ matrices and, by Nakayama's lemma \cite{SP}*{Lemma \href{https://stacks.math.columbia.edu/tag/00DV}{00DV}}, which allows us to test the surjectivity $B$-fiberwise, the locus of surjective homomorphisms corresponds to those matrices whose $n' \times n'$ minors generate the unit ideal in $B$. Thus, since $B^\times \subset B$ is open, the locus of surjective homomorphisms is also open.  
\epf

\begin{bclaim} \label{claim:dense}
For a ring homomorphism $B_0 \ra B$ with a dense image, a finite $B_0$-module $M_0$, and a topological $B$-module structure on $M \ce B \tensor_{B_0} M_0$ \up{for instance, the one given by the canonical topology}, the image of the map $M_0 \ra M$ is dense. 
\end{bclaim}

\bpf
Since $M$ is a topological $B$-module, every $B$-module surjection $B^{\oplus n} \surjects M$ is continuous. Thus, a choice of a surjection $B_0^{\oplus n} \surjects M_0$ reduces us to the evident case when $M_0 = B_0^{\oplus n}$ and its base change $M \cong B^{\oplus n}$ is endowed with the product topology.
\epf

\blem \label{find-affine}
Let $B$ be a topological ring such that $B^\times \subset B$ is open, let $B_0 \ra B$ be a ring map with a dense image, and let $X$ be a $B_0$-scheme. If either
\benumr
\m \label{FA-i}
$X$ is ind-quasi-affine\uscolon or

\m \label{FA-ii}
$X$ is subprojective and the map $\Pic(B_0) \ra \Pic(B)$ is surjective\uscolon
\eenum
then every $B$-point of $X$ factors through an affine open of $X$.
\elem

\bpf
We begin with the ind-quasi-affine case, in which, since every $B$-point factors through a quasi-compact open, we may assume that $X$ is quasi-affine. Then 
\[
\tst X = \Spec(C) \setminus V(\fc) \qxq{for some $B_0$-algebra} C \qxq{and ideal} \fc \subset C,
\]
so a $B$-point of $X$ is given by a map $C \ra B$ such that the images $c_i \in B$ of some $\wt{c}_i \in \fc$ satisfy 
\[
\tst \sum_i b_i c_i = 1  \qxq{for some} b_i \in B.
\]
Since $B^\times \subset B$ is open and $B_0 \ra B$ has dense image, there are $b_i' \in B_0$ whose images are close to the $b_i$ such that $\sum_i b_i' c_i \in B^\times$. Thus, our $B$-point factors through the affine open $\Spec(C[\f{1}{\sum_i b_i' c_i}]) \subset X$.

In the remaining subprojective case, $X$ is open in a closed subscheme of some $\bP^n_{B_0}$, so  the settled ind-quasi-affine case reduces us to $X = \bP^n_{B_0}$. A $B$-point of $\bP^n_{B_0}$ amounts to a $B$-module quotient 
\[
\tst \pi \colon B^{\oplus (n + 1)} \surjects M \qx{with $M$ projective of rank $1$.}
\]
By the surjectivity assumption for the Picard groups, $M \simeq B \tensor_{B_0} M_0$ for a finite projective $B_0$-module $M_0$ of rank $1$, and we consider the composition
\[
\Hom_{B_0}(M_0, B_0^{\oplus(n + 1)}) \xra{\id_B\tensor_{B_0} -} \Hom_B(M, B^{\oplus (n + 1)}) \xra{\pi \circ - } \Hom_B(M, M).
\]
\Cref{lem:LMB,claim:dense} applied to this composition supply a $B_0$-module map
\[
s \colon M_0 \ra B_0^{\oplus(n + 1)} \qxq{such that the composition} \pi \circ (\id_B \tensor_{B_0} s) \colon M \ra M \qx{is surjective.}
\]
Since $M$ is finitely generated, this composition is then an isomorphism, so that, in particular, the map $\id_B \tensor_{B_0} s \colon M \ra B^{\oplus(n + 1)}$ is an inclusion of a direct summand that is complementary to $\Ker(\pi)$. In particular, the map $\Spec(B) \ra \Spec(B_0)$ factors through the maximal open $U \subset \Spec(B_0)$ over which the dual map $s^\vee \colon B_0^{\oplus(n + 1)} \ra M_0^\vee$ is surjective. Over this open, $\Ker(s^\vee)$ is locally a direct summand of $B_0^{\oplus(n + 1)}$, so, by forming duals again, $s$ is locally an inclusion of a direct summand~over~$U$.


The settled ind-quasi-affine case shows that $\Spec(B)$ factors through an affine open of $U$. Thus, we may assume that $U = \Spec(B_0)$, so that $s\colon M_0 \hra B_0^{\oplus(n + 1)}$ is an inclusion of a direct summand. We obtain a hyperplane $H \subset \bP^n_{B_0}$ that parametrizes those projective rank $1$ quotients of $B_0^{\oplus(n + 1)}$ whose composition with $s$ vanishes. By construction, the quotient given by $\pi$ corresponds to a $B$-point of $\bP^n_{B_0}$ that does not meet $H$, so that this point factors through the affine open $\bP^n_{B_0} \setminus H$, as desired. 
\epf

We now take advantage of the topology of \S\ref{xI-top} to present a nonaffine version of \Cref{thm:basic-Elkik}. This generalizes \cite{GR03}*{Proposition 5.4.21}, which, roughly speaking, is part \ref{BA-i} of the following theorem in the case when $B = A[\f1t]$ for a Henselian Gabber--Ramero triple $(A, t, I)$. 

\bthm \label{thm:beyond-affine}
Let $B$ be a topological ring such that $\wh{B}^\times \subset \wh{B}$ is open and the inversion map of $\wh{B}^\times$ is continuous in the subspace topology, assume that there is an  open nonunital subring $B' \subset B$ that is Zariski and whose induced topology has an open neighborhood base of zero consisting of ideals of $B'$, and let $X$ be a smooth $B$-scheme. 
If either 
\benumr
\m \label{BA-i}
$B'$ is Henselian and $X$ is either ind-quasi-affine or subprojective\uscolon or

\m \label{BA-ii}
$X$ is ind-quasi-affine and $X(C) \surjects X(C/J)$ for every Zariski $B$-pair $(C, J)$ with $C/J \cong \wh{B}$\uscolon or 

\m
$\wh{B}$ is local and $B'$ is Henselian\uscolon or

\m \label{banana}
$\wh{B}$ is local and $X(C) \surjects X(C/J)$ for every Zariski $B$-pair $(C, J)$ with $C/J \cong \wh{B}$\uscolon
\eenum
then the pullback map
\[
X(B) \ra X(\wh{B}) \qxq{has dense image for the topology defined in \uS\uref{xI-top}.}
\]
\ethm

\bpf
By \Cref{Hi-algebraize}~\ref{Hi-a}, in the case \ref{BA-i} the map $\Pic(B) \ra \Pic(\wh{B})$ is an isomorphism. Thus, \Cref{find-affine} shows that every $\wh{B}$-point of $X$ factors through an affine open of $X$, which automatically inherits the surjectivity assumption in \ref{BA-ii} and \ref{banana}. Consequently, the desired density follows from the affine case established in \Cref{thm:basic-Elkik} and the definitions of \S\ref{xI-top}. 
\epf

We now turn to extending this approximation theorem beyond affine $B$ in the setting of Gabber--Ramero triples in \Cref{Elk-Gab}. For this, we will use patching techniques: we review a generalization of the Beauville--Laszlo patching in \Cref{BL-glue}, deduce \Cref{BL-cute-cor}, then review the Ferrand patching in \Cref{more-glue}, and deduce \Cref{mod-rN-cute}.



\blem \label{BL-glue}
Let $A$ be a ring, let $t \in A$, and let $A \ra A'$ be a ring map that induces an isomorphism on derived $t$-adic completions \up{concretely, $A/(t^m) \isomto A'/(t^m)$ for $m > 0$ and $A\langle t^\infty \rangle \isomto A' \langle t^\infty \rangle$}.
\benum
\item \lab{BLG-a}
Base change is an equivalence from the category of $A$-modules $M$ such that $M \ra M \tensor_A A'$ induces an isomorphism on derived $t$-adic completions, concretely, such that
\be \label{glue-cond}
\tst M \langle t^\infty \rangle \hra M \tensor_A A',
\ee
to that of triples consisting of an $A[\f{1}{t}]$-module, an $A'$-module, and an isomorphism of their base changes to $A'[\f{1}{t}]$\uscolon an $A$-flat $M$ satisfies \eqref{glue-cond}, and any $M$ is flat \up{respectively, finite\uscolon respectively, finite projective} if and only if the same holds for its base changes to both $A'$ and~$A[\f{1}{t}]$.

\item \lab{BLG-b}
{\upshape (de Jong).} For a flat, quasi-affine $A$-group scheme $G$, base change is an equivalence from the category of $G$-torsors $\cT$ to that of triples 
\be \lab{BLG-eq}
\qq (T, T', \iota\colon T_{A'[\f{1}{t}]} \isomto T'_{A'[\f{1}{t}]})
\ee
consisting of a $G_{A[\f{1}{t}]}$-torsor $T$, a $G_{A'}$-torsor $T'$, and an indicated torsor isomorphism $\iota$.
\eenum
\elem

\bpf 
\hfill
\benum
\item
Before entering the argument, we recall that the case when $t$ is a nonzerodivisor on both $A$ and $M$ and $A'$ is the $t$-adic completion $\wh{A}$ of $A$ amounts to the main result of \cite{BL95}. The case when $t$ is a nonzerodivisor and $A'$ is arbitrary follows from \cite{BD19}*{Theorem 2.12.1}. The nonzerodivisor assumption was removed by de Jong in \cite{SP}*{Section \href{https://stacks.math.columbia.edu/tag/0BNI}{0BNI}}, whose argument was partly inspired by that of Kedlaya--Liu carried out in \cite{KL15}*{Section 2.7}. The proofs of \cite{SP}*{Section \href{https://stacks.math.columbia.edu/tag/0BNI}{0BNI}} turned out to work beyond the case $A' = \wh{A}$, and we have updated \cite{SP}*{Section \href{https://stacks.math.columbia.edu/tag/0BNI}{0BNI}} to accommodate for this. 

In more detail, by \cite{SP}*{Lemma \href{https://stacks.math.columbia.edu/tag/0BNR}{0BNR}}, the pair $(A \ra A', t)$ is ``glueing'' and, by \cite{SP}*{Lemma~\href{https://stacks.math.columbia.edu/tag/0BNW}{0BNW}}, the condition \eqref{glue-cond} is equivalent to $M \ra M \tensor_A A'$ inducing an isomorphism on derived $t$-adic completions and amounts to $M$ being ``glueable for $(A \ra A', t)$''; by \cite{SP}*{Remark \href{https://stacks.math.columbia.edu/tag/0BNX}{0BNX}}, any $A$-flat $M$ satisfies \eqref{glue-cond}. Thus, \cite{SP}*{Theorem \href{https://stacks.math.columbia.edu/tag/0BP2}{0BP2}} gives the claimed equivalence of categories. The assertion about testing properties over $A'$ and $A[\f{1}{t}]$ follows from \cite{SP}*{Lemmas \href{https://stacks.math.columbia.edu/tag/0BP7}{0BP7}, \href{https://stacks.math.columbia.edu/tag/0BNN}{0BNN}, and \href{https://stacks.math.columbia.edu/tag/0BP6}{0BP6}}. 

\item
The full faithfulness follows from \cite{SP}*{Lemma \href{https://stacks.math.columbia.edu/tag/0F9T}{0F9T}}, according to which a similar base change functor is fully faithful even on the category of flat algebraic spaces with affine diagonal. For the essential surjectivity, since $T$ and $T'$ are quasi-affine, by \cite{SP}*{Lemmas \href{https://stacks.math.columbia.edu/tag/0F9U}{0F9U} and \href{https://stacks.math.columbia.edu/tag/0F9R}{0F9R}}, any triple as in \eqref{BLG-eq} arises from a faithfully flat, quasi-compact, separated $A$-algebraic space $\cT$. By the full faithfulness \cite{SP}*{Lemma \href{https://stacks.math.columbia.edu/tag/0F9T}{0F9T}} again, $\cT$ comes equipped with a $G$-action for which the map 
\[
\qq G \times_A \cT \xra{(g,\, \tau)\, \mapsto\, (g\tau,\, \tau)} \cT \times_A \cT
\]
is an isomorphism. Consequently, $\cT$ is a $G$-torsor (and hence is a quasi-affine scheme).
\qedhere
\eenum
\epf

\bprop \label{BL-cute-cor}
For a ring $A$, a $t \in A$, a ring map $A \ra A'$ that induces an isomorphism on derived $t$-adic completions,
a flat $A$-scheme $U$, and a $U$-scheme $X$ that is either $U$-ind-quasi-affine 
or $U$-subprojective,\footnote{Added in proof: Our assumption on $X$ is not optimal, in fact, \eqref{BL-cute-eq} holds for any $U$-scheme $X$ and also for any quasi-compact, quasi-separated $U$-algebraic space $X$. To argue this, first reduce to the case $U = \Spec(A)$ as in the proof. Then the case of a quasi-compact, quasi-separated algebraic space follows from \cite{Bha16}*{Theorem 1.4 1.} (to check its assumption, use \cite{flat-purity}*{Lemma 5.4.3}). In the remaining scheme case, by considering all the quasi-compact opens of $X$ we reduce to the case when $X$ is quasi-compact. As in the proof of \cite{torsors-regular}*{Proposition 6.1.1 (c)}, at least as far as ring-valued points are concerned, \cite{SP}*{Lemma \href{https://stacks.math.columbia.edu/tag/03K0}{03K0}} implies that a quasi-compact scheme is a filtered direct limit of quasi-compact, quasi-separated schemes, so to reduce to the already settled case of the latter it remains to use the fact that filtered direct limits of sets commute with fiber products.} 
we have
\be \label{BL-cute-eq}
\tst X(U) \isomto X(U_{A'}) \times_{X(U_{A'[\f{1}{t}]})} X(U_{A[\f{1}{t}]}).
\ee
\eprop

\bpf
The assumptions are stable upon replacing $A$ (respectively, $A'$) by the coordinate ring of a variable affine open of $U$ (respectively, of its base change to $A'$). Thus, since the functors that underlie both sides of \eqref{BL-cute-eq} are Zariski sheaves on $U$, by passing to such an open we may assume that $U = \Spec(A)$. In this case, which we now assume, \cite{SP}*{Lemma \href{https://stacks.math.columbia.edu/tag/0BNR}{0BNR}} shows that the assumption on $A \ra A'$ implies that
\[
\tst A \isomto A' \times_{A'[\f{1}{t}]} A[\f{1}{t}].
\]
This key identification already gives \eqref{BL-cute-eq} for affine $X$. It also implies that the surjection 
\be \label{sch-dom}
\tst \Spec(A') \bigsqcup \Spec(A[\f{1}{t}]) \surjects \Spec(A)
\ee
is schematically dominant, so an $A$-point of an $A$-scheme $Y$ factors through a given open (respectively, closed) subscheme if and only if the same holds for its pullbacks to $A'$ and $A[\f{1}{t}]$. In particular, \eqref{BL-cute-eq} for quasi-affine $X$ follows from its case for affine $X$. Thus, for ind-quasi-affine $X$, the identification \eqref{BL-cute-eq} holds for every quasi-compact open $X' \subset X$ in place of $X$ and, since every finite collection of $A$-points (respectively, $A'$- and $A[\f{1}{t}]$-points) of $X$ factors through such an $X'$, also for $X$ itself.

We turn to the remaining case when $X$ is open in a projective $A$-scheme. The same  reduction allows us to assume that $X$ is projective, then that $X = \bP(\sE)$ for a quasi-coherent, finite type module $\sE$ on $\Spec(A)$, and finally, by choosing a surjection $\sO^{\oplus (n + 1)} \surjects \sE$, that $X = \bP^n$. By \cite{EGAI}*{Corollaire~9.5.6}, the schematic dominance of \eqref{sch-dom} and the separatedness of $X$ ensure that \eqref{BL-cute-eq} is injective. For the remaining surjectivity, by the description of the functor of points of $\bP^n$ given in \cite{EGAII}*{Th\'{e}or\`{e}me~4.2.4}, we need to show that any pair of compatible surjections 
\[
\tst A'^{\,\oplus (n + 1)} \surjects M' \qxq{and} (A[\f{1}{t}])^{\oplus (n + 1)} \surjects M''
\]
with $M'$ (respectively, $M''$) a finite projective $A'$-module (respectively, $A[\f{1}{t}]$-module) of rank $1$ is a base change of such a surjection of $A$-modules. \Cref{BL-glue}~\ref{BLG-a} supplies the unique candidate $\pi \colon A^{\oplus(n + 1)} \ra M$ and implies that $M$ is finite projective of rank $1$. It remains to observe that $\pi$ is surjective, as may be checked after base change to the residue fields of $A$.
\epf

\blem[\cite{Fer03}*{Th\'{e}or\`{e}me 2.2 iv)}] \label{more-glue}
For a fiber product $R_1 \times_R R_2$ of rings with either $R_1 \ra R$ or $R_2 \ra R$ surjective, pullback is an equivalence from the category of flat $(R_1 \times_R R_2)$-modules to that of triples consisting of a flat $R_1$-module, a flat $R_2$-module, and an isomorphism of their base changes to $R$\uscolon the same holds with `flat' replaced by `finite projective.' Moreover, an $(R_1 \times_R R_2)$-module is flat \up{respectively, finite\uscolon respectively, finite projective} if and only if so are its base changes to $R_1$ and $R_2$.\footnote{Added in proof: The lemma implies that for any flat, finitely presented, affine $(R_1 \times_R R_2)$-group scheme $G$, pullback is an equivalence from the category of $G$-torsors to that of triples consisting of a $G_{R_1}$-torsor, a $G_{R_2}$-torsor, and an isomorphism of their base changes to $R$, see \cite{Sta19}*{Lemma 3.1 and its proof}; the same holds for any flat, affine $(R_1 \times_R R_2)$-group scheme $G$ granted that we consider torsors in the fpqc topology.} \QED
\elem

\beg \label{more-glue-eg}
The case when both $R_1 \ra R$ and $R_2 \ra R$ are surjective corresponds to a ring $A$ and ideals $I_1, I_2 \subset A$ with $I_1 \cap I_2 = 0$: then
\be \label{fibre-prod-0}
A \isomto A/I_1 \times_{A/(I_1 + I_2)} A/I_2.
\ee
For instance, for a Gabber--Ramero triple $(A, t, I)$ with $I\langle t^\infty\rangle = I\langle t^n\rangle$ for an $n > 0$ and $\wt{A} \ce A/I\langle t^\infty \rangle$,
\be \label{fibre-prod}
A \isomto \wt{A} \times_{\wt{A}/(t^nI)} A/t^nI.
\ee
\eeg

\bprop \label{mod-rN-cute}
For a ring $A$ and ideals $I_1, I_2 \subset A$ with $I_1 \cap I_2 = 0$, 
a flat $A$-scheme $U$, and a $U$-scheme $X$ that is either $U$-ind-quasi-affine or $U$-subprojective,\footnote{Added in proof: Our assumptions are not optimal, in fact, \eqref{mod-rN-eq} holds for any $U$-algebraic space $X$ and even with $A$ replaced by any fiber product $R_1 \times_R R_2$ as in \Cref{more-glue} (compare with \Cref{more-glue-eg}). To argue this, first reduce to the case $U = \Spec(A)$ as in the proof. Then use the result of Temkin--Tyomkin \cite{TT16}*{Lemma 4.1 and Theorem 4.3}, according to which $\Spec(-)$ transforms fiber products as in \Cref{more-glue} into pushouts in the category of algebraic spaces. }
\be \label{mod-rN-eq}
X(U) \isomto X(U_{A/I_1}) \times_{X(U_{A/(I_1 + I_2)})} X(U_{A/I_2}). 
\ee
\eprop

\bpf
As in the proof of \Cref{BL-cute-cor}, both sides of \eqref{mod-rN-eq} are Zariski sheaves in $U$, so we may work locally on $U$ to reduce to the case when $U = \Spec(A)$. Then the assumption $I_1 \cap I_2 = 0$ implies the schematic dominance and the surjectivity of the map
\[
\tst \Spec(A/I_1) \bigsqcup \Spec(A/I_2) \ra \Spec(A).
\]
Thus, the arguments of the proof of \Cref{BL-cute-cor} reduce to $X = \bP^n$ and also prove the injectivity of \eqref{mod-rN-eq} in this case. For the surjectivity, we first use \Cref{more-glue} to glue compatible surjections 
\[
\tst (A/I_1)^{\, \oplus (n + 1)} \surjects M' \qxq{and} (A/I_2)^{\oplus (n + 1)} \surjects M''
\]
onto finite projective modules of rank $1$ to a unique $A$-module map $\pi \colon A^{\,\oplus (n + 1)} \ra M$ with $M$ finite projective of rank $1$, and then check on the residue fields of $A$ that $\pi$ is surjective.
\epf

\bpp[$\x{The $(t, I)$-adic topology on $A[\f 1t]$-points}$] \label{pp:nonaffine-top}
Let $(A, t, I)$ be a Zariski Gabber--Ramero triple and let $B$ be a topological ring that is either $A$ or $A[\f 1t]$ (see \S\ref{GR-triples}). Then $B^\times \subset B$ is open and the inversion map of $B^\times$ is continuous for the subspace topology: indeed, for instance, multiplication by any $a \in A[\f 1t]^\times$ is a homeomorphism of $A[\f 1t]$, so, by the Zariski assumption, the $a + a t^n\ov{I} \subset A[\f1t]^\times$ for varying $n > 0$ form an open neighborhood base of $a$. Thus,  \S\ref{extend-top} applies and endows $X(B)$ with a topology for every $B$-scheme $X$ that is either ind-quasi-affine or subprojective. If $(A, t, I)$ is even Henselian, then, although we will not use this, smooth morphisms of quasi-projective $A[\f 1t]$-schemes induce open maps on $A[\f 1t]$-points, see \cite{GR03}*{Proposition 5.4.29 (ii)} for a precise statement.

We use the resulting \emph{$(t, I)$-adic topologies} on the sets of sections $X(A)$ and $X(A[\f 1t])$ to topologize sets of sections beyond the case of an affine base as follows. Fix an open subscheme
\[
\tst \Spec(A[\f{1}{t}]) \subset U \subset \Spec(A)
\]
and a locally of finite type $U$-scheme $X$ that is either $U$-ind-quasi-affine or $U$-subprojective. For every affine open $V \subset U$, consider the Zariskization $A_V$ of the coordinate ring of $V$ along the closed subscheme $V_{A/tI}$ as well as the Zariski Gabber--Ramero triple $(A_V, t, I A_V)$ (which vanishes if $V_{A/tI} = \emptyset$). We endow $X(U)$ with the coarsest topology for which the pullback maps
\[
\tst X(U) \ra X(A[\f 1t]) \qxq{and} X(U) \ra X(A_V) \qxq{for affine opens} V \subset U
\]
are all continuous. By \S\ref{extend-top} and \S\ref{pp:affine-top}~\ref{top-4}, when $U = \Spec(A)$ or $U = \Spec(A[\f 1t])$, this definition agrees with the one given in the previous paragraph. Similarly, \S\ref{extend-top} and \S\ref{pp:affine-top}~\ref{top-1} and \ref{top-4} ensure that a $U$-morphisms $X \ra X'$ induces a continuous map $X(U) \ra X'(U)$ and that the restriction map $X(U) \ra X(U')$ is continuous for every open $\Spec(A[\f 1t]) \subset U' \subset U$. 

It is useful to know that the topology may be constructed using any affine open cover of $U$ as follows.
\epp

\begin{bclaim}\label{other-desc}
For any affine open cover $U = \Spec(A[\f 1t])\cup \bigcup_{i \in I} V_i$, the topology on $X(U)$ defined above is the coarsest one for which the maps $X(U) \ra X(A[\f 1t])$ and $X(U) \ra X(A_{V_i})$ are all continuous. In particular, for a morphism of Zariski Gabber--Ramero triples $(A, t, I) \ra (A', t, I')$ and an open subscheme $\Spec(A'[\f1t]) \subset U' \subset U_{A'}$, the map $X(U) \ra X(U')$ is continuous. 
\end{bclaim}

\bpf
For the first assertion, we need to show that with this a priori coarser topology on $X(U)$ and every affine open $V \subset U$, the pullback $X(U) \ra X(A_V)$ is already continuous. If $V \subset V_i$, then this follows from the continuity of $X(A_{V_i}) \ra X(A_{V})$. If $V \subset \Spec(A[\f 1t])$, then the topology of $A_V$ is the pullback of the discrete topology of $A_V/IA_V$, so $A[\f 1t] \ra A_V$ is continuous, and we analogously use the resulting continuity of $X(A[\f 1t]) \ra X(A_V)$. In general, we cover $V$ by affine opens $V_j$ each one of which lies either in some $V_i$ or in $\Spec(A[\f 1t])$ and we therefore reduce to showing that the topology on $X(A_V)$ is the coarsest one for which the maps $X(A_V) \ra X(A_{V_j})$ are all continuous. 

For showing that every open of $X(A_V)$ is the pullback of an open of $\prod_j X(A_{V_j})$, we may first assume that $X$ is an $A_V$-scheme and then, by \Cref{GR-lemma}, that it is also affine, so that it is a closed subscheme of some $\bA^N$. The subset $\Spec(A_V) \subset V$ consists of the generizations of points in $V_{A/tI}$, and likewise for each $V_j$. Thus, since the $V_j$ cover $V$, the map $A_V \ra \prod_j A_{V_j}$ is injective, so that
\[
\tst X(A_V) = \bA^N(A_V) \cap \prod_j X(A_{V_j}) \qxq{inside} \prod_j \bA^N(A_{V_j}). 
\]
This reduces us to the case $X = \bA^N$, and then further to $X = \bA^1$. In this case, since 
\[
\tst \Spec(A_V/(t^nI)) = \bigcup_{j} \Spec(A_{V_j}/(t^nI)) \qxq{for every} n > 0,
\]
for every $a \in A_V$ the open $a + t^n IA_V \subset A_V$ is the intersection of the preimages of the opens $a + t^n IA_{V_j} \subset A_{V_j}$. The cosets $a + t^n IA_V$ form a base of the topology of $A_V$, so the assertion follows. 

To deduce the continuity of $X(U) \ra X(U')$, we first use the definition of the topology on $X(U')$ to reduce to $U' = \Spec(A')$. We then apply the first part of the claim to $A'$ to reduce to the case when $U' \ra U$ factors through an affine open $V \subset U$, when the continuity of $X(A_V) \ra X(A')$ suffices. 
\epf

We use this more practical description of the topology of $X(U)$ to exhibit the following openness.

\begin{bclaim} \label{claim:local}
If $I\langle t^\infty\rangle = 0$ and $U$ is quasi-compact, then the pullback map 
\[
\tst X(U) \ra X(A[\f 1t]) \qx{is open and a local homeomorphism.}
\]
In particular, if $A\langle t^\infty\rangle = 0$ and $U$ is quasi-compact, then $X(U) \subset X(A[\f 1t])$ is an open embedding. 
\end{bclaim}

\bpf
If $A\langle t^\infty\rangle = 0$, then, by the separatedness of $X$ and \cite{EGAI}*{Corollaire 9.5.6}, the map $X(U) \ra X(A[\f 1t])$ is injective, so the final assertion follows from the rest. By \S\ref{GR-triples}, the assumption $I\langle t^\infty\rangle = 0$ implies that the restriction of the map $A \ra A[\f 1t]$ to the open ideal $I \subset A$ is an open topological embedding, and analogously for $A_{V} \ra A_{V}[\f1t]$ for every affine open $V \subset U$. Thus, \S\ref{extend-top} and \S\ref{pp:affine-top}~\ref{top-new-2} imply that the pullback maps $X(A_{V}) \ra X(A_{V}[\f 1t])$ are open and local homeomorphisms. 

Quasi-compactness of $U$ supplies a finite affine open cover $U = \Spec(A[\f 1t]) \cup \bigcup_{i \in I} V_i$. To exhibit an open neighborhood $W$ of a fixed $x \in X(U)$ such that $X(U) \ra X(A[\f 1t])$ maps $W$ homeomorphically onto an open subset $W_t$, we choose open neighborhoods $W_i \subset X(A_{V_i})$ and $W_{ii'} \subset X(A_{V_{i} \cap V_{i'}})$ of the pullbacks of $x$ such that the restrictions to $W_i$ and $W_{ii'}$ of the pullbacks discussed in the previous paragraph are homeomorphisms onto their open images $W_{i,\,t}$ and $W_{ii',\, t}$. We then let $W \subset X(U)$ (respectively, $W_t \subset X(A[\f 1t])$) be the intersection of the preimages of the $W_i$ and $W_{ii'}$ (respectively, $W_{i,\, t}$ and $W_{ii',\, t}$). It remains to show that the continuous map $W \ra W_t$ induced by pullback is a homeomorphism. 

By fpqc descent, giving a $w \in W$ amounts to giving its image $w_t \in W_t$ together with elements $w_i \in W_i$ that agree with $w_t$ and are compatible under pullbacks to the $X(A_{V_i \cap V_{i'}})$. However, by construction of $W_t$, the element $w_t$ alone gives rise to unique such compatible $w_i$. 
Therefore, the map $W \ra W_t$ is bijective. To conclude that it maps every open $W' \subset W$ to an open of $W_t$, it remains to first recall from \Cref{other-desc} that locally $W'$ is an intersection of an open of $W_t$ and of preimages of opens of $W_i$ and to then note that, by construction, this intersection may equivalently be formed in $W_t$. 
\epf

For the mere openness of the map in \Cref{claim:local}, mere boundedness of $(A, t, I)$ suffices as follows.

\begin{bclaim} \label{claim:const}
If $U$ is quasi-compact, then the pullback maps
\[
X(U) \ra X(U_{A/t^nI}) \qx{have open fibers.}
\]
If, in addition, the Zariski Gabber--Ramero triple $(A, t, I)$ is bounded and $\wt{A} \ce A/I\langle t^\infty\rangle$, then both
\[
\tst X(U) \ra X(U_{\wt{A}}) \qxq{and} X(U) \ra X(A[\f 1t]) \qx{are open maps.}
\]
\end{bclaim}

\bpf
Since $U_{A/t^nI}$ is covered by its intersections with finitely many affine opens $V \subset U$, the claimed openness of the fibers follows from its case $U = \Spec(A)$ supplied by \S\ref{extend-top} and \S\ref{pp:affine-top}~\ref{top-4}. For the rest, due to \Cref{claim:local}, it suffices to establish the assertion about $X(U) \ra X(U_{\wt{A}})$. However, by \Cref{more-glue-eg} and \Cref{mod-rN-cute}, the boundedness assumption implies that for every $n > 0$ such that $I\langle t^\infty \rangle = I \langle t^n\rangle$, we have
\be \label{X-fibre-prod}
X(U) \isomto X(U_{\wt{A}}) \times_{X(U_{\wt{A}/(t^nI)})} X(U_{A/t^nI}).
\ee
In particular, for every such $n$, each fiber of the map $X(U) \ra X(U_{A/t^nI})$ maps bijectively to a fiber of the map $X(U_{\wt{A}}) \ra X(U_{\wt{A}/(t^nI)})$. By the first part of the claim, these fibers are open and, by \Cref{other-desc}, the bijections in question are continuous. For arguing that they are also open, we choose a finite affine open cover $U = \Spec(A[\f 1t]) \cup \bigcup_{i \in I} V_i$ and combine \Cref{other-desc} with the analogues of \eqref{X-fibre-prod} for the $A_{V_i}$ to reduce to the case $U = \Spec(A)$. We may then assume that $X$ is affine, so that it is a closed subscheme of some $\bA^N$. By again working with the fibers of the reduction modulo $t^nI$, we reduce further to $X = \bA^N$, so even to $X = \bA^1$. To then conclude it remains to note that the cosets $a + t^nI$ for varying $a \in A$ and $n > 0$ form an open base for the topology of $A$ and that these cosets map onto their counterparts for $\wt{A}$. 
\epf


We are ready for the following non-Noetherian version of the Elkik approximation theorem. 

\bthm \label{Elk-Gab}\label{Elk-Gab-Zar}
Let $(A, t, I)$ be a Zariski, bounded Gabber--Ramero triple, let 
\[
\tst \Spec(A[\f{1}{t}]) \subset U \subset \Spec(A)
\]
be a quasi-compact open, and let $X$ be a locally of finite type $U$-scheme that is either $U$-ind-quasi-affine or $U$-subprojective and such that $X_{A[\f{1}{t}]}$ is $A[\f{1}{t}]$-smooth. Under one of the following assumptions\ucolon
\benumr
\m
$(A, t, I)$ is Henselian\uscolon or

\m
$X_{A[\f1t]}$ is ind-quasi-affine and $X(C) \surjects X(C/J)$ for Zariski $A[\f{1}{t}]$-pairs $(C, J)$ with $C/J \cong~\!\wh{A}[\f1t]$\uscolon\!\
\eenum
the following pullback map has a dense image for the topology defined in \uS\uref{pp:nonaffine-top}\ucolon
\[
X(U) \ra X(U_{\wh{A}});
\]
in particular, for every $\wh{x} \in X(U_{\wh{A}})$ and $n > 0$, 
there is an 
\[
\tst \qq x \in X(U) \q  \x{whose pullback to} \q X(U_{A/t^nI}) \overset{\eqref{R-Rhat-comp}}{\cong} X(U_{\wh{A}/t^n\wh{I}}) \q \x{agrees with that of $\wh{x}$.}
\]
\ethm

\bpf
The last assertion follows from the rest and \Cref{claim:const}. Moreover, in the key case when $U = \Spec(A[\f{1}{t}])$ the desired density follows from \Cref{thm:beyond-affine}. To deduce the general case, it then remains to recall from \Cref{claim:const} (with \Cref{nzd-example}) that the map $X(U_{\wh{A}}) \ra X(\wh{A}[\f1t])$ is open and to use the following identification supplied by \Cref{BL-cute-cor} (with \eqref{R-Rhat-comp}):
\[
\tst X(U) \isomto X(U_{\wh{A}}) \times_{X(\wh{A}[\f{1}{t}])} X(A[\f{1}{t}]). \qedhere
\]
\epf

\brem \label{etale-eg}
In \Cref{Elk-Gab}, if $(A, t, I)$ is Henselian and $X_{A[\f{1}{t}]}$ \'{e}tale over $A[\f{1}{t}]$, then
\[
X(U) \isomto X(U_{\wh{A}}).
\]
Indeed, the end of the proof reduces us to $U = \Spec(A[\f 1t])$, when the target is discrete because, by \'{e}taleness, an $\wh{A}[\f{1}{t}]$-point of $X_{\wh{A}[\f{1}{t}]}$ is an inclusion of a clopen that maps isomorphically to $\Spec(\wh{A}[\f 1t])$. 
\erem

\beg \label{eg:cute-affine-Gras}
For a ring $R$ and a smooth, ind-quasi-affine $R\{ t \}$-group scheme $G$, we have 
\[
\tst G(R\{t\}[\f 1t])/G(R\{t\}) \isomto G(R\llp t \rrp)/G(R\llb t \rrb).
\]
Indeed, the map is injective by \Cref{BL-cute-cor} and, since $G(R\llp t \rrp)$ is a topological group and $G(R\llb t \rrb) \subset G(R\llp t \rrp)$ is an open subgroup (see \S\ref{pp:affine-top} \ref{top-7} and \ref{top-new}), it is surjective by \Cref{Elk-Gab}.
\eeg

We are ready to deduce the promised Elkik approximation theorem in the Noetherian setting.

\bthm[compare with \cite{Elk73}*{Th\'{e}or\`{e}me 2 bis}] \label{thm:Elkik-reproof}
For a Henselian pair $(A, J)$ with $A$ Noetherian, an open $\Spec(A) \setminus V(J) \subset U \subset \Spec(A)$, and a subprojective $U$-scheme $X$ such that $X_{\Spec(A) \setminus V(J)}$ is smooth, letting $\wh{A}$ denote the $J$-adic completion, we have that
\[
\xq{for every} n > 0 \qxq{and} \wh{x} \in X(U_{\wh{A}}) \qxq{there is an}  x \in X(U) \qxq{with} x = \wh{x} \qxq{in} X(U_{A/J^n}). 
\]
\ethm

\bpf
By \cite{Mat89}*{Theorem 8.12}, the completion of a Noetherian ring $B$ with respect to an ideal $(b_1, \dotsc, b_g)$ is isomorphic to $B\llb t_1, \dotsc, t_g\rrb/(t_1 - b_1, \dotsc, t_g - b_g)$. 
Thus, fixing generators $a_1, \dotsc, a_g$ of our original $J$, we see that the $J$-adic completion of $A$ agrees with the iterated $a_i$-adic completion and 
that the ideal generated by $J$ is still Henselian (respectively, Zariski) in these resulting intermediate completions (see \S\ref{Zar-Hens}). In conclusion, at the expense of applying the statement $g$ times, we have reduced to the case when $J$ is principal, which is a special case of \Cref{Elk-Gab}.
\epf

\brem
In the proof of \Cref{thm:Elkik-reproof}, the only role of Noetherianity is to ensure that $U$ is quasi-compact, that $J$ is finitely generated, that the $J$-adic completion of $A$ agrees with the iterated completion with respect to a generating set $a_1, \dotsc, a_g$ of $J$, and that for $0 \le i < g $ the iterated completion with respect to $a_1, \dotsc, a_i$ is $a_{i + 1}$-Henselian with bounded $a_{i + 1}^\infty$-torsion  (to apply \Cref{Elk-Gab}). Thus, granted that one imposes these assumptions, which 
hold if $U$ is quasi-compact, $A$ is $J$-Henselian, and the sequence $a_1, \dotsc, a_g$ is, for instance, $A$-regular, one obtains a non-Noetherian generalization of \Cref{thm:Elkik-reproof}. It is possible that a weaker condition on $J$ could suffice for this---for some guiding examples in this direction, see \cite{Nak18}.
\erem




\csub[Algebraization beyond the affine case] \label{section:beyond-affine}

As a final goal of \S\ref{chapter-invariance}, we take advantage of our work in \S\ref{section-approximation} to exhibit a nonaffine version of invariance under Henselian pairs in \Cref{G-tor-inv}, 
a nonaffine version of Gabber's affine analogue of proper base change in \Cref{nonaffine-Gab}, and concrete consequences for algebraization in the Noetherian case in \Cref{Br-comp}. To illustrate the method, we begin with a nonaffine variant of \Cref{cor:stays-connected}.

\bthm \label{thm:idem-case}
For  a map $(A, t, I) \ra (A', t, I')$ of Henselian, bounded Gabber--Ramero triples such that 
$A/t^nI \isomto A'/t^nI'$ for all $n > 0$
and an open 
\[
\tst \Spec(A[\f{1}{t}]) \subset U \subset \Spec(A),
\]
the map $U_{A'} \ra U$ induces a bijection on idempotents, in other words, every clopen subscheme of $U_{A'}$ is the base change of a unique clopen subscheme of $U$. 
\ethm

\bpf
By \Cref{nzd-example}, the topological rings $A[\f{1}{t}]$ and $A'[\f{1}{t}]$ have the same completion, 
so \Cref{cor:stays-connected} settles the case $U = \Spec(A[\f{1}{t}])$. Since $A\langle t^\infty \rangle \isomto A'\langle t^\infty \rangle$ and giving an idempotent amounts to giving a map to $\Spec(\bZ[e]/(e^2 - e))$, the general case then follows from \eqref{BL-cute-eq}.
\epf

The following simple lemma will be useful in the proof of \Cref{G-tor-inv}.

\blem \lab{lem-sup-ss}
For a scheme $X$, a closed $Z \subset X$, an abelian fppf sheaf $\sF$ on $X$, and the \'{e}tale sheafification $\cH^j_Z(-, \sF)$ of the functor $X' \mapsto H^j_Z(X', \sF)$, there is a functorial in $X$ and $\sF$ spectral sequence
\be \lab{eq-sup-ss}
E_2^{ij} = H^i_\et(X, \cH^j_Z(-, \sF)) \Ra H^{i + j}_Z(X, \sF). 
\ee

\elem

\bpf
The spectral sequence will be the one associated to the composition of functors 
\[
H^0_\et(X, \cH^0_Z(-, *)) \cong H^0_Z(X, *).
\]
Indeed, the $\cH^j_Z$ are the derived functors of $\sF \mapsto \cH^0_Z(-, \sF)$: they form a $\delta$-functor that, by \cite{SGA4II}*{Expos\'{e} V, Proposition 4.7, Section 4.6}, kills injectives when $j > 0$, so the general criterion \cite{Har77}*{Chapter III, Theorem 1.3A, Corollary 1.4} for being the left derived functors applies. It remains to note that, by \cite{SGA4II}*{Expos\'{e} V, Proposition 4.11 2)}, if $\sI$ is injective, then $\cH^0_Z(-, \sI)$ computed in the fppf site is also injective sheaf and hence, by preservation of flasque sheaves under pushforwards \cite{SGA4II}*{Expos\'{e} V, Proposition 4.9 1), D\'{e}finition 4.1}, is acyclic on the \'{e}tale site.
\epf

The Brauer group aspect of the following result includes the statement announced in \cite{Gab04}*{Theorem~2}. 

\bthm \label{G-tor-inv}
Let $(A, t, I) \ra (A', t, I')$ be a map of bounded Gabber--Ramero triples such that 
$A/t^nI \isomto A'/t^nI'$ for all $n > 0$, 
let  
\[
\tst \Spec(A[\f{1}{t}]) \subset U \subset \Spec(A)
\]
be a quasi-compact open, and let $G$ be a locally of finite type, $U$-ind-quasi-affine, flat $U$-group such that $G_{A[\f{1}{t}]}$ is smooth over $A[\f{1}{t}]$. 
\benum
\m \label{GI-ii}
If $(A, t, I)$ is Zariski 
and $G_{A[\f{1}{t}]}$ is an fpqc inner form of $\GL_{n,\, A[\f{1}{t}]}$, then
\[
\qq H^1(U, G) \hra H^1(U_{A'}, G).
\]

\m \label{GI-i}
If $(A, t, I)$ is Henselian, 
then 
\[
\qq H^1(U, G) \hra H^1(U_{A'}, G).
\]

\item \lab{GTI-a}
If $(A, t, I)$ and $(A', t, I')$ are Henselian and $G$ is $U$-quasi-affine, 
then
\be \label{eqn:H1-inv-nonaffine}
\qq H^1(U, G) \isomto H^1(U_{A'}, G),
\ee
in particular, then $\Br(U) \isomto \Br(U_{A'})$ \up{equivalently, $H^2(U, \bG_m)_\tors \isomto H^2(U_{A'}, \bG_m)_\tors$}\uscolon

\item \lab{GTI-b}
if $(A, t, I)$ is Henselian and $G$ is commutative, $U$-quasi-affine, of finite presentation, 
then
\[
\qq H^2(U, G) \hra H^2(U_{A'}, G)
\]
and the same also holds with a commutative, finite, locally free $U$-group $G'$ in place of $G$.
\eenum

\ethm


\bpf
The assumption implies that the map $A \ra A'$ induces an isomorphism $\wh{A} \isomto \wh{A'}$, so in all parts we lose no generality by assuming that $(A', t, I') = (\wh{A}, t, \wh{I})$ (for which the assumptions are retained by \eqref{R-Rhat-comp}). Moreover, by \cite{SP}*{Lemma \href{https://stacks.math.columbia.edu/tag/0APK}{0APK}}, the functor that parametrizes isomorphisms between two $G$-torsors is representable by an ind-quasi-affine $U$-scheme $X$ that is locally finite type and whose base change to $A[\f{1}{t}]$ is smooth. Thus, \ref{GI-i} follows from \Cref{Elk-Gab} and, to conclude that so does \ref{GI-ii}, we only need to check that in \ref{GI-ii} the pullback map $X(C) \ra X(C/J)$ is surjective for every Zariski $A[\f1t]$-pair $(C, J)$. However, by \cite{Gir71}*{Chapitre III, Proposition 2.6.1~(i)}, up to translation inner forms have identical collections of torsors, so \Cref{sta-free-eg} implies that if $X(C/J) \neq \emptyset$, then also $X(C) \neq \emptyset$. If this happens, then, after choosing an element of $X(C)$, the map $X(C) \surjects X(C/J)$ becomes identified with $G(C) \ra G(C/J)$, and hence is surjective by \cref{foot:surj}.

The injectivity in \ref{GTI-a} is a special case of \ref{GI-i} and the bijectivity when $U = \Spec(A[\f{1}{t}])$ follows from \Cref{Hi-algebraize}~\ref{Hi-a}. The surjectivity in general then follows from the patching result recorded in \Cref{BL-glue}~\ref{BLG-b} (whose assumptions are met by \Cref{nzd-example}). The Brauer group assertion in \ref{GTI-a} follows from \eqref{eqn:H1-inv-nonaffine} applied to $G = \GL_N$ and $G = \PGL_N$ with varying $N$: indeed, by definition, 
\[
\tst \Br(U) \ce \bigcup_{N \ge 1} \im\p{H^1(U, \PGL_N) \ra H^2(U, \bG_m)_\tors}\!
 \]
 and likewise for $U_{A'}$. The parenthetical assertion in \ref{GTI-a} then follows from Gabber's theorem established in \cite{dJ02} that, in particular, identifies $\Br(-)$ with $H^2_\et(-, \bG_m)_\tors$ for quasi-affine schemes. 
 
In \ref{GTI-b}, the finite locally free case follows from the rest applied to the terms of the B\'{e}gueri resolution 
 \[
 0 \ra G' \ra \Res_{G'^*/U}(\bG_m) \ra Q \ra 0,
 \]
 where $G'^*$ denotes the Cartier dual. Moreover, the case $U = \Spec(A[\f{1}{t}])$ follows from \Cref{Hi-algebraize}~\ref{Hi-b}. For general $U$, we use the cohomology with supports along $\{t = 0\}$ sequences, \ref{GTI-a}, the $U = \Spec(A[\f{1}{t}])$ case of \ref{GTI-b}, and the five lemma, to reduce to showing that\footnote{In this step, in order to apply \eqref{eqn:H1-inv-nonaffine} to $G_{A[\f 1t]}$ we use the $A[\f 1t]$-smoothness of $G_{A[\f 1t]}$, which is a not an entirely natural assumption: the case $U = \Spec(A[\f 1t])$ settled in \Cref{Hi-algebraize}~\ref{Hi-b} did not need it. }
\[
H^2_{\{t = 0\}}(U, G) \hra H^2_{\{t = 0\}}(U_{A'}, G).
\]
As in \Cref{lem-sup-ss}, we let $\cH^*_{\{t = 0\}}$ denote the \'{e}tale sheafification of the flat cohomology with supports in $\{t = 0\}$. Since these \'{e}tale sheaves are supported on the closed subscheme cut out by $t$, we then use the spectral sequences \eqref{eq-sup-ss} to reduce to showing that the sheaves
\[
 \cH^0_{\{t = 0\}}(-, G) \q \x{and} \q \cH^1_{\{t = 0\}}(-, G) \q \x{on $(U_{A'})_\et$ are pullbacks of their counterparts on $U_\et$}
\]
and that
\be \lab{GTI-2}
\cH^2_{\{t = 0\}}(-, G) \hra \cH^2_{\{t = 0\}}((-)_{A'}, G).
\ee
These assertions are local, so we may replace $A'$ (respectively, $A$) by its strict Henselization at a variable point of $U_{A'/tA'}$ (respectively, at its image in $U$) and assume that $U = \Spec(A)$ and $U' = \Spec(A')$ with $I = A$ and $I' = A'$. The $\cH^i_{\{ t = 0\}}$ are then simply $H^i_{\{ t = 0\}}$, and the assertion about $\cH^0_{\{t = 0\}}(-, G)$ follows from the following identification supplied by \Cref{BL-cute-cor}:
\[
\tst G(A) \isomto G(A') \times_{G(A'[\f{1}{t}])} G(A[\f{1}{t}])
\]
Letting $\wh{A}$ denote common $t$-adic completion of $A$ (or of $A'$), 
this identification gives the injections 
\[
\tst G(A[\f{1}{t}])/\im(G(A)) \hra G(A'[\f{1}{t}])/\im(G(A')) \hra G(\wh{A}[\f{1}{t}])/\im(G(\wh{A})).
\]
\Cref{Elk-Gab} (with \Cref{claim:const}) 
ensures that this composition is surjective, so both arrows are bijective. The five lemma and \ref{GTI-a} (respectively, and \Cref{BL-glue}~\ref{BLG-b}) then give the desired 
\[
\tst H^1_{\{t = 0\}}(A, G) \isomto H^1_{\{t = 0\}}(A', G)
\]
(respectively, reduce the remaining \eqref{GTI-2} to showing that $H^2(A, G) \hra H^2(A', G)$). For the latter, it suffices to apply \Cref{Hens-pair-inv}~\ref{HPI-b} to conclude that even the following composition is injective:
\[
H^2(A, G) \ra H^2(A', G) \ra H^2(A'/tA', G) \cong H^2(A/tA, G). \qedhere
\]
\epf

The technique we used for \Cref{G-tor-inv} also leads to the following nonaffine generalization of Gabber's affine analogue of proper base change theorem \cite{Gab94}*{Theorem 1}  and of its nonabelian analogue (the case $t = 1$ below recovers these affine versions). 
Related results appear in \cite{Fuj95}*{Corollary 6.6.4, Theorem 7.1.1}, \cite{ILO14}*{Expos\'{e} XX, Section 4.4}, and \cite{BM21}*{Corollary 1.18}, one distinction being the present setting of Gabber--Ramero triples with a possibly nontrivial ideal $I$. In the case when $U = \Spec(A[\f{1}{t}])$ and $t$ is a nonzerodivisor, the finite \'{e}tale aspect of \Cref{fet-algebraize} appears in \cite{GR03}*{Proposition 5.4.53}.

\bthm \label{nonaffine-Gab} \label{fet-algebraize}
For a map $(A, t, I) \ra (A', t, I')$ of Henselian, bounded Gabber--Ramero triples such that 
$A/t^nI \isomto A'/t^nI'$ for all $n > 0$, and an open 
\[
\tst \Spec(A[\f{1}{t}]) \subset U \subset \Spec(A),
\]
pullback gives an equivalence between the categories of finite \'{e}tale schemes over $U$ and $U_{A'}$ and 
\be \label{NAG-Hi}
R\Gamma_\et(U, \sF) \isomto R\Gamma_\et(U_{A'}, \sF) \qxq{for every torsion abelian sheaf} 
\sF\q \x{on $U_\et$},
\ee
in particular, for the closed subsets $Z \ce \Spec(A) \setminus U$ and $Z' \ce \Spec(A') \setminus U_{A'} \cong Z_{A'}$, 
\be \label{NAG-HiZ}
R\Gamma_Z(A, \sF) \isomto R\Gamma_{Z'}(A', \sF) \qxq{for every torsion abelian sheaf} 
\sF\q \x{on $A_\et$}.
\ee
\ethm

\bpf
By \Cref{nzd-example}, we have $\wh{A[\f{1}{t}]} \cong \wh{A'[\f{1}{t'}]}$, so  
\Cref{cor:finite-etale} (with \Cref{filter-eg}) gives the case $U = \Spec(A[\f{1}{t}])$. For the general case, we begin with the claim about finite \'{e}tale schemes and use patching as follows. By \Cref{BL-glue}~\ref{BLG-a} (with \Cref{nzd-example}), base change is an equivalence from the category of flat quasi-coherent $\sO_U$-modules $M$ to that of triples 
\[
\tst \p{M', M_t, \iota\colon M' \tensor_{\sO_{U_{A'}}} A'[\f{1}{t}] \isomto M_t \tensor_{A[\f{1}{t}]} A'[\f{1}{t}]}
\]
 consisting of a quasi-coherent, flat  $\sO_{U_{A'}}$-module $M'$, a flat $A[\f{1}{t}]$-module $M_t$, and the indicated isomorphism $\iota$. 
 Moreover, by \Cref{BL-glue}, the $\sO_U$-module $M$ is of finite type if and only if both $M \tensor_{\sO_U} \sO_{U_{A'}}$ and $M[\f{1}{t}]$ are of finite type as modules over $\sO_{U_{A'}}$ and $A[\f{1}{t}]$, respectively. This implies the same glueing for finite \'{e}tale algebras, which bootstraps the desired equivalence between the categories of finite \'{e}tale schemes over $U$ and $U_{A'}$ from the settled case $U = \Spec(A[\f{1}{t}])$. 

For the rest, we focus on \eqref{NAG-Hi} because, due to the cohomology with supports triangle, it implies \eqref{NAG-HiZ}. 
Moreover, its settled case $U = \Spec(A[\f 1t])$ 
reduces us to showing that
\[
H^i_{\{t = 0\}} (U, \sF) \cong H^i_{\{t = 0\}} (U_{A'}, \sF) \qxq{for every} i \in \bZ.
\]
The corresponding sheafified cohomologies with supports vanish away from the loci $\{t = 0\}$, and $U$ and $U_{A'}$ agree modulo $t$, so, due to the local-to-global spectral sequence \cite{SGA4II}*{Expos\'{e} V, Proposition 6.4} (which is the version of \Cref{lem-sup-ss} for the \emph{\'{e}tale} cohomology with supports), we reduce further to showing that, for every $i \in \bZ$, the sheaf
\be \label{NAG-desire}
 \cH^i_{\{t = 0\}} (-, \sF) \x{ on $(U_{A'})_\et$} \q \x{is the pullback of} \q \cH^i_{\{t = 0\}} (-, \sF) \x{ on $U_\et$}.
\ee
The pullback of the second sheaf maps to the first, 
so \eqref{NAG-desire} may be checked after replacing $U_{A'}$ (respectively, $U$) by its strict  Henselization at a point at which $t$ vanishes (respectively, at its image in $U$). 
This reduces us to the case when $U = \Spec(A)$ and $U' = \Spec(A')$, 
a case in which the invariance of \'{e}tale cohomology under Henselian pairs (a special case of \Cref{cor:finite-etale} and an input to its proof) and the settled case $U = \Spec(A[\f 1t])$ give the desired identifications
\[
H^i_{\{t = 0\}} (A, \sF) \isomto H^i_{\{t = 0\}} (A', \sF). \qedhere
\]
\epf




Parts \ref{BC-d} and \ref{BC-H2}  of the following corollary include the results announced in \cite{Gab93}*{Theorem 2.8~(i)}. For a generalization of the results announced in \cite{Gab93}*{Theorem 2.8~(ii)}, see \cite{flat-purity}*{Corollary~5.6.10}. 

\bcor \label{Br-comp}
For a map $A \ra A'$ of Noetherian rings, an ideal $J \subset A$ such that $(A, J)$ and $(A', JA')$ are Henselian \up{respectively, Zariski} pairs and $A/J^n \isomto A'/J^nA'$ for $n > 0$, and an open 
\[
\Spec(A) \setminus \Spec(A/J) \subset U \subset \Spec(A),
\]
we have
\benum
\m \label{BC-c}
for any quasi-affine, smooth $U$-group $G$,
\[
\qq H^1(U, G) \isomto H^1(U_{A'}, G) \q \x{{\upshape(}respectively,} \q H^1(U, \GL_n) \hra H^1(U_{A'}, \GL_n));
\]

\m \label{BC-d}
a pullback isomorphism 
\[
\qq \Br(U) \isomto \Br(U_{A'});
\]

\m  \lab{BC-H2}
for any commutative $U$-group $G$ such that either 
\benumr
\m
$G$ is quasi-affine and smooth over $U$\uscolon or 

\m
$G$ is finite and locally free over $U$\uscolon
\eenum
we have
\[
\qq H^2(U, G) \hra H^2(U_{A'}, G);
\]

\m \label{BC-b}
an equivalence of categories
\[
\qq U_\fet \isomto (U_{A'})_\fet;
\]

\m \label{BC-a}
for every torsion abelian sheaf $\sF$ on $U_\et$, pullback isomorphisms
\[
\qq H^i_\et(U, \sF) \isomto H^i_\et(U_{A'}, \sF) \qxq{for every}  i \in \bZ.
\]
\eenum
\ecor

\bpf
Since the map $A \ra A'$ is an isomorphism on $J$-adic completions, we lose no generality by assuming that $A'$ is the $J$-adic completion $\wh{A}$ of $A$. Then,
as in the proof of \Cref{thm:Elkik-reproof}, by forming the completion iteratively with respect to a fixed system of generators for $J$, we reduce to $J$ being principal. In this case the assertions amount to special cases of \Cref{G-tor-inv,nonaffine-Gab}.
\epf

\brem
When $A$ is excellent, the injectivity aspects of \Cref{Br-comp} become straight-forward consequences of \Cref{xexc-comp} that is based on the N\'{e}ron--Popescu approximation.
\erem 



\section{Torsors over $R\llp t \rrp$} \label{chapter:loop-torsors}

We enter a more detailed study of torsors over the base $R\llp t \rrp$. In \S\ref{section:H1-description}, we present a complete description of such torsors when the group is a torus defined over $R$. In \S\ref{section:tame-tori}, we complement this with a crucial for our goals vanishing statement when the torus need not descend to $R$. 

\csub[A formula for $H^1(R\llp t \rrp, T)$ for any $R$-torus $T$] \label{section:H1-description}

\numberwithin{equation}{subsection}

By a result of Weibel \cite{Wei91}, for every commutative ring $R$ we have
\be \label{Wei-main}
\addtocounter{subsubsection}{1}
\Pic(R) \oplus \Pic(R[t])_0 \oplus \Pic(R[t\i])_0 \oplus H^1_\et(R, \bZ) \isomto \Pic(R [ t, t\i ] ),
\ee 
where $(-)_0$ denotes the classes that die in $\Pic(R)$ under $t \mapsto 0$ (respectively, $t\i \mapsto 0$) and the map $H^1_\et(R, \bZ) \ra \Pic(R [ t, t\i ] )$ is obtained from  $\bZ \xra{1\, \mapsto\, t} \bG_{m,\, R[t,\, t\i]}$. For applying this formula, it is sometimes useful to compute $H^1_\et(R, \bZ)$ in the Nisnevich topology, which is possible because every $\bZ$-torsor over a Henselian local ring is trivial (see \Cref{Hens-pair-inv}~\ref{HPI-a} and the proof of \Cref{cor:algebraize-loop-Pic}).

Our goal in this section is to expose a similar description for $\Pic(R \llp t \rrp )$ due to Gabber \cite{Gab19}:
\[
\Pic(R) \oplus \Pic(R[t\i])_0 \oplus H^1_\et(R, \bZ) \isomto \Pic(R \llp t \rrp );
\]
in fact,  we present a mild strengthening valid for torsors under any $R$-torus $T$, see \Cref{Gabber}. The results of this section originate in a letter of Gabber \cite{Gab19} to the first named author.

\numberwithin{equation}{subsubsection}

We begin by reviewing the Weierstrass preparation theorem in a Henselian setting. 

\bprop \label{Weierstrass}
For a Henselian local ring $(R, \fm)$, every $f \in R\{ t \} \setminus \fm(R\{t\})$  is of the form 
\[
f = Pu \q \x{for a unique} \q P = t^d + a_{d - 1} t^{d - 1} + \cdots + a_0 \q \x{with} \q a_i \in \fm  \q \x{and a unique} \q  u \in R\{t \}^\times;
\]
moreover, the natural map $R[t]/(P) \ra R\{t\}/(f)$ is an isomorphism. 
\eprop

\bpf
The functor $R \mapsto R\{t\}$ commutes with filtered direct limits (see \S\ref{hens-series}), so we lose no generality by assuming that $R$ is Noetherian. The uniqueness can then be seen over the $(\fm, t)$-adic completion $\wh{R}\llb t \rrb$ of $R\{ t \}$, which satisfies the usual Weierstrass preparation theorem \cite{BouAC}*{Chapitre~VII, section 3, num\'{e}ro 8,~Proposition 6}. 

For the existence, we begin by setting $k \ce R/\fm$, so that 
\[
k \{ t \} \cong R\{t\}/\fm(R\{ t \})
\]
(see \S\ref{hens-series}). Every nonempty closed subscheme of $\Spec(R\{t\})$ meets $\Spec(k \{ t \})$ and $\Spec(k\{t\}/f k \{t\})$ is Artinian local supported along $\{t = 0\}$. Thus, isolation of quasi-finite parts of finite type schemes over Henselian local bases \cite{EGAIV4}*{Th\'{e}or\`{e}me 18.5.11 c)} applied to $R'/(f)$, where $R[t] \ra R'$ is a sufficiently small \'{e}tale neighborhood of the zero section of $R[t]$, implies that $R\{t\}/(f)$ is a finite local $R$-algebra that is isomorphic to a factor of $R'/(f)$. In particular, the $\fm$-adic topology of $R\{t\}/(f)$ agrees with its $(\fm, t)$-adic topology, to the effect that
\[
R\{t\}/(f) \tensor_R \wh{R} \isomto R\{t\}/(f) \tensor_{R\{t\}} \wh{R}\llb t \rrb \cong \wh{R}\llb t \rrb/f\wh{R}\llb t \rrb. 
\]
However, $\wh{R}\llb t \rrb/f \wh{R}\llb t \rrb$ is $\wh{R}$-flat by Weierstrass preparation \cite{BouAC}*{Chapitre VII, section 3, num\'{e}ro~8, Proposition 5}, so $R\{t\}/(f)$ must be $R$-flat, and hence even finite free of some rank $d \ge 0$ as an $R$-module. Let $P \in R[t]$ be the characteristic polynomial of the $R$-linear scaling by $t$ action on $R\{t\}/(f)$. By construction, $P$ is of the form claimed in the statement, and the Cayley--Hamilton theorem \cite{BouA}*{Chapitre III, section 8, num\'{e}ro 11, Proposition 20} supplies a map $i \colon R[t]/(P) \ra R\{t\}/(f)$. By Weierstrass preparation \cite{BouAC}*{Chapitre VII, section 3, num\'{e}ro 8, Proposition 5} again, $i \tensor_R \wh{R}$ is an isomorphism, so $i$ is too. In the resulting expression $P = fv$ one must have $v \in R\{ t \}^\times$: indeed, since $i$ is an isomorphism, $v$ must reduce to a  unit modulo $\fm$, that is, to a unit of the discrete valuation ring $k\{t\}$. 
\epf

The following proposition and the subsequent \ref{GL-c} and \ref{hens-seminormal} capture some concrete geometric consequences of the Henselian Weierstrass preparation established in \Cref{Weierstrass}.

\bprop \label{Gab-lem}
Let $(R, \fm)$ be a Henselian local ring with the residue field $k \ce R/\fm$.
\benum
\item \label{GL-a}
Base change and schematic image give inverse bijections between the set of closed subschemes of $\bP^1_R$ that do not meet $\Spec(k[t\i])$ and the set of closed subschemes of $\Spec(R[t]_{1 + tR[t]})$ \up{respectively, of $\Spec(R\{t \})$} that do not contain the $k$-fiber. Moreover, for any such closed $Z \subset \bP^1_R$,
\[
\qq Z \isomto Z_{R[t]_{1 + tR[t]}} \isomto  Z_{R\{ t \}}.
\]

\item \label{GL-b}
Base change and schematic image give inverse bijections between the set of closed subschemes of $\Spec(R[t\i])$ that do no meet $\Spec(k[t\i])$ and the set of closed subschemes of $\Spec((R[t]_{1 + tR[t]})[\f{1}{t}])$ \up{respectively, of $\Spec(R\{t \}[\f{1}{t}])$} that do not meet the $k$-fiber.\footnote{For the sake of concreteness, we recall from \Cref{Lam-cor} that the $t$-localized Zariskization $(R[t]_{1 + tR[t]})[\f{1}{t}]$ of $R[t]$ along $\{t = 0\}$ is also isomorphic to the localization of $R[t]$ with respect to the monic polynomials.} 

\item \label{GL-d}
The bijections of \ref{GL-a} and \ref{GL-b} respect the property of being an effective Cartier divisor.
\eenum
\eprop

\bpf 
All the morphisms for which we consider schematic images are quasi-compact and quasi-separated, so the formation of these images exhibits no pathologies, for instance, it commutes with flat base change, see \cite{EGAI}*{Corollaire 9.5.2} and \cite{EGAIV1}*{num\'{e}ro 1.7.8}.

\benum
\item
Firstly, the maps in question supplied by base change are injective: indeed, $\Spec(R[t]_{1 + tR[t]})$ is the Zariskization of $\bP^1_R$ along $\{t = 0\}$, so it contains all the generizations in $\bP^1_R$ of the origin of $\bP^1_k$, and $R\{t \}$ is faithfully flat over $R[t]_{1 + tR[t]}$. 
Moreover, the ideal of $R\{t \}$ that cuts out a closed subscheme $Z' \subset \Spec(R\{t \})$ not containing the $k$-fiber must contain an $f \in R\{ t \} \setminus \fm(R\{t\})$, so \Cref{Weierstrass} implies that $Z'$ is the base change of a closed subscheme $Z \subset \bP^1_R$ that does not meet $\Spec(k[t\i])$ (note that $R[T]/(P)$ is $R$-finite for any monic $P$ and hence defines a closed subscheme of $\bP^1_R$). Consequently, the maps supplied by base change are even bijective, their inverses are given by forming schematic images, and 
\[
\qq Z \isomto Z_{R[t]_{1 + tR[t]}} \isomto Z_{R\{ t \}}.
\]

\item
Via schematic images in $\bP^1_R$, the closed subschemes $Z' \subset \Spec(R[t\i])$ in question correspond to those closed subschemes $Z \subset \bP^1_R$ as in \ref{GL-a} for which $t$ is a nonzerodivisor on $Z$, and likewise for $\Spec((R[t]_{1 + tR[t]})[\f{1}{t}])$ or $\Spec(R\{t \}[\f{1}{t}])$. Thus, \ref{GL-a} gives the claim by using its isomorphy aspect to handle the nonzerodivisor requirement.

\item
The isomorphy aspect of \ref{GL-a} ensures that the bijections respect the property of being of finite presentation over the respective ambient schemes because the latter is equivalent to finite presentation over $R$. Thus, by faithfully flat descent, it then also respects the further property of being cut out by a nonzerodivisor at every stalk of the ambient space. 
 \qedhere
\eenum
\epf

\bcor \label{GL-c}
For a Henselian local ring $(R, \fm)$, pullback maps induce isomorphisms 
\[
\tst \Pic(R[t\i]) \isomto \Pic((R[t]_{1 + tR[t]})[\f{1}{t}]) \isomto \Pic(R\{t \}[\f{1}{t}]) \isomto \Pic(R\llp t \rrp).
\]
\ecor

We will globalize \Cref{GL-c} in \eqref{eqn:Pic-loop-global} below.

\bpf
The last map is an isomorphism by \Cref{Hi-algebraize}, so we focus on the other two. Moreover, by limit arguments, we may assume that $R$ is Noetherian. By \Cref{BL-glue}, any line bundle on $R[t\i]$ that trivializes over $R\{t\}[\f{1}{t}]$ extends to a line bundle on $\bP^1_R$, and so it must be trivial by, for instance, \cite{Lam06}*{Chapter IV, Horrocks' Theorem 2.2}. Thus, the first map and also its composition with the second one are injective. Since every line bundle on $R\{t\}[\f{1}{t}]$ or $(R[t]_{1 + tR[t]})[\f{1}{t}]$ is associated to an effective Cartier divisor that may be chosen to not meet the special fiber (or any finite set of points, see \cite{SP}*{Lemma \href{https://stacks.math.columbia.edu/tag/0AYM}{0AYM}}), the remaining surjectivity assertion follows from \Cref{Gab-lem}~\ref{GL-d}.
\epf

\bcor \label{hens-seminormal}
For a seminormal, Henselian local ring $(R, \fm)$, we have $\Pic(R\llp t \rrp) = 0$.
\ecor

\bpf
It suffices to combine \Cref{GL-c} with the seminormality criterion \cite{Swa80}*{Theorem~1}.
\epf

\Cref{hens-seminormal} is a special case of a general formula for $H^1(R\llp t \rrp, T)$, which we present in \Cref{Gabber}. For this, we begin with the following basic description of the units of $R\{t \}[\f{1}{t}]$ that simultaneously records the basic structure of the affine Grassmannian of a torus.

\blem \label{lem:unit-description}
For a reduced ring $R$ and an $R$-torus $T$, we have compatible isomorphisms
\be \label{eqn:unit-desc}
\ba\tst (X_*(T))(R) \times T(R) &\isomto  T(R[t, t\i]), \tst\qq (X_*(T))(R) \times T(R\{t \}) \isomto T(R\{t \}[\f{1}{t}]),\\
\tst &\q (X_*(T))(R) \times T(R\llb t \rrb) \isomto T(R\llp t \rrp),  
\ea
\ee
where on $X_*(T) \cong \Hom(\bG_m, T)$ the maps are obtained from $\bZ \xra{1\, \mapsto\, t} \bG_{m,\, R[t,\, t\i]}$ via $\Hom(\bZ, T) \cong T$.
\elem

\bpf
Since $T$ is affine and $R[t, t\i] \subset R\{t\}[\f{1}{t}] \subset R\llp t \rrp$ (see \eqref{eqn:hens-series}), we have compatible inclusions
\[
\tst T(R[t]) \subset T(R\{t\}) \subset T(R\llb t \rrb) \qxq{and} T(R[t, t\i]) \subset T(R\{t\}[\f{1}{t}]) \subset T(R\llp t \rrp). 
\]
Moreover, by \eqref{BL-cute-eq}, 
\[
\tst T(R[t]) \isomto T(R[t, t\i]) \cap T(R\llb t \rrb) \qxq{and}T(R\{t\}) \isomto T(R\{t\}[\f{1}{t}]) \cap T(R\llb t \rrb) \qxq{in} T(R\llp t \rrp). 
\]
Thus, since $T(R) \isomto T(R[t])$ (as may be checked \'{e}tale locally on $R$), we reduce to the case of $R\llp t \rrp$. 

For an \'{e}tale cover $R \ra R'$, we see coefficientwise that both
\[
R\llb t \rrb \ra R'\llb t \rrb \rightrightarrows (R'\tensor_R R')\llb t \rrb \qxq{and} R\llp t \rrp \ra R'\llp t \rrp \rightrightarrows (R'\tensor_R R')\llp t \rrp
\]
are equalizer diagrams. Consequently, again since $T$ is affine, both functors $R \mapsto T(R\llb t\rrb)$ and $R \mapsto T(R\llp t\rrp)$ are \'{e}tale sheaves in $R$. Thus, we may work \'{e}tale locally on $R$ and assume that~$T = \bG_m$. 

For an $f \in R\llp t \rrp^\times$, the function $\ord_f \colon \Spec(R) \ra \bZ$ that maps a prime $\fp \subset R$ to the $t$-adic valuation of the image of $f$ in $k(\fp)\llp t \rrp$ is upper semicontinuous: each $\fp$ has an open neighborhood in $\Spec(R)$ on which $\ord_f$ is $\le \ord_f(\fp)$ (concretely, a neighborhood on which the coefficient of $t^{\ord_f(\fp)}$ is a unit). Moreover, $\ord_f + \ord_{f\i}$ is identically zero. Thus, $\ord_f$ is in fact locally constant, so there is a unique $R$-point $n$ of $X_*(\bG_m) \cong \underline{\bZ}$ such that $\ord_{t^{n}f}$ is identically zero. Consequently, at every residue field of $R$ the coefficient of $t^m$ of $t^nf$ vanishes for $m < 0$ (respectively, is a unit for $m = 0$). Thus, this coefficient lies in the nilradical and hence vanishes (respectively, is a unit) in $R$, so that $ t^{n}f \in R\llb t \rrb^\times$.
\epf

\bthm \label{Gabber}
For a ring $R$ and an $R$-torus $T$, we have the identification
\be \label{eqn:Gabber}
\tst H^1(R[t\i], T) \oplus H^1(R, X_*(T)) \isomto H^1(R \{ t \}[\f{1}{t}], T ) \overset{\x{\uref{Hi-algebraize}}}{\cong} H^1(R\llp t \rrp, T),  
\ee
where the map on $H^1(R, X_*(T))$ is obtained from $\bZ \xra{1\, \mapsto\, t} \bG_{m,\, R[t,\, t\i]}$ as in Lemma~\uref{lem:unit-description}\uscolon 
thus, 
\[
H^1(R, T) \hra H^1(R\llp t \rrp, T)
\]
and, for reduced $R$,
\be \label{eqn:H1-support-zero}
\xymatrix{(X_*(T))(R) \ar[r]_-{\sim}^-{\ref{lem:unit-description}} & H^1_{\{t = 0\}}(R\{t\}, T) \cong H^1_{\{t = 0\}}(R[t], T) } \q \x{functorially in $R$ and $T$.}
\ee
\ethm

\bpf
We focus on \eqref{eqn:Gabber}---the rest will follow: for $H^1(R, T) \hra H^1(R\llp t \rrp, T)$, one decomposes
\[
\tst H^1(R[t\i], T) \cong H^1(R, T) \oplus \Ker\p{H^1(R[t\i], T) \xra{t\i\, \mapsto\, 0} H^1(R, T)},
\]
and for \eqref{eqn:H1-support-zero}, one first deduces that $H^1(R\{t\}, T) \hra H^1(R \{ t \}[\f{1}{t}], T )$ (see \Cref{Hens-pair-inv}), and then combines the cohomology with supports sequence with \Cref{lem:unit-description} and excision \cite{Mil80}*{Chapter~III, Proposition 1.27}.

For \eqref{eqn:Gabber} itself, by limit arguments based on \cite{SGA4II}*{Expos\'{e} VII, Corollaire 5.9}, we may assume that $R$ is Noetherian. The topological invariance of the \'{e}tale site and the insensitivity of torsors under smooth groups to the nilradical \cite{SGA3IIInew}*{Expos\'{e} XXIV, Lemme 8.1.8} then allow us to replace $R$ by $R^\red$. 
We consider the functor $\sF$ defined on the category of \'{e}tale $R$-algebras by
\[
\sF\colon R' \mapsto \{ (X, \iota\colon X|_{\{t\i = 0\}} \isomto T_{R'})\}/\sim,
\]
where $X$ is a $T_{R'[t\i]}$-torsor and $\iota$ is a trivialization of its pullback to $R'$ along $t\i \mapsto 0$. Since $R'$ is reduced, we see \'{e}tale locally on $R'$ that $T(R') \isomto T(R'[t\i])$, to the effect that the pairs $(X, \iota)$ have no nonidentity automorphisms. Thus, by descent, $\sF$ is an \'{e}tale sheaf on $R$ with global sections 
\[
\tst \sF(R) \cong \Ker\p{H^1(R[t\i], T) \xra{t\i\, \mapsto\, 0} H^1(R, T)}.
\]
By \Cref{GL-c} and limit arguments, for any $T_{R \{ t \}[\f{1}{t}]}$-torsor $\cX$ there is an \'{e}tale cover $R \ra R'$ such that the isomorphism class of $\cX_{R' \{ t \}[\f{1}{t}]}$ lifts to a unique $x \in \sF(R')$. The unique lift continues to exist over any refinement of the cover, so, since $\cX$ begins life over $R \{ t \}[\f{1}{t}]$ and $\sF$ is an \'{e}tale sheaf, the local lifts glue to a unique global 
\[
x \in \sF(R) \subset H^1(R[t\i], T).
\]
By adjusting $\cX$ by the pullback of the $T_{R[t\i]}$-torsor determined by $x$, we reduce to the case when $x = 0$. Such $\cX$ are precisely the $T_{R \{ t \}[\f{1}{t}]}$-torsors that trivialize over $(R_\fp^\sh)\{t\}[\f{1}{t}]$ for every prime $\fp \subset R$. By viewing $\fp$ as a prime of $R\{t\}$ that contains $t$, we have 
\be\label{eqn:stalk-id}
(R_\fp^\sh)\{t\} \isomto R\{t\}_\fp^\sh, \q \x{compatibly with $t$ on both sides}
\ee
(both sides are initial among the strictly Henselian local $R_\fp$-algebras $(A, \fm)$ equipped with a local map $R_\fp \ra A$ and a $t \in \fm$). Thus, the $\cX$ as above are precisely the $T_{R \{ t \}[\f{1}{t}]}$-torsors that trivialize over the pullback of some \'{e}tale cover of $R\{t \}$, in other words, by the characterization of torsors under pushforward groups \cite{Gir71}*{Chapitre V, Proposition 3.1.3}, they are precisely the $(j_*(T_{R \{ t \}[\f{1}{t}]}))$-torsors where 
\[
\tst j \colon \Spec(R \{ t \}[\f{1}{t}]) \hra \Spec(R\{t \}).
\]
We have reduced to showing that
\be \label{eqn:H1-jstar}
\tst H^1(R, T) \oplus H^1(R, X_*(T)) \isomto H^1(R\{t\}, j_*(T_{R \{ t \}[\f{1}{t}]})) \qxq{inside} H^1(R\{t\}[\f{1}{t}], T).
\ee
For this,  we first exhibit the following short exact sequence on $R\{t\}_\et$: 
\be \label{eqn:want-SES}
0 \ra T_{R\{t\}} \ra  j_*(T_{R \{ t \}[\f{1}{t}]}) \ra i_*(X_*(T)) \ra 0,
\ee
where $i \colon \Spec(R) \hra \Spec(R\{t \})$ is the closed immersion complementary to $j$. To explain the third map, we first note that, by the insensitivity of the \'{e}tale site to nonreduced structure, the $t$-adic completion of any \'{e}tale $R\{ t\}$-algebra $A$ is canonically and functorially isomorphic to $(A/(t))\llb t \rrb$. Thus, the third map is well-defined by letting it on $A$-points be the composition 
\[
\tst T(A[\f 1t]) \ra T((A/(t))\llp t \rrp) \xra{\eqref{eqn:unit-desc}} (X_*(T))(A/(t)).
\]
This makes \eqref{eqn:want-SES} short exact because, by \Cref{Elk-Gab} (with \Cref{claim:const}), an element of $(X_*(T))(A/(t))$ is in the image of $T(A[\f 1t])$ after pullback along an \'{e}tale $A \ra A'$ with $A/(t) \isomto A'/(t)$. 

We claim that \eqref{eqn:want-SES} induces split short exact sequences on $H^0$ and $H^1$. For $H^0$, this is part of \Cref{lem:unit-description}. For $H^1$, the composition of the evaluation at $t$ map $X_*(T)_{R\{t\}} \ra j_*(T_{R \{ t \}[\f{1}{t}]})$ used there with the third map of \eqref{eqn:want-SES} is, by construction, the adjunction map $X_*(T)_{R\{t\}} \ra i_*(X_*(T))$, and it then suffices to argue that the latter induces an isomorphism 
\[
H^1(R\{t\}, X_*(T)) \isomto H^1(R, X_*(T)).
\]
Indeed, this last map is surjective because it has a section, and it is injective by \Cref{Hens-pair-inv}~\ref{HPI-a}. Finally, to deduce that the split exact sequence on $H^1$ gives \eqref{eqn:H1-jstar}, we use \Cref{Hens-pair-inv} again.
\epf

\brem
We recall from \cite{CTS87}*{Lemma 2.4} (with a limit argument that eliminates the Noetherianity hypothesis) that if $R$ is a normal integral domain, then 
\[
H^1(R, T) \isomto H^1(R[t\i], T),
\]
so that
\[
H^1(R, T) \oplus H^1(R, X_*(T)) \overset{\eqref{eqn:Gabber}}{\cong} H^1(R\llp t \rrp, T).
\]
\erem

In the following consequence of \Cref{Gabber} the ring $R$ could, for instance, be a normal domain.

\bcor \label{cor:algebraize-loop-Pic}
For a geometrically unibranch ring $R$ whose spectrum is irreducible, 
\[
\Pic(R[t\i]) \isomto \Pic(R\llp t \rrp) \qxq{{\upshape(}respectively,} \Pic(R) \isomto \Pic(R\llp t \rrp) \q \x{if $R$ is also seminormal}).
\]
In particular, for a local normal domain $R$, we have $\Pic(R\llp t \rrp) \cong 0$. 

\ecor

\bpf
\Cref{Gabber} (respectively, and the characterization of seminormality given in \cite{Swa80}*{Theorem 1}) reduces us to showing that $H^1_\et(R, \bZ) = 0$. However, by the characterization of components of \'{e}tale schemes over normal bases \cite{EGAIV4}*{Proposition 18.10.7}, the sheaf $\bZ$ on $R_\et$ is the pushforward of $\bZ$ from the generic point of $\Spec(R)$. Thus, letting $K$ denote the residue field at this generic point, we get that 
\[
H^1_\et(R, \bZ) \subset H^1_\et(K, \bZ) \cong \Hom_\cont(\Gal(\ov{K}/K), \bZ) = 0. \qedhere
\]
\epf

We conclude the section by generalizing the formula of Weibel \cite{Wei91} for $\Pic(R[t, t\i])$ recalled in \eqref{Wei-main} to a formula for $H^1(R[t, t\i], T)$ valid for any $R$-torus $T$.

\bthm \label{thm:Wei-torus}
For a ring $R$ and an $R$-torus $T$, we have
\be \label{eqn:Wei-torus}
H^1(R, T) \oplus H^1(R[t], T)_0 \oplus H^1(R[t\i], T)_0 \oplus H^1(R, X_*(T)) \isomto H^1(R [ t, t\i ], T),
\ee
where $(-)_0$ denotes the kernel of the evaluation $t \mapsto 0$ \up{respectively, $t\i \mapsto 0$} and the map on $H^1(R, X_*(T))$ is obtained from $\bZ \xra{1\, \mapsto\, t} \bG_{m,\, R[t,\, t\i]}$ analogously to the map in Lemma \uref{lem:unit-description}.
\ethm

\bpf
\Cref{Hens-pair-inv,Gabber} with cohomology sequences supply a commutative diagram
\[
\xymatrix@C=10pt@R=12pt{
0 \ar[r] & H^1(R[t], T) \ar[r] \ar@{->>}[d]^-{t\,\mapsto\, 0} & H^1(R[t, t\i], T) \ar[d] \ar@{->>}[r] & H^1(R[t, t\i], T)/H^1(R[t], T) \ar@{^(->}[r] \ar[d] & H^2_{\{t = 0 \}}(R[t], T) \ar[d]^{\sim} \\
0 \ar[r] & H^1(R, T) \ar[r] & H^1(R\{t\}[\f{1}{t}], T) \ar@{->>}[r] & H^1(R[t\i], T)_0 \oplus H^1(R, X_*(T)) \ar@{^(->}[r] & H^2_{\{t = 0 \}}(R\{t\}, T),
}
\]
where we used excision \cite{Mil80}*{Chapter III, Proposition 1.27} for the last vertical arrow. Consequently, the third vertical arrow is injective and, since the splitting 
\[
\tst H^1(R[t\i], T)_0 \oplus H^1(R, X_*(T)) \ra H^1(R\{t\}[\f{1}{t}], T)
\]
supplied by \Cref{Gabber} naturally factors through $H^1(R[t, t^{-1}], T)$, it is also surjective. The factorization then ensures that the top short exact sequence splits, and the desired \eqref{eqn:Wei-torus} follows. 
\epf

With the argument of \Cref{thm:Wei-torus}, we now also describe $T$-torsors over $R(t)\overset{\ref{Lam-cor}}{\simeq} (R [ t]_{1 + tR[t]})[\f{1}{t}]$.

\begin{variant} \label{var:H1-torus}
For a ring $R$ and an $R$-torus $T$, 
\[
\tst H^1(R[t]_{1 + tR[t]}, T) \oplus H^1(R[t\i], T)_0 \oplus H^1(R, X_*(T)) \isomto H^1((R [ t]_{1 + tR[t]})[\f{1}{t}], T);
\]
in particular, 
\be \label{eqn:Pic-loop-global}
\tst \Pic(R[t\i]) \oplus H^1(R, \bZ) \isomto \Pic((R [ t]_{1 + tR[t]})[\f{1}{t}]) \cong \Pic(R\{t\}[\f{1}{t}]) \cong \Pic(R\llp t \rrp). 
\ee
\end{variant}

\bpf
Example \uref{sta-free-eg} ensures that $\Pic(R) \isomto \Pic(R[t]_{1 + tR[t]})$, so the second assertion follows from the first and \Cref{Gabber}. For the first assertion, the argument is as for \Cref{thm:Wei-torus}. 
\epf


\csub[Torsors under tame isotrivial tori over $R\llp t \rrp$] \label{section:tame-tori}

We wish to extend the vanishing that we saw in \Cref{hens-seminormal} to tori over $R\llp t \rrp$ that need not come from $R$, see \Cref{twisted} below for a precise statement. This will be of central importance for the product formula for the Hitchin fibration in the proof of \Cref{thm:product-formula}. 
The arguments in this section build on the ones explained by Gabber in a conversation with the first named author.

\blem \label{lem:Abhyankar-input}
For a strictly Henselian local ring $(R, \fm)$, the finite \'{e}tale Galois covers 
\[
R\llp t \rrp \ra R\llp t^{1/d} \rrp \qx{with $d$ invertible in $R$}
\]
are cofinal among the tamely ramified relative to $R$ finite \'{e}tale $R\llp t \rrp$-algebras\uscolon\!\!\footnote{We say that a finite \'{e}tale $R\llp t \rrp$-algebra $S$ is \emph{tamely ramified relative to $R$} if for every residue field $k$ of $R$, the \'{e}tale algebra $S \tensor_{R\llp t \rrp} k\llp t \rrp$ is tamely ramified over the discretely valued field $k\llp t \rrp$, compare with the definition of tameness given in \cite{SGA1new}*{Expos\'{e} XIII, D\'{e}finition 2.1.1}. 
} 
in particular, the $R$-tame \'{e}tale fundamental group of $R\llp t \rrp$ is $\prod_{\ell \neq \Char(R/\fm)} \bZ_\ell(1)$ and, via base change, the category of tamely ramified relative to $R$ finite \'{e}tale $R\llp t \rrp$-algebras is equivalent to its counterpart for $R/\fm$.
\elem

\bpf
\Cref{cor:stays-connected} ensures that $\Spec(R\llp t \rrp)$ is connected, so the degrees of its finite \'{e}tale covers are constant. Moreover, $R\llp t \rrp \ra R\llp t^{1/d} \rrp$ is indeed finite \'{e}tale Galois with group $\mu_d$. By \Cref{cor:finite-etale}, we may replace $R\llp t \rrp$ and $R\llp t^{1/d} \rrp$ by $R\{ t \}[\f 1t]$ and $R\{ t^{1/d} \}[\f 1t]$. By the relative Abhyankar's lemma \cite{SGA1new}*{Expos\'{e} XIII, Proposition 5.5}, for any tamely ramified relative to $R$ finite \'{e}tale cover, its base change to $R\{ t^{1/d} \}[\f 1t]$ extends to a finite \'{e}tale cover of $R\{ t^{1/d} \}$ for some $d$ invertible in $R$. However, $R\{ t^{1/d} \}$ is strictly Henselian local, so its \'{e}tale covers split. All the assertions follow.
\epf

\Cref{lem:Abhyankar-input} leads to the following result on cohomology.

\bprop \label{prop:pass-to-k}
For a strictly Henselian local ring $R$ with residue field $k$ and a finite, \'{e}tale $R\llp t\rrp$-group scheme $G$ of order invertible in $R$ such that $G$ becomes constant over a finite \'{e}tale cover of $R\llp t \rrp$ that is  tamely ramified relative to $R$, we have
\[
H^0(R\llp t \rrp, G) \isomto H^0(k\llp t \rrp, G) \qxq{and} H^1(R\llp t \rrp, G) \isomto H^1(k\llp t \rrp, G);
\]
if $G$ is commutative, then, in addition,
\[
R\Gamma(R\llp t \rrp, G) \isomto R\Gamma(k\llp t \rrp, G) \qxq{and} H^{i}(R\llp t \rrp, G) \cong 0 \qxq{for} i \ge 2.
\]
\eprop

\bpf
Since $G$ is tamely ramified relative to $R$, the statement about $H^0$ follows from the equivalence of categories of \Cref{lem:Abhyankar-input}. For $H^1$, we first note that if $G$ is constant, then $G$-torsors amount to $G$-conjugacy classes of homomorphisms $\pi_1^\et(R\llp t \rrp) \ra G(R\llp t \rrp)$, and likewise over $k\llp t\rrp$. Due to the assumption on the order of $G$, every such homomorphism factors through the maximal tame relative to $R$ quotient of $\pi_1^\et(R\llp t \rrp)$, and likewise over $k\llp t \rrp$, so \Cref{lem:Abhyankar-input} gives the statement about $H^1$ for constant $G$. For general tamely ramified $G$, this argument shows that every $G$-torsor splits over some tame relative to $R$ cover $R\llp t \rrp \ra R\llp t^{1/d}\rrp$ that may be chosen so that $G$ is constant over it, and likewise over $k\llp t \rrp$. Thus, \Cref{lem:Abhyankar-input} also gives the general case of the statement about $H^1$: now $G$-torsors amount to twisted conjugacy classes of twisted homomorphisms $\pi_1^\et(R\llp t \rrp) \ra G(R\llp t^{1/d}\rrp)$, and likewise over $k\llp t \rrp$, compare with \cite{Ser02}*{Chapter I, Sections 5.1--5.2} or \cite{Gir71}*{Chapitre III, Section 3.7}. 

For the rest of the proof, we assume that $G$ is commutative. Since $k$ is separably closed, the \'{e}tale cohomological dimension of the field $k\llp t \rrp$ is $\le 1$, so the vanishing assertion follows from the assertion about $R\Gamma$. For the latter, as we saw, by \Cref{lem:Abhyankar-input}, our group $G$ becomes constant over the finite \'{e}tale cover $R\llp t \rrp \ra R\llp t^{1/d} \rrp$ for some $d$ that is invertible in $R$. Thus, we may find a resolution 
\[
0 \ra G \ra G_1 \ra G_2 \ra \dotsc 
\]
of $G$ in which each $G_i$ is a finite \'{e}tale group obtained from some constant group over $R\llp t^{1/d} \rrp$ by restriction of scalars. The cohomology spectral sequence associated to this resolution allows us to replace $G$ by $G_i$, so we only need to show that
\[ 
 H^j(R\llp t^{1/d}\rrp, \bZ/n\bZ) \isomto H^j(k\llp t^{1/d} \rrp, \bZ/n\bZ) \qxq{for all} j \in \bZ \qx{and $n$ invertible in $R$.}
\]
The case $j \le 1$ follows from the previous paragraph. 
For $j \ge 2$, both sides 
vanish by  the relative cohomological purity \cite{SGA4III}*{Expos\'{e} XVI, Th\'{e}or\`{e}me 3.7} (with \Cref{cor:finite-etale} to replace $\llp -\rrp$ by $\{ -\}[\f 1t]$).
\epf

\beg
Let $R$ be a strictly Henselian local ring, let $n \in \bZ_{> 0}$ be invertible in $R$, and consider $G \ce \bZ/n\bZ \simeq \mu_n$. \Cref{prop:pass-to-k} implies that $H^1(R\llp t \rrp, \mu_n)$ is cyclic of order $n$, with a generator given by the image of $t$ under the connecting Kummer map $R\llp t \rrp^\times\!/R\llp t \rrp^{\times n} \hra H^1(R\llp t \rrp, \mu_n)$. It follows that for such $R$ the Picard group $\Pic(R\llp t \rrp)$ has no torsion of order invertible in $R$.
\eeg

\bthm \label{twisted}\label{galoisi}
For a seminormal, strictly Henselian, local ring $R$ and a $R\llp t\rrp$-torus $T$ that splits over $R\llp t^{1/d}\rrp$ for some $d \in \bZ_{> 0}$ that is invertible in $R$ \up{for instance, a $T$ that splits over some $W$-torsor over $R\llp t \rrp$ for a finite group $W$ whose order is invertible on $R$}, we have
\[
H^1(R\llp t \rrp, T) = 0.
\]
\ethm

\bpf
\Cref{lem:Abhyankar-input} implies the parenthetical example.  
Moreover, by \Cref{hens-seminormal}, we have 
\[
\tst H^1(R\llp t^{1/d} \rrp, T) \cong 0.
\]
Thus, by using the trace map constructed in \cite{SGA4III}*{Expos\'{e} XVII, Section 6.3.13 up to Proposition~6.3.15, especially, Proposition 6.3.15 (iv)}, we conclude that $H^1(R\llp t \rrp, T)$ is killed by $d$, so that 
\[
H^1(R\llp t \rrp, T) \hra H^2(R\llp t \rrp, T[d]).
\]
Since $\mu_d \simeq \bZ/d\bZ$ over $R$, \Cref{prop:pass-to-k} implies that the target of this injection vanishes. 
\epf

\section{Reductive group schemes and the Hitchin fibration}

We turn to the geometry and arithmetic of reductive group schemes. In \S\ref{section:Chevalley-iso}, we present a new proof of the Chevalley isomorphism valid over arbitrary bases. In \S\ref{section:Chevalley-Kostant}, we record several improvements that concern the geometry of the Chevalley morphism. Both \S\ref{section:Chevalley-iso} and \S\ref{section:Chevalley-Kostant} simultaneously build up the setup for the product formula for the Hitchin fibration, which we finally take up in \S\ref{section:Hitchin-fibration}.

\csub[The Chevalley isomorphism for root-smooth reductive group schemes] \label{section:Chevalley-iso}

The goal of this section is to show in \Cref{adj} that for a reductive group $G$ over a scheme $S$ equipped with a maximal $S$-torus $T \subset G$ and its Lie algebra $\ft \subset \fg$, the map
\[
\kt/W \ra \kg\!\sslash\! G
\]
is an isomorphism whenever $G$ is root-smooth in the sense of \S\ref{pp:root-smoothness} below (for instance, whenever either $2$ is invertible on $S$ or $\sR(G)$ has no contributions of type $C_n$). The statement seems new already when $S$ is an algebraically closed field of small positive characteristic with $G$ nonsemisimple, and it includes results of Springer--Steinberg \cite{SS70}*{Chapter II, Section 3.17$'$} and Chaput--Romagny \cite{CR10}*{Theorem 1.1} as special cases.\footnote{Chaput--Romagny point out in \cite{CR10}*{footnote on page 692} that the Springer--Steinberg proof has an unclear point, which seems to be inherited by several other references that claim this argument in root-smooth settings over a field.} 
The proof seems new already over $\bC$ and uses the Grothendieck alteration $\wt{\fg} \surjects \fg$ reviewed in \S\ref{pp:GS-map}. The idea is to extend the $W$-torsor structure on the regular semisimple locus $\wt{\fg}^\rs$ to a $W$-action on the locus $\wt{\fg}^\fin$ over which the alteration is finite (that is, the regular locus $\wt{\fg}^\reg$, although we do not use the identification $\wt{\fg}^\fin = \wt{\fg}^\reg$ that we later review in \Cref{prop:characterization-of-greg}). To control the $W$-invariants of the extended action, we use a reduction to finite fields trick (see the proof of \Cref{prop:extend-W-action-to-greg}), and then conclude by using the fact the image $\fg^\fin \subset \fg$ of $\wt{\fg}^\fin$ is large enough to cover all the points of $\fg$ that are $S$-fiberwise of height $\le 1$.

\bpp[Root-smoothness] \label{pp:root-smoothness}
A reductive group $G$ over a scheme $S$ is \emph{root-smooth} if for every geometric point $\ov{s}$ of $S$ and every maximal $\ov{s}$-torus $T \subset G_{\ov{s}}$, each root $T \ra \bG_{m,\, \ov{s}}$ is a smooth morphism. For each $s \in S$, it suffices to verify this  for a single $\ov{s}$ over $s$ and a single $T \subset G_{\ov{s}}$: the condition only depends on the $\ov{s}$-isomorphism class of $(G_{\ov{s}}, T)$, and $G(k(\ov{s}))$-conjugation acts transitively on the possible $T$. Thus, by the fibral criterion \cite{EGAIV3}*{Corollaire 11.3.11}, a $G$ with a maximal $S$-torus $T \subset G$ is root-smooth if and only if each root $T \ra \bG_m$ is smooth over the \'{e}tale cover over which it is defined. The smoothness of $T \ra \bG_m$ amounts to the surjectivity of $\Lie(T) \ra \Lie(\bG_m)$, so root-smoothness is an open condition: if $G_s$ is root-smooth for an $s \in S$, then so is $G_U$ for some open neighborhood $U \subset S$ of $s$; in particular, there is the unique largest open of $S$ over which $G$ is root-smooth.

Root-smoothness amounts to a concrete combinatorial condition: indeed, a $T \ra \bG_{m,\, \ov{s}}$ is smooth if and only if it not divisible by $\Char(k(\ov{s}))$ in $X^*(T)$. For instance, by \cite[Section 13.3]{Jan04}, the reductive group $G$ is root-smooth whenever
\bitem
\m
$2$ is invertible on $S$; or
\m
$G_\der$ is adjoint (for instance, whenever $G$ is adjoint semisimple); or  

\m
the root data of the geometric $S$-fibers of $G$ have no contributions of type $C_n$ with $n \ge 1$;
\eitem
so certainly whenever the order of the Weyl group of $G$ is invertible on $S$ (of course, it even suffices that each $S$-fiber of $G$ satisfy one of the above conditions). A basic example of a $G$ that is not root-smooth is $\SL_2$ in characteristic $2$ (note that $C_1 = A_1$). 
\epp

Root-smooth groups admit a Lie-theoretic characterization of maximal tori. We record this in \Cref{eqC} because it generalizes \cite{SGA3II}*{Expos\'{e} XIV, Th\'{e}or\`{e}me 3.18} that exhibited such a description for adjoint semisimple $G$ (see also \cite{SGA3II}*{Expos\'{e} XIV, Corollaire 3.19} where root-smoothness is already~visible).

\bpp[Subgroups of type (C)] \label{pp:type-C-review}
For a reductive group $G$ over a scheme $S$, we recall from \cite{SGA3II}*{Expos\'{e} XIII, Proposition 4.4 and what follows; Expos\'{e} XIV, D\'{e}finition 2.4} that a Lie $S$-subalgebra $\fc \subset \Lie(G)$ is \emph{Cartan} if $\fc$ is Zariski locally on $S$ a module direct summand of $\Lie(G)$ and, for every geometric point $\ov{s}$ of $S$, the Lie subalgebra $\fc_{\ov{s}} \subset \Lie(G)_{\ov{s}}$ is nilpotent and equal to its own normalizer. Every Cartan $\fc$ arises from an $S$-subgroup: in fact, by \cite[Expos\'{e} XIV, Th\'{e}or\`{e}me 3.9]{SGA3II}, the functor $H \mapsto \Lie(H)$ induces a bijection between the closed, smooth, fiberwise connected subgroups $H \subset G$ such that $\Lie(H) \subset \Lie(G)$ is a Cartan subalgebra (such $H$ are the subgroups \emph{of type {\upshape (C)}}) and the Cartan subalgebras of $\Lie(G)$. 
 By \cite{SGA3II}*{Expos\'{e} XIV, Corollaire 3.5, Th\'{e}or\`{e}me 3.9}, any two subgroups of $G$ of type (C) (respectively, any two Cartan subalgebras of $\Lie(G)$) are $G$-conjugate \'{e}tale locally on $S$. If $S = \Spec(k)$ for a field $k$, then, by \cite{SGA3II}*{Expos\'{e} XIV, Lemme~1.2}, the group $G$ has a subgroup of type (C), equivalently, then the Lie algebra $\Lie(G)$ has a Cartan subalgebra.
\epp

\begin{prop}\label{eqC}
A reductive group $G$ over a scheme $S$ is root-smooth if and only if \'{e}tale locally on $S$ it has a subgroup of type {\upshape(C)} that is a maximal torus, in which case the subgroups of type {\upshape(C)} of base changes of $G$ are precisely the maximal tori\uscolon in particular, for a root-smooth $G$, the Cartan subalgebras of $\Lie(G)$ are precisely the Lie algebras of the maximal $S$-tori of $G$. 
\end{prop} 

\begin{proof}
By \cite{SGA3II}*{Expos\'{e} X, Corollaire 4.9}, being a torus is a fibral condition on a closed, smooth, fiberwise connected $S$-subgroup $H \subset G$, so we lose no generality by assuming that $S$ is a geometric point. 
We then let $T \subset G$ be a maximal torus with its Lie algebra $\ft \subset \fg$. Since $\kt$ is abelian, by the criterion \cite[Expos\'{e} XIII, Proposition 4.4 and what follows]{SGA3II}, it is a Cartan subalgebra of $\fg$ (that is, $T$ is a subgroup of type (C)) if and only if there exists a $\gamma\in\kt$ for which $\ad(\gamma)|_{\kg/\kt}$ is injective. However, the root decomposition 
\[
\tst \kg=\kt\oplus\bigoplus_{\gA}\kg_{\gA},
\]
shows that $\ad(\gamma)|_{\kg/\kt}$ is semisimple with eigenvalues $d\gA(\gamma)$, so the injectivity of $\ad(\gamma)|_{\kg/\kt}$ amounts to the nonvanishing of each $d\gA(\gamma)$. Thus, $G$ is root-smooth if and only if $T$ is a subgroup of type (C). To show that in this case every subgroup $H \subset G$ of type (C) is a maximal torus, we first recall from \cite[Expos\'{e} XIV, Th\'{e}or\`{e}me 1.1, Corollaire 5.6 (with Expos\'{e} XII, Section 1.0)]{SGA3II} that any such $H$ contains a maximal torus $T$ of $G$. By \cite[Expos\'{e} XIII, Proposition 4.4 and what follows]{SGA3II}, Cartan subalgebras are maximal nilpotent Lie subalgebras of $\fg$, so the inclusion $\Lie(T) \subset \Lie(H)$ is an equality. It then follows from \S\ref{pp:type-C-review} that $T = H$.
\end{proof}

\brem \label{rem:SL2-fail}
In general, Cartan subalgebras differ from Lie algebras of maximal tori. For instance, we recall from \cite[Expos\'{e} XIII, Remarques 6.6 c)]{SGA3II} that the Lie algebra of $G\ce (\SL_2)_{\ov{\bF}_2}$ is nilpotent, to the effect $\Lie(G)$ is its own unique Cartan subalgebra (and $(\SL_2)_{\ov{\bF}_2}$ is its own unique subgroup of type (C)).
\erem

Root-smooth groups possess a well-behaved regular semisimple locus $\fg^\rs \subset \fg$ that will be used below.


\bpp[The regular semisimple locus $\fg^\rs \subset \fg$] \label{pp:regular-semisimple}
For a reductive group $G$ over a scheme $S$ and a variable $S$-scheme $S'$, we recall from \cite{SGA3II}*{Expos\'{e} XIII, Proposition 4.2 and what follows; Expos\'{e} XIV, D\'{e}finition 2.5} that a section $\gamma \in \fg(S')$ is \emph{regular semisimple} (called simply `regular' in \emph{loc.~cit.}) if it lies in some Cartan subalgebra $\fc \subset \fg_{S'}$ (see \S\ref{pp:type-C-review}) and for every geometric point $\ov{s}$ of $S'$ we have 
\[
\tst \fc_{\ov{s}} = \bigcup_{n \ge 0} \Ker(\ad(\gamma_{\ov{s}})^n).
\]
By \cite{SGA3II}*{Expos\'{e} XIV, Proposition 2.6}, for a regular semisimple $\gamma$ we have 
\[
\tst \fc = \bigcup_{n \ge 0} \Ker(\ad(\gamma)^n),
\]
so the Cartan $\fc$ is uniquely determined. Conversely, if $S$ is a geometric point and $\gamma \in \fg(S)$ lies in at least one and at most finitely many Cartans $\fc \subset \fg$, then, by \cite{SGA3II}*{Expos\'{e} XIII, Th\'{e}or\`{e}me 6.1}, it is regular semisimple, and so lies in a unique $\fc$.  
This immediately reveals pathologies: for instance, by \Cref{rem:SL2-fail}, the entire Lie algebra of $(\SL_2)_{\ov{\bF}_2}$ is regular semisimple. To avoid them, we will only discuss regular semisimple sections for root-smooth reductive $G$.

For example, if $G$ is reductive and root-smooth, then, by \Cref{eqC} and the characterization of being regular semisimple in terms of the injectivity of the adjoint action modulo a Cartan \cite[Expos\'{e} XIII, Proposition 4.6]{SGA3II}, a $\gamma \in \fg(S)$ is regular semisimple if and only if it lies in the Lie algebra $\ft \subset \fg$ of a maximal $S$-torus $T \subset G$ 
and $(d\gA)(\gamma_{\ov{s}}) \neq 0$ for every geometric point $\ov{s}$ of $S$ and every root $\gA \colon T_{\ov{s}} \ra \bG_m$. 
In particular, in the root-smooth case a regular semisimple section is semisimple (and is also regular as we will see in \Cref{lem:description-of-centralizer} and \S\ref{reg-part}, which will also show the converse: sections that are both regular and semisimple are regular semisimple).  



When $G$ is reductive and root-smooth, \cite{SGA3II}*{Expos\'{e} XIV, Proposition 2.9, Corollaire 2.10} (with \Cref{eqC}) gives an $S$-fiberwise dense, stable under the adjoint action of $G$ and under the scaling action of $\bG_m$ open
\[
\kg^{\rs}\subset\kg \q\x{that represents the subfunctor of regular semisimple sections.} 
\]
\epp

For later use, we review the description of centralizers of semisimple and regular semisimple sections. For a description of this sort over an arbitrary base, see \Cref{adjfond2} below.


\blem \label{prop:rs-reg-comparison} \label{lem:description-of-centralizer}
Let $G$ be a reductive group  over an algebraically closed field $k$, let $T \subset G$ be a maximal $k$-torus with  Lie algebra $\ft \subset \fg$ and Weyl group $W \ce N_G(T)/T$, and fix a $\gamma \in \ft(k)$. The group $C_G(\gamma)^0$ is reductive with the root system formed by those $T$-roots of $G$ that vanish on $\gamma$, 
\[\ba
&(C_G(\gamma)^0)(k) \qxq{is generated by} T(k) \qxq{and the $T$-root groups} U_\gA(k) \qxq{with} d\gA(\gamma) = 0,  \\
&(C_G(\gamma))(k) \qxq{is generated by}  (C_G(\gamma)^0)(k) \qxq{and the elements of}  W \qxq{centralizing} \gamma, \qx{and}
\ea\]
$\gamma \in \Lie(\mathrm{Cent}(C_G(\gamma)^0))$. In particular, if $G$ is root-smooth and $\gamma \in \fg^\rs(k) \cap \ft(k)$, then 
\[
C_G(\gamma)^0=T.
\]
\elem

\bpf
Since $\gamma$ is semisimple, $C_G(\gamma)$ is $k$-smooth by \cite{Bor91}*{Chapter III, Section 9.1, Proposition}. Thus, \cite{Ste75}*{Lemma 3.7} gives the claims about $C_G(\gamma)^0$ and the descriptions of $(C_G(\gamma)^0)(k)$ and  $(C_G(\gamma))(k)$. Root decompositions of $\Lie(C_G(\gamma)^0)$ show that $\gamma$ lies in the Lie algebra of every maximal torus of $C_G(\gamma)^0$, so 
\[
\gamma \in \Lie(\mathrm{Cent}(C_G(\gamma)^0))
\]
by \cite{SGA3II}*{Expos\'{e} XII, Proposition 4.10 and what follows}. The last claim follows from the rest and \S\ref{pp:regular-semisimple}.
\epf

We turn to reviewing our main tool for studying the morphism $\ft/W \ra \fg\!\sslash\!G$ in \Cref{adj}.

\bpp[The Grothendieck alteration] \label{GS-res} \label{pp:GS-map}
Let $G$ be a reductive group over a scheme $S$, let $\fg$ be its Lie algebra, let $\sB$ be the smooth, projective $S$-scheme that parametrizes Borel subgroups of $G$ and that was constructed in \cite{SGA3IIInew}*{Expos\'{e} XXII, Corollaire 5.8.3}, and consider the closed $G$-stable $S$-subscheme of $\sB \times_S \fg$ defined as a subfunctor~by
\[
\widetilde{\kg} := \{(B,\gamma)\in \sB \times_S \kg~\vert~\gamma \in \Lie(B) \subset \fg \} \subset \sB \times_S \fg,
\]
in other words, $\wt{\fg}$ is the Lie algebra of the universal Borel subgroup of $G$. The formation of $\widetilde{\kg}$ commutes with base change and the projecting to $\fg$ gives a projective morphism, the \emph{Grothendieck~alteration}
\[
\widetilde{\kg} \surjects \kg,
\]
which is $G$-equivariant and surjective as indicated: indeed, if $S$ is a geometric point, then \cite{SGA3II}*{Expos\'{e} XIV, Th\'{e}or\`{e}me 4.11} ensures that $\fg$ is the union of the Lie algebras of Borel subgroups of $G$.
The projection 
\[
\wt{\fg} \surjects \sB
\]
is also $G$-equivariant and Zariski locally on $\sB$ isomorphic to a relative affine space (namely, to the total space of the Lie algebra of the universal Borel subgroup of $G$). Thus, $\wt{\fg}$ is $S$-smooth and its geometric $S$-fibers are integral, of the same dimension as those of $G$. In particular, by the dimension formula \cite{EGAIV2}*{Corollaire 5.6.6}, the open $\fg^{\fin} \subset \fg$ over which the Grothendieck alteration is finite is $S$-fiberwise dense in $\fg$ (\Cref{prop:characterization-of-greg} reviews a group-theoretic description of $\fg^\fin$).
 In fact,\footnote{If $G$ is root-smooth, then this claim is a special case of a much more general property of the Grothendieck alteration, namely, of its smallness \cite{Jan04}*{Section 12.17; Section 13.2, Lemma}.}
\be \label{eqn:gfin-fiberwise-large}
\fg^{\fin}  \qxq{contains all the points of}  \fg \q \x{that are of height $\le 1$ in their $S$-fiber,}
\ee
indeed, for all $s \in S$, the preimage of $\fg_s\setminus \fg_s^\fin$  is a proper closed subscheme of $\wt{\fg}_s$, so the dimension bound \cite{EGAIV3}*{\'{e}quation 10.6.1.2} gives $\dim(\fg_s) - 1 \ge \dim(\fg_s \setminus \fg_s^\fin) + 1$. Due to the $S$-smoothness of its source and target, the fibral criterion \cite{EGAIV3}*{Corollaire 11.3.11}, and miracle flatness \cite{EGAIV2}*{Proposition 6.1.5},
\be \label{eqn:over-gfin}
\x{the base change} \q \widetilde{\kg}^\fin \surjects \kg^\fin \qxq{of} \widetilde{\kg} \surjects \kg \q \x{is finite locally free.}
\ee
By transport of structure, the open $\fg^\fin$ is $G$-stable.

For any Borel $S$-subgroup $B \subset G$, the quotient $B \surjects T$ by the unipotent radical defined in \cite{SGA3IIInew}*{Expos\'{e} XXII, Proposition 5.6.9} with its Lie algebra map $\mathrm{pr}\colon\fb \surjects \ft$ gives a morphism
\be \label{eqn:gtilde-to-t}
\wt{\fg} \surjects \ft \qxq{defined by} (B', \gamma) \mapsto \mathrm{pr}(\Ad(g)\gamma) 
\ee
where $g$ is an \'{e}tale local section of $G$ such that $g B'g\i = B$: the map is well defined because $T$ is abelian and the self-normalizing property of Borel subgroups \cite{SGA3IIInew}*{Expos\'{e} XXII, Corollaire~5.8.5} ensures that $g$  is unique up to left multiplication by a section of $B$. 
Moreover, if $B$ has an $S$-Levi $T \subset B$, then, since the maps $\ft \ra \fg\!\sslash\!G$ and $\fb \ra \fg\!\sslash\!G$ induced by inclusion are compatible with $\mathrm{pr} \colon \fb \surjects \ft$ (as is seen from \cite{AFV18}*{Step 2 of the proof of Lemma 13}), we have a commutative diagram
\be \ba \label{eqn:cool-diagram}
\xymatrix@R=20pt@C=10pt{
\wt{\fg} \ar@{->>}[rr]^-{\eqref{eqn:gtilde-to-t}} \ar@{->>}[d] && \ft \ar[d] \\
\fg \ar[r] & \fg\!\sslash\!G  & \ft/W. \ar[l]
}
\ea\ee
\epp

We are ready to describe a $W$-action on the preimage $\wt{\fg}^\rs$ of $\fg^\rs$ and then to extend it to $\wt{\fg}^\fin$.

\begin{prop} \label{prop:describe-restriction-to-rs}
Let $G$ be a root-smooth reductive group over a scheme $S$, let $T \subset G$ be a maximal $S$-torus with Lie algebra $\ft \subset \fg$, set $\ft^\rs \ce \ft \cap \fg^\rs $, let $W \ce N_G(T)/T$ be the Weyl group, and let $T \subset B \subset G$ be a Borel $S$-subgroup. The restriction 
\[
\wt{\fg}^\rs \surjects \fg^\rs
\]
of the Grothendieck alteration has the structure of a $W$-torsor that is compatible with the $G$-actions on $\wt{\fg}^\rs$ and $\fg^\rs$ and is supplied by the $W$-action on $G/T\times_S\kt$ over $\fg$ given by 
\[
w\cdot (gT,\gamma)\mapsto (gw^{-1}T,\, \Ad(w)\gamma)
\]
and by the commutative diagram
\be\label{eqn:grs-diagram} \ba
\xymatrix@C=45pt{
G/T\times_S\kt^{\rs} \ar@{->>}[rd]_<<<<<<<<{(gT,\,\gamma)\,\mapsto\,\Ad(g)\gamma\ \ } \ar[rr]_-{\sim}^-{(gT,\, \gamma)\, \mapsto\, (gBg\i,\, \Ad(g)\gamma)} && \wt{\fg}^\rs \ar@{->>}[ld] \\
& \fg^\rs.
} 
\ea\ee
In terms of this diagram the map $\wt{\fg}^\rs \xra{\eqref{eqn:gtilde-to-t}} \ft$ corresponds to the projection onto $\ft^\rs$, and so is $W$-equivariant. In particular, for root-smooth $G$, we have
\[
\fg^\rs \subset \fg^\fin.
\]
\end{prop}

\begin{proof}
We may focus on the claim about \eqref{eqn:grs-diagram} because it implies the other assertions.
A $W$-action makes $\wt{\fg}^\rs$ a $W$-torsor over $\fg^\rs$ if and only if the maps 
\[
W_{\fg^\rs} \times_{\fg^\rs} \wt{\fg}^\rs \xra{(w,\, x)\, \mapsto\, (wx,\, x)} \wt{\fg}^\rs \times_{\fg^\rs} \wt{\fg}^\rs
\]
and $\wt{\fg}^\rs/W \ra \fg^\rs$ are isomorphisms. Thus, thanks to the fibral criterion \cite{EGAIV4}*{Corollaire 17.9.5}, we may pass to $S$-fibers and assume that $S$ is a geometric point. Moreover, 
\[
\xymatrix@C=20pt{G/T\times_S\kt \ar[rrrr]_-{\sim}^-{(gT,\, \gamma)\, \mapsto\, (gT,\, \Ad(g)\gamma) } &&&& \wt{\fg}_T \ce \{ (gT, \gamma)\, \vert\, \gamma \in \Ad(g)(\ft) \} \subset G/T \times_S \fg,}
\]
and in terms of the target of this isomorphism the map to $\wt{\fg} \subset \sB \times_S \fg \cong  G/B \times_S \fg$ becomes $(gT, \gamma) \mapsto (gB, \gamma)$. Thus, \cite[Section 13.4, Lemma]{Jan04} shows that the top horizontal arrow in \eqref{eqn:grs-diagram} is an isomorphism and the vertical ones are finite \'{e}tale. It then remains to show that $W$ acts simply transitively on the $S$-fibers of $\wt{\fg}^\rs \surjects \fg^{\rs}$. In terms of $\wt{\fg}_T$, the $W$-action is 
\[
(w, (gT, \gamma))\mapsto (gw\i T, \gamma),
\]
and hence commutes with the evident $G$-action, for which the projection $\wt{\fg}_T \ra \fg$ is $G$-equivariant. By \S\ref{pp:regular-semisimple}, the $G(S)$-translates of $\ft^\rs(S)$ exhaust $\fg^\rs(S)$, so we only need to consider the $S$-fibers above $\ft^\rs \subset \fg^\rs$. For the latter, we fix a $\gamma \in \ft^\rs(S)$ and apply the last assertion of \Cref{lem:description-of-centralizer} to get $C_G(\gamma)^0 = gTg\i$ for every $S$-point $(gT, \gamma)$ of $\wt{\fg}_T$ above $\gamma$. This means that that the $S$-fiber above $\gamma$ is precisely $\{(wT, \gamma) \, \vert \, w \in W(S)\}$, and the desired simple transitivity follows.
\end{proof}


\begin{prop} \label{prop:extend-W-action-to-greg}
Let $S$, $G$, $T \subset B \subset G$, $\ft \subset \fg$, and $W$ be as in Proposition \uref{prop:describe-restriction-to-rs}. The $W$-action on $\wt{\fg}^\rs$ over $\fg^\rs$ constructed there extends uniquely to a $W$-action on $\wt{\fg}^\fin$ over $\fg^\fin$ that commutes with the $G$-actions on $\wt{\fg}^\fin$ and $\fg^\fin$, this extension satisfies 
\[
\wt{\fg}^\fin / W \isomto \fg^\fin, \qxq{and the map} \wt{\fg}^\fin \xra{\eqref{eqn:gtilde-to-t}} \ft \qx{is $W$-equivariant.}
\]
\end{prop}

\begin{proof}
We may work \'{e}tale locally on $S$, so we lose no generality by assuming that $G$ is split. Moreover, since $\wt{\fg}$ is $S$-smooth with integral geometric fibers (see \S\ref{pp:GS-map}) and $\wt{\fg}^\rs \subset \wt{\fg}$ is $S$-fiberwise dense (see \S\ref{pp:regular-semisimple}), the uniqueness aspect follows from the relative characterization of support \cite{EGAIV4}*{Proposition 19.9.8}. Similarly, the compatibility with the $G$-actions and the $W$-equivariance of $\wt{\fg}^\fin \ra \ft$ will follow from \Cref{prop:describe-restriction-to-rs}. For the existence and the identification $\wt{\fg}^\fin/W$, granted that we show that the formation of the coarse space $\wt{\fg}^\fin / W$ commutes with any base change, we reduce to $S$ being open in $\Spec(\bZ)$ (passage to an open ensures root-smoothness, see \S\ref{pp:root-smoothness}).

In the case when $S$ is a localization of $\Spec(\bZ)$, the scheme $\wt{\fg}^\fin$ (and also $\wt{\fg}$) is normal, so, by \eqref{eqn:over-gfin}, it is the normalization of $\fg^\fin$ in $\wt{\fg}^\rs$. The $W$-action then extends by the functoriality of normalization. By the normality of rings of invariants \cite{AM69}*{Proposition 7.8}, the coarse quotient $\wt{\fg}^\fin/W$ inherits normality and finiteness over $\fg^\fin$ from $\wt{\fg}^\fin$. In addition, by \Cref{prop:describe-restriction-to-rs}, the map $\wt{\fg}^\fin/W \ra \fg^\fin$ is an isomorphism over $\fg^\rs$. Since the latter is dense in $\fg^\fin$ (see \S\ref{pp:regular-semisimple}), it follows that this map is the normalization morphism for $\fg^\fin$ in its function field, so, by normality, it must be an isomorphism. The same argument works if the base is a finite field instead of a localization of $\bZ$, so the formation of $\wt{\fg}^\fin/W$ commutes with base change to every finite field. It then follows from the $\bZ$-fibral criterion \cite{modular-description}*{Lemma 3.3.1} that its formation commutes with any base change, as promised.
\end{proof}


We turn to the Chevalley isomorphism for root-smooth groups. Even though root-smoothness is a very mild condition, it is  
not always necessary: for instance, by \cite[Theorem 1.2]{CR10}, for $G = \Sp_{2n}$ (type $C_n$) the Chevalley map is an isomorphism if and only if $S$ has no nonzero $2$-torsion, but such $G$ is root-smooth if and only if $2$ is invertible on $S$. 




\begin{theorem}\label{adj} \label{chev}
For a reductive group $G$ over a scheme $S$, the adjoint action of $G$ on its Lie algebra $\fg$, and a maximal $S$-torus $T \subset G$ with Lie algebra $\ft \subset \fg$ and Weyl group $W \ce N_G(T)/T$, 
\[
\kt/W \ra \kg\!\sslash\! G \ \ \x{is a schematically dominant map that is an isomorphism if $G$ is root-smooth.}
\]
\end{theorem}

\begin{proof}
The schematic dominance was settled in \cite{AFV18}*{Remark 14} (and in \cite{CR10}*{Theorem 3.6} in a special case), so we assume that $G$ is root-smooth and seek to show the isomorphism.  The formation of $\kt/W \ra \kg\!\sslash\! G$ commutes with flat base change (see \S\ref{conv}), so we work \'{e}tale locally on $S$ to assume that $S = \Spec(A)$ is affine and $G$ is split with respect to $T$, equipped with a Borel $T \subset B \subset G$. 
Our task is to show that the injection $A[\kg]^G \hra A[\ft]^W$ is surjective. The idea is to consider the diagram
\[
\xymatrix@C=20pt{
G/T \times_S \ft^\rs \ar@{->>}[d] \ar@{^(->}[r] & \wt{\fg}^\fin \ar@{->>}[d]\ar[rr]^{\eqref{eqn:gtilde-to-t}} & & \ft \\
\fg^\rs \ar@{^(->}[r] & \fg^\fin
}
\]
supplied by \Cref{prop:describe-restriction-to-rs,prop:extend-W-action-to-greg}. The square is Cartesian and the top horizontal maps compose to a projection and are $W$-equivariant. Pullback of an $\gA \in A[\ft]^W$ along this projection is a $G$-invariant and $W$-invariant global section $f$ of $G/T \times_S \ft^\rs$. Since the left vertical map is a $W$-torsor, the $W$-invariance means that $f$ comes from a $G$-invariant global section of $\fg^\rs$ whose restriction to $\ft^\rs$, by construction, agrees with the restriction of $\gA$. Thus, by the relative characterization of support \cite{EGAIV4}*{Proposition 19.9.8}, all that remains is to extend $f$ to a global section of $\fg$, which will then necessarily be unique and $G$-invariant. In fact, since the open 
\[
\fg^\fin \subset \fg
\]
contains all the points of $\fg$ that have height $\le 1$ in their $S$-fiber (see \eqref{eqn:gfin-fiberwise-large}), \emph{loc.~cit.}~even ensures that it suffices to extend $f$ to a global section of $\fg^\fin$. For this, the pullback of $\gA$ along the second top horizontal map extends $f$ to a $W$-invariant global section of $\wt{\fg}^\fin$, which, by \Cref{prop:extend-W-action-to-greg}, descends to a desired extension of $f$ to $\fg^\fin$.
\end{proof}

\brem
As we already mentioned, in general the map $\kt/W \ra \kg\!\sslash\! G$ is not an isomorphism: for instance, by \cite[Section 6.1]{CR10}, for $G = (\SL_2)_{\bF_2}$ equipped with its diagonal torus,
\[ (\bF_2[\fg])^{G} \hra (\bF_2[\kt])^{W} \qxq{may be identified with the inclusion}  \bF_2[x^2]  \hra \bF_2[x].
\]
In contrast, the group version of \Cref{adj} does hold in general: for any reductive $G$ over any base scheme $S$ and any maximal $S$-torus $T \subset G$, by \cite{Lee15}*{Section 1, Theorem},
\[
 T/ W \isomto  G \!\sslash\! G, \q\x{where $G$ acts on itself by conjugation.}
\]
\erem

The formation of the adjoint quotient $\fg\sslash G$ need not commute with nonflat base change, see \cite{CR10}*{Theorem 1.3} for (root-smooth) counterexamples with $G$ of type $B_n$ or $D_n$ in characteristic $2$. Nevertheless, \Cref{adj} relates this base change to its analogue for the \emph{a priori} simpler quotient $\ft/W$. For instance, it implies that for root-smooth $G$, if the order of $W$ is invertible on $S$, then the formation of $\fg\!\sslash\! G$ commutes with arbitrary base change. In fact, by \Cref{prop:base-change-for-g//G} below, a significantly weaker condition based on the following notion of torsion primes suffices for this.

\bpp[Torsion primes for a root datum] \label{pp:torsion-primes}
Fix a reduced root datum  
\[
\sR = (X, \Phi, X^\vee, \Phi^\vee)
\]
and consider its associated semisimple, simply-connected root datum
\[
\sR^{\mathrm{sc}} = (X^{\mathrm{sc}}, \Phi, (X^{\mathrm{sc}})^\vee, \Phi^\vee),
\]
so that $(X^{\mathrm{sc}})^\vee = \bZ\Phi^\vee \subset X^\vee$. By \cite{SGA3IIInew}*{Expos\'{e} XXI, Corollaires 7.1.6 et 7.4.4, Remarque~7.4.6}, the datum
 $\sR^{\mathrm{sc}}$ is a product of semisimple, simply-connected root data of one of the well-known Dynkin types: $A_n$ with $n \ge 1$, or $B_n$ with $n \ge 2$, or $C_n$ with $n \ge 3$, or $D_n$ with $n \ge 4$, or $E_6$, or $E_7$, or $E_8$, or $F_4$, or $G_2$. Following \cite{Dem73}*{Proposition 6; Proposition 8 and what follows}, we say that a prime $p$ is a \emph{torsion prime} for $\sR$ if either
\benuma
\item \label{item:torsion-prime-type-1}
$p \mid \# \Coker(X \ra X^{\mathrm{sc}})$ (for an $\sR$ associated to a split reductive group $G$, by \cite{SGA3IIInew}*{Expos\'{e} XXII, Propositions 4.3.1 et 6.2.7; Expos\'{e} XXIII, Th\'{e}or\`{e}me 4.1}, this cardinality is the degree of the isogeny $(G_\der)^{\mathrm{sc}} \ra G_\der$);~or

\item \label{item:torsion-prime-type-2}
$p = 2$ and $\sR^{\mathrm{sc}}$ has a factor of one of the following types: $B_n$, $D_n$, $E_6$, $E_7$, $E_8$, $F_4$, $G_2$;

\item \label{item:torsion-prime-type-3}
$p = 3$ and $\sR^{\mathrm{sc}}$ has a factor of one of the following types: $E_6$, $E_7$, $E_8$, $F_4$;

\item \label{item:torsion-prime-type-4}
$p = 5$ and $\sR^{\mathrm{sc}}$ has a factor of type  $E_8$.
\eenum
By \cite{Dem73}*{Lemme 7}, torsion primes for $\sR$ divide the order of the Weyl group of $\sR$. The converse fails: for instance, a semisimple, simply connected $\sR$ of type $C_n$ has no torsion primes. Being a torsion prime for a root datum (as above) is different than being a torsion prime for the associated root \emph{system} in the sense sometimes used in the literature (compare with \cref{foot:torsion-for-root-system} below).
\epp

\bd \label{def:no-torsion}
For a reductive group $G$ over a scheme $S$, we say that \emph{$\sR(G)$ has no torsion residue characteristics} if $\Char(k(s))$ is not a torsion prime for $\sR(G_{\ov{s}})$ for every $s \in S$. Conversely, we say that a prime $p$ is a \emph{torsion residue characteristic for $\sR(G)$} if there is an $s \in S$ with $p = \Char(k(s))$ such that $p$ is a torsion prime for $\sR(G_{\ov{s}})$.
\ed

\begin{prop} \label{prop:base-change-for-g//G}
For a reductive group $G$ over a scheme $S$ with Lie algebra $\fg$, if $\sR(G)$ has no torsion residue characteristics, 
then $\fg\!\sslash\! G$ is of formation compatible with base change and \'{e}tale locally on $S$ an affine space of relative dimension $\rk(G)$.
\end{prop}

We will refine the compatibility with base change aspect of \Cref{prop:base-change-for-g//G} in \eqref{eqn:greg-gerbe} below.

\begin{proof}
We work \'{e}tale locally to assume that $G$ is split. The claims then follow from \cite{AFV18}*{Proposition 10 (with Remark 8)}, which provides an isomorphism between $\fg\!\sslash\! G$ and a Kostant section $\cS \subset \fg$ that is an affine space of relative dimension $\rk(G)$ and commutes with base change (see also \S\ref{pp:Kostant-section}). 
\end{proof}

\brem
The combination of \Cref{adj} and \Cref{prop:base-change-for-g//G} implies that for a root-smooth reductive group $G$ over a scheme $S$ such that $\sR(G)$ has no torsion residue characteristics and a maximal $S$-torus $T \subset G$ with Lie algebra $\ft \subset \fg$ and Weyl group $W \ce N_G(T)/T$, the formation of $\ft/W$ commutes with arbitrary base change---this reproves \cite{Dem73}*{Corollaire sur la page 296}.
\erem



\bq \label{q:univ-homeo}
For a reductive group $G$ over a scheme $S$, its Lie algebra $\fg$, and an $S$-scheme $S'$, is the map $\kg_{S'}\!\sslash\! G_{S'} \ra (\kg\!\sslash\! G)_{S'}$ always a universal homeomorphism?
\eq

For $\ft/W$, the positive answer is a general property of coarse moduli spaces, so \Cref{adj} gives a positive answer whenever $G$ is root-smooth. \Cref{prop:base-change-for-g//G} does the same whenever $\sR(G)$ has no torsion residue characteristics; in fact, this also follows from the following variant of \Cref{adj}.

\bcor \label{cor:always-univ-homeo}
For a reductive group $G$ over a scheme $S$ and a maximal $S$-torus $T \subset G$ with Lie algebra $\ft \subset \fg$ and $W \ce N_G(T)/T$, if $\sR(G)$ has no torsion residue characteristics, then the schematically-dominant morphism $\kt/W \ra \kg\!\sslash\! G$ is a universal homeomorphism.
\ecor

\bpf
We use the fpqc local on the base nature of being a universal homeomorphism \cite{SP}*{Lemma~\href{https://stacks.math.columbia.edu/tag/0CEX}{0CEX}} to work \'{e}tale locally on $S$ and assume that $G$ is split with respect to $T$. 
By \Cref{prop:base-change-for-g//G}, the formation of $\kg\!\sslash\! G$ commutes with any base change. For the coarse space $\ft/W$, the same holds up to a universal homeomorphism, so we may assume that $S = \Spec(\bZ)$. Then, since $T$ is split, 
\[
T \cong \underline{\Hom}(X^*(T), \bG_m), \qxq{so} \Lie(T) \cong (X^*(T))^\vee, \qxq{and hence} \ft \cong \Spec(\Sym(X^*(T)))
\]
(see \S\ref{conv}). By \cite{Dem73}*{Th\'{e}or\`{e}me 3}, the assumption on the residue characteristics implies that the $W$-invariants of $\Sym(X^*(T))$ form a polynomial algebra, so $\ft/W$ is an affine space of relative dimension $\rk(G)$. 
Moreover, by \cite{Jan04}*{Section 7.13, Claim}, the map $\kt/W \ra \kg\!\sslash\! G$ is bijective on points valued in every algebraically closed $S$-field. 
Consequently, since $\kg\!\sslash\! G$ is also an affine space of relative dimension $\rk(G)$ (see \Cref{prop:base-change-for-g//G}), the fibral criterion \cite{EGAIV3}*{Corollaire 11.3.11} and \cite{EGAIV2}*{Proposition 6.1.5} ensure that $\kt/W \ra \kg\!\sslash\! G$ is flat. 
Since this map is also of finite presentation, it is an open, continuous bijection on topological spaces, and remains so after any base change, so it is a universal homeomorphism. 
\epf

The preceding corollary allows us to describe the basic geometric properties of the map $\fg \ra \kg\!\sslash\!G$.

\bcor\label{prop:describe-Chevalley}
Let $G$ be a reductive group over a scheme $S$ such that either $G$ is root-smooth, or $2$ is not a torsion residue characteristic for $\sR(G)$, or $S$ is the spectrum of a field. 
The affine map
\be \label{eqn:topological-properties-of-Chevalley}
\fg \ra \kg\!\sslash\!G \ \ \x{is surjective, with irreducible geometric fibers of dimension $\dim(G) - \rk(G)$.} 
\ee
These geometric fibers consist \up{set-theoretically} of finitely many $G$-orbits with exactly one semisimple orbit, those above the points in $0 \in (\kg\!\sslash\!G)(S)$ consist precisely of the nilpotent sections of $\fg$, and if $G$ is root-smooth, then the geometric fibers of $\fg \ra \kg\!\sslash\!G$ that meet $\fg^\rs$ consist of a single $G$-orbit.
\ecor

\bpf
We may work \'{e}tale locally on $S$, so we fix a maximal $S$-torus $T \subset G$. Moreover, since $G$ becomes root-smooth after inverting $2$ on $S$, we may weaken the case `$2$ is not a torsion residue characteristic for $\sR(G)$' to `$\sR(G)$ has no torsion residue characteristics.' Then, by \Cref{adj} and \Cref{cor:always-univ-homeo}, the assumed conditions ensure that up to a universal homeomorphism the formation of $\kg\!\sslash\!G$ commutes with base change. Thus, \eqref{eqn:topological-properties-of-Chevalley} and all the subsequent claims except for the one about $\fg^\rs$ follow from the geometric point case supplied by \cite{Jan04}*{Section 7.13}. 
For the assertion about $\fg^\rs$, suppose that $G$ is root-smooth, let $\ov{s}$ be a geometric $S$-point, and let $\gamma \in \fg^\rs(\ov{s})$. By \cite{Bor91}*{Chapter III, Section 9.2, Theorem},~the adjoint orbit of $\gamma$ is closed in $\fg$ and, by \Cref{lem:description-of-centralizer} its dimension is $\dim(G) - \rk(G)$. Thus, due to \eqref{eqn:topological-properties-of-Chevalley}, this orbit sweeps out the entire $\ov{s}$-fiber of $\fg \ra \kg\!\sslash\!G$.
\epf

\brem
If the answer to \Cref{q:univ-homeo} is positive, then, except for the claim about $\fg^\rs$, \Cref{prop:describe-Chevalley} holds without any assumptions on the reductive group $G$ over a scheme $S$. 
\erem


\csub[The geometry of the Chevalley morphism on the regular locus] \label{section:Chevalley-Kostant}


As we discuss in this section, the Chevalley map $\kg\rightarrow\kg\!\sslash\!G$ is particularly well behaved when restricted to the regular locus $\fg^\reg \subset \fg$ (equivalently, to $\fg^\fin$, see \Cref{prop:characterization-of-greg}). 
The main point is that the results below hold under weaker assumptions than known previously and over an arbitrary base---roughly, it suffices to assume that the torsion primes for the root datum of $G$ are invertible on $S$. Under this assumption, we show that the map $\fg^\reg \ra \fg\!\sslash\!G$ is smooth (see \Cref{thm:describe-chi}), construct a canonical descent $J$ to $\kg\!\sslash\!G$ of the centralizer of the universal regular section of $\fg$ (see \Cref{pla}), and review a Galois-theoretic description of $J|_{\kg^\rs\sslash G}$ that will be crucial for \S\ref{section:Hitchin-fibration} (see \Cref{prop:galois-theoretic}). This generalizes and improves various statements in the literature, notably from \cite{Ngo10} and \cite{Ric17}, builds the setup for studying the Hitchin fibration in \S\ref{section:Hitchin-fibration}, and leads to a concrete result about the conjugacy class of the Kostant section that is presented in \Cref{conj}.  


We begin by recalling the definition of the regular locus $\fg^\reg \subset \fg$ and analyzing its nilpotent sections.


\bpp[The regular locus $\fg^\reg \subset \fg$] \label{reg-part} 
For a reductive group $G$ over a scheme $S$, its Lie algebra $\fg$, and an $S$-scheme $S'$, the centralizer $C_G(\gamma) \subset G_{S'}$ of a section $\gamma \in \fg(S')$ under the adjoint action of $G$ on $\fg$ is a closed $S'$-subgroup of $G_{S'}$ whose formation commutes with base change. We have\footnote{\label{foot-centralizer}For the sake of completeness, we recall the argument. Letting $\ov{k}$ be the algebraic closure of the residue field at $s$, we choose a Borel subgroup $B \subset G_{\ov{k}}$ with $\gamma \in \Lie(B)$ (see \S\ref{GS-res}), and we let $U \subset B$ be the unipotent radical. Since $B/U$ is commutative, the $B$-orbit of $\gamma_{\ov{k}}$ has a constant image in $\Lie(B/U)$, so the dimension of this orbit is $\le \dim(U)$, and hence $\dim(C_B(\gamma_{\ov{k}})) \ge \rk(G_s)$. Consequently, since $C_B(\gamma_{\ov{k}}) \subset C_G(\gamma)_{\ov{k}}$, we obtain the desired $\dim(C_G(\gamma)_s) \ge \rk(G_s)$.}
\[
\dim(C_G(\gamma)_s) \ge \rk(G_s) \qxq{for every} s \in S',
\]
and we say that $\gamma$ is \emph{regular} 
if the equality holds:
\be \label{reg-ineq}
\dim(C_G(\gamma)_s) = \rk(G_s) \qxq{for every} s \in S'. 
\ee
The function $s \mapsto \rk(G_s)$ is locally constant on $S'$, so the Chevalley semicontinuity theorem \cite[Th\'{e}or\`{e}me 13.1.3]{EGAIV3} (applied along the identity section of $C_G(\gamma)$) ensures that the $s \in S'$ at which \eqref{reg-ineq} holds form an open subscheme of $S'$ whose formation commutes with base change. By considering the universal case $S' = \fg$, we find an open subscheme
\[
\fg^\reg \subset \fg \q\x{that represents the subfunctor of regular sections.}
\]
By transfer of structure, the adjoint action of $G$ on $\fg$ preserves $\fg^\reg$, and so does the scaling action of $\bG_m$ on $\fg$. \Cref{lem:description-of-centralizer} implies that
\[
\fg^\rs \subset \fg^\reg\qx{for root-smooth $G$}
\]
and gives a converse: if each geometric fiber of $\fg^\reg$ has semisimple sections, then $G$ is root-smooth.
\epp


\bpp[Regular nilpotent elements] \label{pp:regular-nilpotent}
The regular locus $\fg^\reg$ is $S$-fiberwise dense in $\fg$: indeed, if $G$ has semisimple rank $\ge 1$ and $S = \Spec(k)$ is a geometric point, then $\fg^\reg(k)$ has nilpotent elements, see \cite{Jan04}*{Sections 6.3--6.4} (for $k$ of low characteristic, this relies on the type-by-type analysis of nilpotent $G$-orbits in $\fg$). Moreover, regular nilpotent elements of $\fg(k)$  form a single $G(k)$-conjugacy class (\emph{loc.~cit.}). We may use \cite{Spr66}*{proof of Lemma 5.8} to describe such elements explicitly: granted that we choose a splitting of $G$, then its pinning $\{e_\gA \in \fg_\gA\}_{\gA \in \Delta}$, and then extend to a system of bases $\{e_\gA \in \fg_\gA\}_{\gA > 0}$ (for instance, to the positive part of a Chevalley system whose existence is ensured by \cite{SGA3IIInew}*{Expos\'{e} XXIII, Proposition 6.2}), 
a nilpotent element of $\fg(k)$ is regular if and only if it is $G(k)$-conjugate to some
\[
\tst \sum_{\al\ssup 0}\xi_{\al}e_{\al} \qxq{with} \xi_\alpha \neq 0 \qxq{for} \alpha \in \Delta.
\]
Thus, by \cite{Spr66}*{proof of Lemma 5.3}, every regular nilpotent element of $\fg(k)$ lies in the Lie algebra of a unique Borel subgroup of $G$.  Conversely, a nilpotent element $\gamma \in \fg(k)$ that is not regular lies in the Lie algebras of infinitely many Borel subgroups: indeed, by conjugation we may assume that $x = \sum_{\al\ssup 0}\xi_{\al}e_{\al}$ with $\xi_{\beta} = 0$ for some $\beta \in \Delta$, so that 
\[
x \in \Lie(R_u(P_{\beta})) \subset \Lie(B_\Delta)
\]
for the Borel $B_\Delta \subset G$ associated to $\Delta$ and the minimal parabolic  $B_\Delta \subset P_{\beta} \subset G$ whose Lie algebra contains $\fg_{-\beta}$ (for the construction of $P_\beta$ and a characterization of its unipotent radical, see \cite{SGA3IIInew}*{Expos\'{e} XXVI, Propositions 1.4 et 1.12~(i)}); since 
\[
\Ad(g)x \in \Lie(R_u(P_{\beta})) \subset \Lie(B_\Delta) \qxq{for} g \in P_{\beta}(k),
\]
the Lie algebras of the infinitely many Borels $\{g\i B_\Delta g\}_{g \in P_{\beta}(k)/B_\Delta(k)}$ all contain $x$. 
\epp


We are ready to relate the regular locus $\fg^\reg$ to the Grothendieck alteration $\wt{\fg} \ra \fg$ reviewed in \S\ref{pp:GS-map}.

\bprop \label{prop:characterization-of-greg}
For a reductive group $G$ over a scheme $S$ and its Lie algebra $\fg$, the locus $\fg^\fin$ over which the Grothendieck alteration $\wt{\fg} \surjects \fg$ has finite fibers is precisely the regular locus $\fg^\reg \subset \fg$, that is,
\[
\fg^\fin = \fg^\reg \qxq{inside} \fg.
\]
\eprop


\begin{proof}
By passing to fibers, we assume that $S$ is the spectrum of an algebraically closed field $k$, and we need to show that a $\gamma \in \fg(k)$ is regular if and only if there are only finitely many Borel $k$-subgroups of $G$ whose Lie algebras contain $\gamma$. For nilpotent $\gamma$, we saw this in \S\ref{pp:regular-nilpotent}. In general, \cite{Bor91}*{Chapter I, Section 4.4, Theorem} supplies the Jordan decomposition $\gamma = \gamma_s + \gamma_n$ with $\gamma_s$ semisimple, $\gamma_n$ nilpotent, and $[\gamma_s, \gamma_n] = 0$, so that 
\[
\gamma_n \in \Lie(C_G(\gamma_s)).
\]
The functoriality of the decomposition under isomorphisms implies that a $g\in G(k)$ centralizes $\gamma$ if and only if it centralizes both $\gamma_s$ and $\gamma_n$, so
\[
C_G(\gamma) \qxq{and} C_{C_G(\gamma_s)}(\gamma_n) \qxq{agree set-theoretically in} G.
\]
Thus, since, by \Cref{prop:rs-reg-comparison}, the group $(C_G(\gamma_s))^0$ is reductive, of the same rank as $G$,
\be \label{eqn:unique-Borel}
\gamma \in \fg^\reg(k) \qxq{if and only if}  \gamma_n \in (	\Lie(C_G(\gamma_s)^0))^\reg(k),
\ee
equivalently, by \S\ref{pp:regular-nilpotent}, 
\be\label{eqn:unique-Borel-re}
\gamma \in \fg^\reg(k) \qxq{if and only if}  \gamma_n \in(\Lie(C_G(\gamma_s)^0))^\fin(k).
\ee
On the other hand, by the uniqueness of the Jordan decomposition, the Lie algebra of a Borel $B \subset G$ contains $\gamma$ if and only if it contains $\gamma_s$ and $\gamma_n$. In this case, $ C_G(\gamma_s) \cap B = C_B(\gamma_s)$, so \cite{Bor91}*{Chapter~III, Section 9.1, Proposition} ensures that this intersection is smooth. Moreover, $\gamma_s$ lies in the Lie algebra of a maximal torus of $B$, so \Cref{lem:description-of-centralizer} ensures that the smooth, solvable group $( C_G(\gamma_s) \cap B)^0$ is of sufficiently large dimension to be a Borel of $C_G(\gamma_s)^0$. Since $\gamma_n \in \Lie(C_G(\gamma_s)\cap B)$, we find that if $\gamma \in \fg^\reg(k)$, then $B$ contains a fixed Borel of $C_G(\gamma_s)$ (see \eqref{eqn:unique-Borel} and \S\ref{pp:regular-nilpotent}), so also a fixed maximal torus of $C_G(\gamma_s)$. The latter is also maximal for $G$, so $B$ belongs to a finite list of Borels, and hence $\gamma \in \fg^\fin(k)$.


Conversely, by its solvability, every Borel of $C_G(\gamma_s)^0$ lies in a Borel of $G$. Since $\gamma_s$ lies in its Lie algebra (see \Cref{lem:description-of-centralizer}), this Lie algebra contains $\gamma_n$ if and only if it contains $\gamma$. Thus, if $\gamma \in \fg^\fin(k)$, then $\gamma_n \in (\Lie(C_G(\gamma_s)^0))^\fin(k)$, and so, by \eqref{eqn:unique-Borel-re}, also $\gamma \in \fg^\reg(k)$, as desired.
\end{proof}

The promised analysis of the map $\fg^\reg \ra \fg\!\sslash\!G$ will rest on the construction of its Kostant section.

\bpp[A Kostant section] \label{pp:Kostant-section}
Let $(\sR, \Delta)$ be a based root datum, as defined in \cite{SGA3IIInew}*{Expos\'{e}~XXIII, Section 1.5}, let $\tau_\sR$ be the product of the torsion primes for $\sR$ (see \S\ref{pp:torsion-primes}), let $G$ be a split, pinned reductive group over $\bZ[\f{1}{\tau_\sR}]$ associated to $(\sR, \Delta)$, let $\{e_\gA \in \fg_\gA\}_{\gA \in \Delta}$ be the basis given by the pinning, set 
\[
\tst e \ce \sum_{\gA \in \Delta} e_\gA,
\]
let $B \subset G$ be the Borel \emph{opposite} to the one associated to $\Delta$, and let $U \subset B \subset G$ be its unipotent radical with Lie algebra $\fu \subset \fb \subset \fg$. By \cite{Ric17}*{proof of Lemma 3.1.2},\footnote{\label{footnote:Riche-argument}The statement of \cite{Ric17}*{Lemma 3.1.2} is  weaker---there one inverts $2$ for type $C_n$ and $3$ for type $G_2$, whereas we do not---but its argument still works as follows (this improvement was also observed in \cite{AFV18}*{Section 2.3, especially, Remark 8}). By the definition of $\tau_\sR$, especially, by \S\ref{pp:torsion-primes}~\ref{item:torsion-prime-type-1}, we may replace $G$ by $(G_{\der})^{\mathrm{sc}}$ (to keep track of Borel subgroups under this one uses \cite{SGA3IIInew}*{Expos\'{e} XXVI, Proposition 1.19}). Then we consider the grading $\fg = \bigoplus_{i \in \bZ} \fg^i$ by the heights of roots (or by the half of the cocharacter given by the sum of positive coroots, as is equivalent \cite{SGA3IIInew}*{Expos\'{e} XXI, Proposition 3.5.1}), so that
\[
\tst \fu = \bigoplus_{i < 0} \fg^i \qxq{and} \fb = \bigoplus_{i \le 0} \fg^i.
\]
The formulas for the Lie bracket given in \cite{SGA3IIInew}*{Expos\'{e} XXII, Remarque 5.4.10} ensure that the map $[ e, -]$ restricts to a map $\fg^i \ra \fg^{i + 1}$. The desired injectivity and freeness of the cokernel of the latter for $i < 0$ follow from \cite{Spr66}*{Propositions 2.2 and 2.4, Theorem 2.6}, which classifies the possible torsion in the cokernel: indeed, in \cite{Spr66}*{Theorem 2.6}, for type $C_n$ with $p = 2$ and type $G_2$ with $p = 3$ the elementary divisors only occur for $i > 0$, whereas the definition of $\tau_\sR$, especially, \S\ref{pp:torsion-primes}~\ref{item:torsion-prime-type-2}--\ref{item:torsion-prime-type-4}, rules out the other occurrences.} the map
\[
\tst [e, -] \colon \fu \ra \fb \qxq{is injective and its cokernel} \fb/[e, \fu] \qx{is $\bZ[\f{1}{\tau_\sR}]$-free.}
\]
We may therefore find a $\bZ[\f{1}{\tau_\sR}]$-splitting $\fb = \fs \oplus  [e, \fu]$ such that $\fs \subset \fg$ is graded for the filtration $\fg = \bigoplus_{i \in \bZ} \fg^i$ reviewed in \cref{footnote:Riche-argument}. A \emph{Kostant section} of $(\sR, \Delta)$ is the closed $\bZ[\f{1}{\tau_\sR}]$-subscheme
\[
\tst \cS \ce e + \fs \subset \fg^\reg  \subset \fg \qx{that depends on the choice of $\fs$,} 
\]
where \cite{Ric17}*{Equation (3.1.1)} ensures the factoring through $\fg^\reg$. By \cite{AFV18}*{Proposition 10}, the scheme $\cS$ is indeed a section to the map $\fg \ra \fg \!\sslash\! G$: for any $\bZ[\f{1}{\tau_\sR}]$-scheme $S$,
\[
\cS_S \isomto \fg_S \!\sslash\! G_S \overset{\ref{prop:base-change-for-g//G}}{\cong} (\fg \!\sslash\! G)_S.
\]
By construction, letting $2\rho \colon \bG_m \ra G$ denote the sum of all the positive coroots, the action of $\bG_m$ on  $\fg$ given by
\be \label{eqn:action-on-Kostant}
(t, \gamma) \mapsto t^{-2}\Ad((2\rho)(t))(\gamma)
\ee
preserves $\cS \subset \fg$. If we endow $\fg \!\sslash\! G$ with the $\bG_m$-action induced from $(t, \gamma) \mapsto t^{-2} \gamma$, then the isomorphism $\cS \isomto \fg \!\sslash\! G$ becomes $\bG_m$-equivariant because $\gamma$ and $\Ad((2\rho)(t))(\gamma)$ have the same image in $\fg \!\sslash\! G$. 
\epp

\brem \label{rem:canonical-S}
Keeping the setup \S\ref{pp:Kostant-section}, we assume that our split, pinned $G$ is instead defined over a scheme $S$ on which $\#(\pi_1(\sR^\ad))$ and every prime that appears in the expression of some root of $\sR$ in terms of a base of simple roots is invertible,\footnote{\label{foot:torsion-for-root-system}In terms of the appearing types, the following primes divide some coefficient in the expression of some root (equivalently, the highest root) in terms of a base of simple roots (compare with the shorter list in \S\ref{pp:torsion-primes}~\ref{item:torsion-prime-type-2}--\ref{item:torsion-prime-type-4}): 
\begin{itemize}
\item
$p = 2$ for types $B_n$, $C_n$, $D_n$, $E_6$, $E_7$, $E_8$, $F_4$, and $G_2$;

\item
$p = 3$ for types $E_6$, $E_7$, $E_8$, $F_4$, and $G_2$;

\item
$p = 5$ for type $E_8$.
\end{itemize}
The primes not in this list are often called ``good,'' as is done, for instance, in \cite{Spr66}*{Section 2.10}. The ``bad'' primes, that is, the ones in the list divide the order of the Weyl group, as does $\#(\pi_1(\sR^\ad))$.
}
and that one has a $G$-equivariant $S$-isomorphism $\iota\colon \fg \isomto \fg^*$ (equivalently, a perfect $G$-equivariant $S$-pairing on $\fg$).\footnote{\label{foot-Let}If $S$ is over an algebraically closed field $k$, then $\iota$ exists: indeed, the other assumptions imply that $\Char(k)$ must be ``very good for $G$'' in the sense of \cite{Let05}*{Definition 2.5.5}, so the criterion \cite{Let05}*{Proposition 2.5.12} for the existence of $\iota$ applies.} We claim that then 
\[
\tst \xq{there are canonical choices for}  \fs \qxq{and} \cS \qxq{over}\tst S \qxq{that depend on $\iota$.}
\]
Indeed, by \cite{Spr66}*{Proposition 2.4 and Theorem 2.6}, the cokernel of $[e, -] \colon \fg \ra \fg$ is locally free, so, since this map respects the grading discussed in \cref{footnote:Riche-argument}, it suffices to exhibit a canonical graded complement to $[e, \fg] \subset \fg$. In fact, $\Ker([e, -])$ is such a complement: by the flatness of $\fg/[e,\fg]$, this kernel is locally free, of formation compatible with base change, so we may assume that the base is a field and need to argue that $[e, \fg] \cap \Ker([e, -]) = 0$. However, $\Ker([e, -])$ is even orthogonal to $[e, \fg]$ under the perfect pairing $(\cdot, \cdot)_\iota$ on $\fg$ determined by $\iota$ because, by deriving the $G$-equivariance of $(\cdot, \cdot)_\iota$, we get 
\[
([\gamma, x], y)_\iota = (x, [-\gamma, y])_\iota \qxq{and may choose} \gamma = e, \q x \in \fg, \q y \in \Ker([e, -]).
\]
\erem

We now improve \cite{Ric17}*{Proposition 3.3.3} by showing that $\fg^\reg \ra \fg \!\sslash\! G$ is smooth under weaker hypotheses. 

\bprop \label{thm:describe-chi}
For a reductive group $G$ over a scheme $S$ such that $\sR(G)$ has no torsion residue characteristics, the map $\fg \ra \fg \!\sslash\! G$ \up{already discussed in Proposition \uref{prop:describe-Chevalley}} is finitely presented, flat, with reduced, local complete intersection geometric fibers, its restriction
\[
\fg^\reg \surjects \fg \!\sslash\! G \qxq{is smooth of relative dimension $\dim(G) - \rk(G)$ and surjective,}
\]
and  $\fg^\reg \subset \fg$ forms a single nonempty $G$-orbit over each geometric $(\fg \!\sslash\! G)$-fiber.
\eprop

\begin{proof}
By \Cref{prop:base-change-for-g//G}, the $S$-scheme $\fg \!\sslash\! G$ is finitely presented, of formation compatible with base change. Thus, by \cite{SP}*{Lemma \href{https://stacks.math.columbia.edu/tag/02FV}{02FV}}, the map $\fg \ra \fg \!\sslash\! G$ is also finitely presented and, by the fibral criterion \cite{EGAIV3}*{Corollaire 11.3.11}, the remaining claims are reduced to the case when $S$ is the spectrum of algebraically closed field $k$. The surjectivity of $\fg^\reg \ra \fg \!\sslash\! G$ then follows from the existence of a Kostant section $\cS \subset \fg^\reg$ (see \S\ref{pp:Kostant-section}) and the $G$-orbit claim is part of \cite{Jan04}*{Section~7.13, Proposition}. Thus, by $G(k)$-conjugation, for the smoothness of $\fg^\reg \ra \fg \!\sslash\! G$ it suffices to consider  the points of $\fg^\reg$ that lie on $\cS$, at which the map is indeed smooth by \cite{EGAIV4}*{Th\'{e}or\`{e}me 17.11.1} because it is surjective on tangent spaces due to $\cS \isomto \fg \!\sslash\! G$. Moreover, for dimension reasons, that is, by \Cref{prop:describe-Chevalley} and \cite{EGAIV2}*{Proposition 6.1.5}, the map $\fg \ra \fg \!\sslash\! G$ between smooth $k$-schemes is flat, and hence, by \cite{SP}*{Lemma \href{https://stacks.math.columbia.edu/tag/0E9K}{0E9K}}, it is automatically a relative complete intersection. In particular, the geometric fibers of $\fg \ra \fg \!\sslash\! G$ are Cohen--Macaulay and, due to the relative density of $\fg^\reg$, also generically smooth, so the (R$_0$)$+$(S$_1$) criterion \cite{SP}*{Lemma \href{https://stacks.math.columbia.edu/tag/031R}{031R}} ensures that they are reduced. 
\end{proof}

The smoothness of the map $\fg^\reg \surjects \fg\!\sslash\!G$ supplied by \Cref{thm:describe-chi} leads to the description of the basic properties of the centralizer of the universal section of $\fg^\reg$ in \Cref{pla}. The latter generalizes \cite{Ric17}*{Corollary 3.3.8, Proposition 3.3.9} and \cite{Ngo10}*{Lemma 2.1.1, Proposition 2.2.1} and rests on the following lemma. 


\blem \label{adjfond2}
For a root-smooth reductive group $G$ over a scheme $S$ such that $\Char(k(s))$ is prime to $\# (\pi_1(G_{\der, \ov{s}}))$ for $s \in S$,\footnote{This condition holds if $\sR(G)$ has no torsion residue characteristics, see \uS\uref{pp:torsion-primes}} a maximal $S$-torus $T \subset G$, and its Lie algebra $\ft \subset \fg$, the Weyl group 
\[
W \ce N_G(T)/T \qxq{acts freely on} \ft^{\rs} \ce \ft \cap \fg^{\rs} \qxq{and}  T \isomto C_G(\gamma)    \qxq{for any} \gamma \in \ft^\rs(S).
\] 
\elem

\bpf
The freeness of the $W$-action amounts to the following map being an isomorphism:
\[
\ft^\rs \xra{\gamma\, \mapsto\, ((1, \gamma), \gamma)} (W \times_S \ft^\rs) \times_{\ft^\rs \times_S \ft^\rs,\, \Delta} \ft^\rs,  \qxq{where} W \times_S \ft^\rs \xra{(w,\, \gamma)\, \mapsto\, (w\gamma,\, \gamma)} \ft^\rs \times_S \ft^\rs.
\]
Thus, the fibral criterion \cite{EGAIV4}*{Corollaire 17.9.5} reduces us to the case when $S = \Spec(k)$ for an algebraically closed field $k$. For such $S$, the assertion about the $W$-action follows from \cite{Ric17}*{Lemma~2.3.3}, and, once the freeness of the $W$-action is known, $T \isomto C_G(\gamma)$ follows from \Cref{lem:description-of-centralizer}. 
\epf

\bthm\label{pla} 
Let  $G$ be a reductive group over a scheme $S$ such that $\sR(G)$ has no torsion residue characteristics and let $C \subset G_\fg$ be the centralizer of the universal section of $\fg\ce \Lie(G)$. The $\fg^\reg$-group $C_{\fg^\reg}$ is commutative, flat, affine, a relative local complete intersection \up{so of finite presentation}, and there is a unique $(\fg\!\sslash\! G)$-group scheme $J$ equipped with a $G$-equivariant isomorphism 
\[
J_{\kg^{\reg}}\simeq C_{\kg^{\reg}}, \qxq{which then extends to a unique $G$-equivariant $\fg$-group map} J_\fg \ra C,
\]
where the $G$-actions on scheme valued points are described by
\[
g \cdot (j, \gamma) = (j, \Ad(g)\gamma) \qxq{and} g \cdot (c, \gamma) \mapsto (gcg\i, \Ad(g)\gamma).
\]
Moreover, $\fg\!\sslash\! G$ is the coarse moduli space of the algebraic stack quotient $[\fg^\reg/G]$, more precisely,
\be \label{eqn:greg-gerbe}
[\fg^\reg/G] \ra \fg\!\sslash\! G \qxq{is a gerbe bound by} J \qx{for the fppf topology},
\ee
and if $G$ is root-smooth, then $J$ \up{respectively, $C$} is a torus over the image $(\fg\!\sslash\!G)^\rs$ of $\fg^\rs$ \up{respectively, over $\fg^\rs$}. 
\ethm

\begin{proof}
We begin with the claims about $C_{\fg^\reg}$ and note that the parenthetical claim will follow from the local complete intersection aspect \cite{SP}*{Lemma \href{https://stacks.math.columbia.edu/tag/069H}{069H}}. For these claims, by the fpqc local on the base nature of local complete intersection morphisms and their stability under base change \cite{SP}*{Lemmas \href{https://stacks.math.columbia.edu/tag/069K}{069K}, \href{https://stacks.math.columbia.edu/tag/02VK}{02VK}, and \href{https://stacks.math.columbia.edu/tag/01UI}{01UI}}, we may first work \'{e}tale locally $S$ and then reduce to  $S$ being a localization of $\Spec(\bZ)$. In the Cartesian square
$$\qq \xymatrix@C=40pt{C_{\fg^\reg}\ar[r]\ar[d]&G\times_S\kg^{\reg}\ar[d]^-{(g,\,\g)\,\mapsto\,(\g,\,\Ad(g)\g)}\\\kg^{\reg}\ar[r]^-{\gamma\, \mapsto\, (\gamma,\, \gamma)}&\kg^{\reg}\times_{\fg\sslash G}\kg^{\reg}}$$
the $S$-scheme $\kg^{\reg}\times_{\fg\sslash G}\kg^{\reg}$ is smooth by \Cref{thm:describe-chi}, so the bottom horizontal map is a regular closed immersion. Due to this smoothness and an $S$-fibral dimension count that uses Propositions \ref{prop:base-change-for-g//G} and \ref{thm:describe-chi}, the miracle flatness \cite{EGAIV2}*{Proposition 6.1.5}, and the fibral criterion \cite{EGAIV3}*{Corollaire 11.3.11}, the finitely presented map 
\be \label{eqn:orbit-map-fflat}
\xymatrix{G\times_S\kg^{\reg}\ar@{->>}[rrr]^-{(g,\,\g)\,\mapsto\,(\g,\,\Ad(g)\g)} &&& \kg^{\reg}\times_{\fg\sslash G}\kg^{\reg}} \qx{is faithfully flat.}
\ee
Consequently, $C_{\fg^\reg}$ is flat over $\fg^\reg$ and is locally cut out by a regular sequence in the smooth $S$-scheme $G \times_S \fg^\reg$. Thus, by, for instance, \cite{SGA2new}*{Expos\'{e} XIV, Proposition 5.2}, the ideal of any closed immersion $C_{\fg^\reg} \hra \bA^n_{\fg^\reg}$ is locally generated by a regular sequence, so $C_{\fg^\reg}$ is a relative local complete intersection over $\fg^\reg$. Due to the $\fg^\reg$-flatness of $C_{\fg^\reg}$, its commutativity can be tested after base change to a geometric generic point $\ov{\eta}$ of $\fg^\reg$. Since $\ov{\eta}$ lies over $\bQ$ and is regular semisimple (see \S\ref{pp:regular-semisimple}), \Cref{adjfond2} ensures that the centralizer $C_{\ov{\eta}}$ is a torus, and hence is commutative.

Next we turn to the claims about $J$, so we revert to $S$ being arbitrary.  Since $C$ is $\fg$-affine, by descent, the existence and uniqueness of $J$ amount to the group $C_{\fg^\reg}$ admitting a unique $G$-equivariant descent datum $\iota$ with respect to the smooth surjection $\fg^\reg \surjects \fg\!\sslash\!G$ of \Cref{thm:describe-chi}. For this, we let $S'$ be a variable $S$-scheme and view $\iota$ as a  collection of group isomorphisms
\[
\iota_{\gamma_1,\, \gamma_2} \colon C_G(\gamma_1) \isomto C_G(\gamma_2) \qxq{for all} \gamma_1, \gamma_2 \in \fg^\reg(S') \qxq{that have the same image in} (\fg\!\sslash\!G)(S') 
\]
such that the $\iota_{\gamma_1,\, \gamma_2}$ are compatible with pullback in $S'$ and are subject to the cocycle condition 
\[
\iota_{\gamma_2,\, \gamma_3} \circ \iota_{\gamma_1,\, \gamma_2} = \iota_{\gamma_1,\, \gamma_3},
\]
compare with this kind of description of descent data in \cite{BLR90}*{Section 6.1, before Lemma 3}; concretely, 
\[
\iota_{\gamma_1,\, \gamma_2}  \qxq{is} C_G(\gamma_1) \isomto J \times_{\fg\sslash G,\, \gamma_1} S' \cong J \times_{\fg\sslash G,\, \gamma_2} S' \overset{\sim}{\longleftarrow} C_G(\gamma_2).
\]
The $G$-equivariance amounts to
\[
\iota_{\gamma,\, \Ad(g)\gamma}\colon C_G(\gamma) \isomto C_G(\Ad(g)\gamma) \qxq{being given by conjugation by $g$ for} \gamma \in \fg^\reg(S'),\, g \in G(S'). 
\]
By \eqref{eqn:orbit-map-fflat}, fppf locally on $S'$ every $(\gamma_1, \gamma_2)$ is of this form, so the $G$-equivariance determines $\iota_{\gamma_1,\, \gamma_2}$ uniquely, the cocycle condition is automatic, and the existence of $\iota$ reduces to $\iota_{\gamma,\, \Ad(g)\gamma}$ not depending on $g$. However, if $\Ad(g)\gamma = \Ad(g')\gamma$, then $g' = gc$ with $c \in (C_G(\gamma))(S')$, so the independence follows from the fact that the conjugation by $c$ has no effect because the group $C_G(\gamma)$ is commutative.

By descent, the group $J$ inherits properties from $C_{\fg^\reg}$: it is commutative, flat, affine, and a relative complete intersection over $\fg\!\sslash\!G$. Moreover, \eqref{eqn:gfin-fiberwise-large} and \Cref{prop:characterization-of-greg} ensures that $J_{\fg^\reg}$ is $S$-fiberwise of codimension $\ge 2$ in $J_\fg$. Thus, since $C$ is affine, the relative version of Hartogs' extension principle \cite{EGAIV4}*{Proposition 19.9.8} supplies the unique morphism $J_\fg \ra C$ extending $J_{\kg^{\reg}}\simeq C_{\kg^{\reg}}$ and, by allowing one to check the commutativity of all the relevant diagrams over $\fg^\reg$, shows that it is a $G$-equivariant $\fg$-group homomorphism.


The algebraicity of the stack $[\fg^\reg/G]$ is ensured by the general criterion \cite{SP}*{Theorem \href{https://stacks.math.columbia.edu/tag/06FI}{06FI}}. By \Cref{thm:describe-chi}, the map $[\fg^\reg/G] \ra \fg\!\sslash\! G$ is an fppf (even smooth) surjection and, by \eqref{eqn:orbit-map-fflat}, two objects of $[\fg^\reg/G]$ above the same scheme-theoretic point of $\fg\!\sslash\! G$ are fppf locally isomorphic. Thus, 
\[
[\fg^\reg/G] \ra \fg\!\sslash\! G
\]
is a gerbe by the criterion \cite{LMB00}*{Remarques 3.16}, and it is bound by $J$ by the construction of the latter.

By \Cref{prop:describe-Chevalley}, if $G$ is root-smooth, then the preimage of $(\fg\!\sslash\! G)^\rs$ in $\fg$ is precisely $\fg^\rs$. Thus, the last assertion follows from \S\ref{pp:regular-semisimple} and \Cref{adjfond2}.
\end{proof}

\brem 
Of course, the gerbe $[\fg^\reg/G] \ra \fg\!\sslash\! G$ is neutral whenever the map $\fg^\reg \ra \fg\!\sslash\! G$ has a section, such as a Kostant section when $G$ is split and pinned   (see \S\ref{pp:Kostant-section}, as well as \Cref{kosgm} below for a quasi-split case).
\erem

\brem \label{rem:J-with-an-L}
Both $\fg^\reg$ and $\fg$ carry compatible scaling by $\bG_{m,\, S}$ actions, which lift to a $\bG_{m,\, S}$-action  on $C$. The descent datum $\iota$ that gives $J$ is equivariant for this action, so $\bG_{m,\, S}$ also acts on $J$ compatibly with its action on $\fg\!\sslash\! G$. Consequently, for every $\bG_{m,\, S}$-torsor $\sL$ with the twist 
\[
\fg_\sL \ce \fg \times^{\bG_m}_S \sL,
\]
we obtain a unique ($\fg_\sL \!\sslash\! G$)-group $J$ equipped with a $G$-equivariant isomorphism 
\be \label{eqn:J-L-descent}
J_{\fg^\reg_\sL} \simeq C_{\fg^\reg_\sL}, \qxq{which extends to a $G$-equivariant $\fg_\sL$-group map} J_{\fg_\sL} \ra C_{\fg_\sL},
\ee
and, by \Cref{pla}, the map $[\fg^\reg_\sL/G] \ra \fg_\sL\!\sslash\! G$ is a gerbe bound by $J$.
\erem

Under the following stronger assumption on the residue characteristics, $C_{\fg^\reg}$ is even $\fg^\reg$-smooth. 

\bprop  \label{prop:J-is-smooth}
Let $G$ be a reductive group over a scheme $S$ such that $\Char(k(s))$ for $s\in S$ divides neither $\#\pi_1(((G_{\ov{s}})_\der)^\ad)$ nor any coefficient in the expression of a root of $G_{\ov{s}}$ in terms of a base of simple roots. The $(\fg\!\sslash\! G)$-group $J$ defined in Theorem \uref{pla} is smooth and the centralizer 
\[
C_{\fg^\reg} \subset G_{\fg^\reg} \qx{of the universal regular section of the Lie algebra $\fg$ of $G$ is $\fg^\reg$-smooth.}
\]
\eprop

\bpf
The assumptions imply that $\sR(G)$ has no torsion residue characteristics (see \S\ref{pp:torsion-primes} and \cref{foot:torsion-for-root-system}). Thus, the claim about $J$ follows from the rest by descent (see \Cref{thm:describe-chi}). 

For the rest, the $\fg^\reg$-flatness of $C_{\fg^\reg}$ established in \Cref{pla} allows us to assume that $S$ is the spectrum of an algebraically closed field $k$. Moreover, the open locus of $\fg^\reg$ over which $C_{\fg^\reg}$ is smooth is open and stable  both under the $\bG_m$-scaling and the adjoint action of $G$. Thus, since the map $\fg^\reg \ra \fg\!\sslash\! G$ is open and its geometric fibers consist of single $G$-orbits (see \Cref{thm:describe-chi}), the open of $\fg^\reg$ over which $C_{\fg^\reg}$ is smooth is a preimage of a $\bG_m$-stable open $U \subset \fg\!\sslash\! G$. The $\bG_m$-action on $\fg\!\sslash\! G$ extends to a map 
\[
\bA^1 \times_k \fg\!\sslash\! G \ra \fg\!\sslash\! G \qxq{that maps} \{ 0 \} \times_k \fg\!\sslash\! G \qxq{to} \{ 0\} \subset \fg\!\sslash\! G.
\]
Thus, once we show that $\{ 0 \} \subset U$ in $\fg\!\sslash\! G$, it will follow that the preimage $V \subset \bA^1 \times_k \fg\!\sslash\! G$ of $U$ under the action map $\bG_m \times_k \fg\!\sslash\! G \ra \fg\!\sslash\! G$ contains $\{ 0 \} \times_k \fg\!\sslash\! G$. Then $V$ will meet every geometric fiber of the projection $\bG_m \times_k \fg\!\sslash\! G \ra \fg\!\sslash\! G$, and hence, being stable under the $\bG_m$-action on the first coordinate, $V$ will contain every such fiber. The desired $U = \fg\!\sslash\! G$ will follow.

To argue the remaining inclusion $\{0 \} \subset U$, since the fiber of $\fg^\reg \ra \fg\!\sslash\! G$ above $\{0\}$ consists of the regular nilpotent sections (see \Cref{prop:describe-Chevalley}), we need to show that the centralizer $C_G(\gamma)$ of every regular nilpotent $\gamma \in \fg^\reg(k)$ is smooth. Such $\gamma$ form a single $G(k)$-conjugacy class (see \S\ref{pp:regular-nilpotent}), so it suffices to show that $C_G(\gamma)$ is smooth when 
\[
\tst \gamma \ce \sum_{\al\in\Delta}e_{\al} \qxq{for some pinning} \{e_\gA\}_{\gA \in \Delta} \qxq{of} G.
\]
 Then, letting $Z \subset G$ be the maximal central torus, we have the group homomorphism 
\[
Z \times C_{G_\der}(\gamma) \ra C_G(\gamma),
\]
 which is surjective and hence also flat. It remains to note that $C_{G_\der}(\gamma)$ is smooth by \cite[Theorem~5.9~a)]{Spr66} (which even shows that our assumption on $\Char(k(s))$ is sharp).
\epf

The ``Galois-theoretic'' description of the torus $J|_{(\fg\sslash G)^\rs}$ presented in \Cref{prop:galois-theoretic} is key for \S\ref{section:Hitchin-fibration} and mildly generalizes \cite{Ngo10}*{Proposition 2.4.2} by weakening its assumption that the cardinality of the Weyl group be invertible on the base (in turn, \emph{loc.~cit.}~builds on \cite{DG02}*{Section 11} that considered a situation over $\bC$). The proof of this description will use the following lemma, which improves \cite{Ric17}*{Lemma 3.5.3}.


\blem \label{lem:Chevalley-fiber-prod}
Let $G$ be a root-smooth reductive group over a scheme $S$ such that $\sR(G)$ has no torsion residue characteristics, let $T \subset G$ be a maximal $S$-torus with Lie algebra $\ft \subset \fg$ and let $W \ce N_G(T)/T$ be the Weyl group scheme. The map $\ft \surjects \ft/W$ is finite locally free of degree $\#W$ and is a $W$-torsor over $\ft^\rs/W$ and for any Borel $S$-subgroup $T \subset B \subset G$ the maps 
\[
\wt{\fg} \xra{\eqref{eqn:gtilde-to-t}} \ft \qxq{and} \fg \ra \fg\!\sslash\!G \overset{\ref{adj}}{\cong} \ft/W
\]
induce an isomorphism
\be \label{describe-wt-greg}
\widetilde{\kg}^{\reg}\isomto \kg^{\reg}\times_{\ft/W}\kt \qxq{over} \fg^\reg.
\ee
\elem

\bpf
The map \eqref{describe-wt-greg} is well-defined by \eqref{eqn:cool-diagram}.
By \Cref{adj} and \Cref{prop:base-change-for-g//G}, we have $\ft/W \cong \fg\!\sslash\! G$ compatibly with base change, so, after first reducing to the split, pinned case, we may assume that $S$ is a localization of $\Spec(\bZ)$. Over the geometric $S$-fibers the finite map $\ft \ra \ft/W$ is a morphism of affine spaces of the same dimension, so, by the fibral criterion \cite{EGAIV3}*{Corollaire~11.3.11}, it is flat; its degree may be read off over $\ft^\rs/W$, over which it is a $W$-torsor by \Cref{adjfond2}. By \Cref{thm:describe-chi}, the $S$-scheme $\kg^{\reg}\times_{\ft/W}\kt$ is smooth,  so it is the normalization of $\fg^\reg$ in $(\kg^{\reg}\times_{\ft/W}\kt)|_{\fg^\rs}$. The same holds for $\wt{\fg}^\reg$ (see \S\ref{pp:GS-map} and \Cref{prop:characterization-of-greg}), so we need to argue that \eqref{describe-wt-greg} is an isomorphism over $\fg^\rs$. By \Cref{prop:describe-restriction-to-rs} and \Cref{adjfond2}, both sides of \eqref{describe-wt-greg} are $W$-torsors over $\fg^\rs$, so it suffices to note that the map is $W$-equivariant: in terms of the isomorphism $G/T \times_S \ft^\rs\cong \wt{\fg}^\rs$ of \Cref{prop:describe-restriction-to-rs}, the map to $\ft$ is simply the projection.
\epf

\begin{prop} \label{prop:galois-theoretic}
Let $G$ be a root-smooth reductive group over a scheme $S$ such that $\sR(G)$ has no torsion residue characteristics, let $T \subset G$ be a maximal $S$-torus with Lie algebra $\ft \subset \fg$, set $\ft^\rs \ce \ft\cap \fg^\rs$, let $W \ce N_G(T)/T$ be the Weyl group, let $T \subset B \subset G$ be a Borel $S$-subgroup, and let $\pi \colon \ft \surjects \ft/W$ be the indicated quotient map. The group $J$ 
defined in Theorem \uref{pla} admits a 
\be \label{eqn:galois-theoretic-map}
\xq{$(\ft/W)$-group homomorphism} J \ra (\pi_*(T _ \ft))^W \qxq{that is an isomorphism over} \ft^\rs/W,
\ee
where  and $W$ acts on the restriction of scalars $\pi_*(T _ \ft)$ via its actions on $T$ and $\ft$. 
In particular, $J|_{\ft^\rs/W}$ becomes isomorphic to a base change of $T$ over the finite \'{e}tale cover $\ft^\rs \surjects \ft^\rs/W$. 
\end{prop}

\begin{proof}
The proof is similar to that of \cite{Ngo10}*{Proposition 2.4.2}, but we include it since our assumptions are slightly weaker. Firstly, $\ft \ra \ft/W$ is finite locally free and a $W$-torsor over $\ft^\rs/W$ by \Cref{lem:Chevalley-fiber-prod}, so, by \cite{BLR90}*{Section 7.6, Theorem 4, Proposition 5}, the restriction of scalars $\pi_*(T \times_S \ft)$ is a smooth, affine $(\ft/W)$-group and it suffices to settle \eqref{eqn:galois-theoretic-map}. For the latter, we will use the Cartesian square supplied by \Cref{lem:Chevalley-fiber-prod} using our fixed Borel $B$:
\[
\xymatrix@C=35pt{
\wt{\fg}^\reg \ar@{->>}[r]^-{\eqref{eqn:gtilde-to-t}} \ar@{->>}[d]^-{\pi'} & \ft \ar@{->>}[d]^-{\pi} \\
\fg^\reg \ar@{->>}[r] & \ft/W
}
\]
in which the vertical arrows are finite locally free and the horizontal ones are smooth by \Cref{thm:describe-chi}. In particular, since the map \eqref{eqn:gtilde-to-t} is $W$-equivariant (see \Cref{prop:extend-W-action-to-greg}) and $\ft/W \cong \fg\!\sslash\!G$ (see \Cref{adj}), by arguing by descent as in the proof of \Cref{pla} it suffices to produce a $G$-equivariant homomorphism 
\[
\tst J_{\fg^\reg} \ra (\pi'_*(T _ {\wt{\fg}^\reg}))^W
\]
that is an isomorphism over $\fg^\rs$ (that is, over the preimage of $\ft^\rs/W$, see \Cref{prop:describe-Chevalley}). The group $J_{\fg^\reg}$ is $G$-equivariantly identified with the centralizer $C_{\fg^\reg}$ of the universal regular section of $\fg$ (see \Cref{pla}), the restriction of $\pi'$ to $\fg^\rs$ is a $W$-torsor, so that 
\[
\fg^\rs \isomto (\pi'_*(\wt{\fg}^\rs))^W
\]
and likewise for $C_{\fg^\rs}$, and the $G$- and $W$-actions on $\wt{\fg}^\reg$ commute (see \Cref{prop:extend-W-action-to-greg}). Thus, all we need to do is to exhibit a $G$-equivariant and $W$-equivariant homomorphism $C_{\wt{\fg}^\reg} \ra T_{\wt{\fg}^\reg}$
that is an isomorphism over $\wt{\fg}^\rs$. 

The group $G_{\wt{\fg}^\reg}$ comes equipped with two subgroups: the base change $C_{\wt{\fg}^\reg}$ of the universal centralizer and the base change $\cB_{\wt{\fg}^\reg}$ of the universal Borel, and we claim that $C_{\wt{\fg}^\reg} \subset \cB_{\wt{\fg}^\reg}$. For this we may first work \'{e}tale locally on $S$ and then reduce to $S$ being a localization of $\Spec(\bZ)$. Moreover, by \Cref{pla}, the $\wt{\fg}^\reg$-group $C_{\wt{\fg}^\reg}$ is flat, so it suffices to check that $C_{\wt{\fg}^\rs} \subset \cB_{\wt{\fg}^\rs}$. By \Cref{pla} again, $C_{\wt{\fg}^\rs}$ is reduced, so this inclusion may be checked fiberwise and amounts to the assertion that for a geometric $S$-point $\ov{s}$, the $G$-centralizer of a $\gamma \in \fg^\rs(\ov{s})$ lies in every Borel whose Lie algebra contains $\gamma$. This, however, follows from the Jordan decomposition and \Cref{adjfond2}.

The universal Borel $\cB$ may be identified with the universal $(G/B)$-conjugate of our fixed Borel $T \subset B$ and, in particular, the quotient by its unipotent radical is canonically a base change of $T$. Thus, the inclusion $C_{\wt{\fg}^\reg} \subset \cB_{\wt{\fg}^\reg}$ gives a desired homomorphism $C_{\wt{\fg}^\reg} \ra T_{\wt{\fg}^\reg}$ that is $G$-equivariant and $W$-equivariant by construction: for instance, to check the $W$-equivariance, one may work over $\wt{\fg}^\rs$, use \Cref{prop:describe-restriction-to-rs} to identify $\wt{\fg}^\rs \cong G/T \times_S \ft$, and then note that on the quotient $\cB_{\wt{\fg}^\reg} \surjects T_{\wt{\fg}^\reg}$ the difference between the ``conjugating back'' by $g\i$ and by $(gw\i)\i$ is the action by $w$ on $T$. Finally, the map $C_{\wt{\fg}^\reg} \ra T_{\wt{\fg}^\reg}$ is an isomorphism over $\wt{\fg}^\rs$ because it is so fiberwise by \Cref{adjfond2}.
\end{proof}

The following result about the conjugacy class of the Kostant section 
illustrates the utility of the Galois-theoretic description of $J$ supplied by \Cref{prop:galois-theoretic} and complements the fact that the geometric fibers of $\fg \ra \fg\!\sslash\! G$ above points in $(\fg\!\sslash\! G)^\rs$ consist of single $G$-orbits (see \Cref{prop:describe-Chevalley}). 


\bthm\label{conj}
For a seminormal, strictly Henselian, local ring $R$ and a reductive $R$-group $G$ with Lie algebra $\fg$, if the order of the Weyl group of $G$ is invertible in $R$, then the
\[
\x{fibers of}\ \ \fg^\rs(R\llp t \rrp) \surjects (\fg \!\sslash\! G)^\rs(R\llp t \rrp)\ \ \x{are precisely the $G(R\llp t\rrp)$-conjugacy classes in $\fg^\rs(R\llp t \rrp)$.}
\] 
\ethm

\begin{proof}
Since the order of the Weyl group is invertible in $R$, the group $G$ satisfies all the ``no small residue characteristics'' assumptions that appear earlier in this chapter (see \S\ref{pp:root-smoothness}, \S\ref{pp:torsion-primes}, and \cref{foot:torsion-for-root-system}), so we may freely apply the preceding results in this proof. In the statement, we let 
\[
(\fg \!\sslash\! G)^\rs \subset \fg \!\sslash\! G
\]
denote the open image of $\fg^\rs$ (see \Cref{thm:describe-chi}); by \Cref{prop:describe-Chevalley}, its preimage in $\fg$ is precisely $\fg^\rs$. Thanks to a Kostant section (see \S\ref{pp:Kostant-section}), the map is surjective as indicated.

The map $\fg^\rs \ra \fg \!\sslash\! G$ is invariant under $G$-conjugation, so $G(R\llp t \rrp)$-conjugate elements of $\fg^\rs(R\llp t \rrp)$ agree in $(\fg \!\sslash\! G)(R\llp t \rrp)$. Conversely, fix $\gamma_1, \gamma_2 \in \fg^\rs(R\llp t \rrp)$ that have a common image 
\[
\ov{\gamma} \in (\fg \!\sslash\! G)^\rs(R\llp t \rrp) \subset (\fg \!\sslash\! G)(R\llp t \rrp).
\]
\Cref{pla} ensures that 
\[
[\fg^\rs/G] \surjects (\fg \!\sslash\! G)^\rs
\]
is a gerbe bound by $J_{(\fg \sslash G)^\rs}$, so the functor that parametrizes isomorphisms between the images of $\gamma_1$ and $\gamma_2$ in $[\fg^\rs/G]$ is a torsor under $J_{R\llp t \rrp}$ (pullback of $J$ along $\ov{\gamma}$). By \Cref{prop:galois-theoretic}, this $J_{R\llp t \rrp}$ is an $R\llp t \rrp$-torus that splits over a $W$-torsor for some finite group $W$ whose order is invertible in $R$. Thus, by \Cref{twisted}, the torsor in question is trivial and the images of $\gamma_1$ and $\gamma_2$ in $([\fg^\rs/G])(R\llp t \rrp)$ may be identified. The fiber of the map 
\[
\fg^\rs \ra [\fg^\rs/G]
\]
over this common image is a $G_{R\llp t \rrp}$-torsor trivialized both by $\gamma_1$ and $\gamma_2$, so the latter are indeed $G(R\llp t \rrp)$-conjugate, as desired.
\end{proof}


\csub[The product formula for the Hitchin fibration beyond the anisotropic locus] \label{section:Hitchin-fibration}

Our final goal is Ng\^{o}'s product formula over the entire $\cA^\heart$ stated precisely in \Cref{thm:product-formula} below. The construction of the stack morphism that encodes this formula amounts to glueing torsors with the help of a twisted Kostant section for quasi-split groups, so our first goal is to construct this section in \Cref{kosgm}. We begin with a brief review of quasi-splitness in order to remind that beyond semilocal bases this is a more stringent condition than the existence of a Borel.

\bpp[Quasi-split reductive groups] \label{pp:quasi-split-review}
We recall from \cite{SGA3IIInew}*{Expos\'{e} XXIV, Section 3.9} that a reductive group $G$ over a scheme $S$ is \emph{quasi-split} if it has a maximal $S$-torus and a Borel $S$-subgroup $T \subset B \subset G$ such that on the $S$-scheme $\underline{\mathrm{Dyn}}(G)$ of Dynkin diagrams the line bundle $\fg_*$ given by the universal root space that is simple with respect to $B$ is trivial. In this case, a choice of $T \subset B \subset G$ and a trivialization $e$ of $\fg_*$ constitutes a \emph{quasi-pinning} of $G$. For example, when $G$ is split with respect to $T$, a choice of $B$ amounts to that of a system $\Delta$ of positive simple roots, 
\[
\tst \underline{\mathrm{Dyn}}(G) \cong \bigsqcup_{\gA \in \Delta} S,
\]
the line bundle $\fg_*$ is given by the $T$-root space $\fg_\gA$ on the copy of $S$ indexed by $\gA$, and $e$ amounts to a trivialization of each $\fg_\gA$ (equivalently, to a nilpotent section of $\fg \ce \Lie(G)$ that is principal with respect to $\Delta$). Thus, in the split case the datum of a quasi-pinning amounts to that of a pinning.
\epp

We are ready to build a twisted Kostant section under more general conditions than in \cite[Lemme~2.2.5]{Ngo10}.

\begin{prop}\label{kos2}\label{kosgm}
Let $T \subset B \subset G$ be a quasi-split reductive group over a scheme $S$ with Lie algebras $\ft \subset \fb \subset \fg$ and Weyl group $W \ce N_G(T)/T$, let $\bG$ be the reductive $S$-group that is the form of $G$ that is split Zariski locally on $S$,\footnote{We recall from \cite{SGA3IIInew}*{Expos\'{e} XXII, D\'{e}finition 1.13} that the root datum of a split group is necessarily constant on the base.} and suppose that there is a $\bG$-isomorphism 
\[
\iota\colon \Lie(\bG)\isomto \Lie(\bG)^{*}
\]
and that $\Char(k(s))$ for $s \in S$ divides neither~$\#\pi_1(((G_{\ov{s}})_\der)^\ad)$ nor any coefficient in the expression of a root of $G_{\ov{s}}$ in terms of a base of simple roots. The Chevalley morphism
\be \label{eqn:quasi-split-Kostant}
 \fg \ra \ft/W \qxq{admits a section} \eps\colon \ft/W \ra \fg^\reg
\ee
and for any $\bG_m$-torsor $\sL$ on $S$ with $\ft_\sL \ce \ft \times^{\bG_m}_S \sL$ and $\fg_\sL \ce \fg \times^{\bG_m}_S \sL$, the induced $S$-morphism
\be \label{eqn:Kostant-Hitchin}
[\fg_{\sL^{\tensor 2}}/G] \ra \ft_{\sL^{\tensor 2}}/W \qxq{admits a section} \eps_\sL \colon \ft_{\sL^{\tensor 2}}/W \ra [\fg^\reg_{\sL^{\tensor 2}}/G].
 \ee
\end{prop}


\begin{proof}
The assumption on $\#\pi_1(((G_{\ov{s}})_\der)^\ad)$ implies that $G$ is root-smooth (see \S\ref{pp:root-smoothness}), 
so the Chevalley morphism makes sense (see \Cref{adj}). We are assuming that $G$ is quasi-split with respect to $T \subset B$, so we choose a trivialization $e$ that extends them to a quasi-pinning of $G$ (see \S\ref{pp:quasi-split-review}). By passing to the clopens of $S$ on which $\sR(G_{\ov{s}})$ is constant, we assume that $\bG$ is split.  

For \eqref{eqn:quasi-split-Kostant}, by descent, we may assume that $G = \bG$ and that $G$ is split with respect to $T$ at the expense of needing to check that the canonical Kostant section $\cS$ built in \Cref{rem:canonical-S} using $e$ and $\iota$ is invariant under any automorphism of $\bG$ that preserves the pinning. The construction of  $\cS$ gives this invariance because it only involves structures respected by every automorphism of $\bG$ that preserves $T$, $B$, and $e$ (in particular, the formula for $\cS$ does not use $\iota$, which need not be preserved by such an automorphism). To deduce \eqref{eqn:Kostant-Hitchin}, we first note that $\fg_{\sL^{\tensor 2}} \cong \fg'_\sL$ and $\ft_{\sL^{\tensor 2}} \cong \ft'_\sL$, where $(-)'$ means that the contracted product is formed with $\bG_m$ acting on $\fg$ and $\ft$ by the square of its scaling action (in terms of a local trivialization $\ell$ for $\sL$, the identifications are induced by $(\gamma, \ell^{\tensor 2}) \mapsto (\gamma, \ell)$). To then construct the desired section
\[
\eps_\sL\colon \ft'_{\sL}/W \ra [\fg'^{\,\reg}_{\sL}/G] \qxq{of} [\fg'_{\sL}/G] \ra \ft'_{\sL}/W,
\]
it suffices to find a functorial isomorphism given by $G$-conjugation that transforms the map
\[
\bG_m \times_S \ft/W \times_S \sL \ra \fg'^{\,\reg}_{\sL} \qxq{given by} (t, (\tau, \ell)) \mapsto (\eps(\tau), \ell) \qxq{into} (t, (\tau, \ell)) \mapsto (\eps(t^{-2}\tau), t\ell).
\]
Letting $2\rho \colon \bG_m \ra G$ be the sum of those coroots that are positive with respect to $B$, we have
\[
(\eps(t^{-2}\tau), t\ell) = (t^2\eps(t^{-2}\tau), \ell) \overset{\eqref{eqn:action-on-Kostant}}{=} (\Ad((2\rho)(t))\eps(\tau), \ell) \qxq{in} \fg'^{\,\reg}_{\sL},
\]
so the desired functorial conjugation is $t\mapsto \Ad((2\rho)(t))$. 
\end{proof}

The subsequent \S\S\ref{pp:Hitchin-fibration}--\ref{pp:symmetries} review the main actors that appear in the product formula.

\bpp[The Hitchin fibration]\label{pp:Hitchin-fibration}
Let $S$ be a scheme, let $\pi \colon X\ra S$ be a proper, smooth scheme morphism with connected geometric fibers of dimension $1$, let $G$ be a reductive $X$-group with Lie algebra $\fg$, and let $\sL$ be a $\bG_m$-torsor on $X$. The \emph{total Hitchin space} associated to this data is the restriction of scalars
\[
\cM_{\sL} \ce \pi_* ([\kg_{\sL}/G]), \qxq{where} \fg_\sL \ce \fg \times^{\bG_m}_X \sL \qxq{and $G$ acts via its adjoint action on $\fg$}
\]
(in spite of abusive notation, we stress that $G$ acts on $\fg_\sL$ on the left, compare with \S\ref{conv}). 

By the general criterion for the algebraicity of restrictions of scalars \cite{HR19}*{Corollary 9.2~(ii)} (or by the more specific \cite[Section 4.2.2]{Ngo10}), the stack $\cM_\sL$ is algebraic, locally of finite presentation over $S$, and has an affine diagonal. 
Concretely, for an $S$-scheme $S'$, the groupoid $\cM_\sL(S')$ consists of left $G_{X_{S'}}$-torsors $E$ equipped with a $G$-equivariant $X$-map $E \ra \fg_\sL$, equivalently, it is the groupoid of left $G_{X_{S'}}$-torsors  $E$  equipped with a section
\be \label{eqn:Hitchin-concrete}
 \phi \in H^0(X_{S'}, (E \times_{X_{S'}} (\fg_\sL)_{X_{S'}})/G_{X_{S'}}),
\ee
where $G_{X_{S'}}$ acts on both factors of $E \times_{X_{S'}} (\fg_\sL)_{X_{S'}}$. For a maximal $X$-torus $T \subset G$ with its Lie algebra $\ft \subset \fg$ and Weyl group $W \ce N_G(T)/T$, the \emph{Hitchin base} is the restriction of scalars of $\ft_\sL/W$:
\[
\cA_{\sL}:=\pi_*(\ft_{\sL}/W), \qxq{where} \ft_\sL \ce \ft \times^{\bG_m}_X \sL.
\]
If $G$ is root-smooth, then we have the regular semisimple locus 
\[
\ft^\rs_\sL/W \subset \ft_\sL/W, \qxq{where} \ft^\rs_\sL \ce \ft_\sL \cap \fg_\sL^\rs.
\]
Its preimage under the universal section $X_{\cA_\sL} \ra  (\ft_{\sL}/W)_{X_{\cA_\sL}}$ is an open $(X_{\cA_\sL})^\rs \subset X_{\cA_\sL}$ whose image is the open
\[
\cA^\heart_\sL \subset \cA_\sL \qx{over which the $(X_{\cA_\sL})^\rs$ is fiberwise dense in $X_{\cA_\sL}$.}
\]
The map $[\kg_{\sL}/G]\rightarrow\ft_{\sL}/W$ supplied by \Cref{adj} induces the \emph{Hitchin fibration} morphism
\be \label{eqn:HF}
f_{\sL}:\cM_{\sL}\rightarrow\cA_{\sL}.
\ee
If $G$ is quasi-split with respect to $T$ and, say, $\#W$ is invertible on $S$, then for every algebraically closed $S$-field $\ov{k}$, \Cref{kosgm} (with \Cref{rem:canonical-S}, especially, \cref{foot-Let}, to obtain the pairing $\iota$) supplies a Kostant section 
\[
\eps_{\sL}\colon \ft_{\sL^{\tensor 2}}/W\rightarrow[\kg_{\sL^{\tensor 2}}^{\reg}/G],
\]
which induces a \emph{Kostant--Hitchin section}
\[
\eps_{\sL} \colon \cA_{\sL^{\tensor 2},\, \ov{k}}\rightarrow \cM_{\sL^{\tensor 2}, \, \ov{k}} \qxq{of the $\ov{k}$-fiber of the Hitchin fibration} f_{\sL^{\tensor 2},\, \ov{k}} \colon \cM_{\sL^{\tensor 2},\, \ov{k}}\rightarrow\cA_{\sL^{\tensor 2},\, \ov{k}}.
\]
By construction, the Kostant--Hitchin section factors through the open substack
\[
\cM_{\sL^{\tensor 2}}^{\reg}:=\pi_*([\kg_{\sL^{\tensor 2}}^{\reg}/G]) \subset \cM_{\sL^{\tensor 2}}.
\]
\epp

\bpp[The affine Springer fibers]\label{pp:affine-Springer-fibration}
In the setting of \S\ref{pp:Hitchin-fibration}, suppose that $G$ is quasi-split and $\#W$ is invertible on $S$, let $\ov{k}$ be an algebraically closed $S$-field, let $a \in \cA_{\sL^{\tensor 2}}(\ov{k})$ be a $\ov{k}$-point as indicated, let the corresponding section of the map $(\ft_{\sL^{\tensor 2}}/W)_{X_{\ov{k}}} \ra X_{\ov{k}}$ be
\[
a\colon X_{\ov{k}} \ra (\ft_{\sL^{\tensor 2}}/W)_{X_{\ov{k}}} \qxq{and let} X^\rs_a \subset X_{\ov{k}} \qxq{be the $a$-preimage of} (\ft^\rs_{\sL^{\tensor 2}}/W)_{X_{\ov{k}}} \subset (\ft_{\sL^{\tensor 2}}/W)_{X_{\ov{k}}},
\]
so that, by definition, $a \in \cA^\heart_{\sL^{\tensor 2}}(\ov{k})$ if and only if $X^\rs_a \neq \emptyset$. Moreover, consider a $\ov{k}$-point $v$ of $X_{\ov{k}}$, let $\wh{\cO}_{v} \simeq \ov{k}\llb t_v \rrb$ be the completed local ring of $X_{\ov{k}}$ at $v$, and let 
\[
\tst a|_{\wh{\cO}_v} \in (\ft_{\sL^{\tensor 2}}/W)(\wh{\cO}_v)
\]
be the resulting $\wh{\cO}_v$-point. The \emph{affine Springer fiber} at $v$ is the functor $\cM_{\sL^{\tensor 2},\, a,\, v}$ that sends a $\ov{k}$-algebra $R$ to the groupoid of lifts
\[
\xymatrix{
\Spec(R\llb t_v\rrb) \simeq \Spec(\wh{\cO}_v \wh{\tensor}_{\ov{k}} R) \ar@{-->}[r] \ar[]!<1.8ex>;[rd]_-{a|_{\wh{\cO}_v}} &[\fg_{\sL^{\tensor 2}}/G] \ar[d] \\
& \ft_{\sL^{\tensor 2}}/W
}
\]
equipped with an isomorphism between the restriction to $(\wh{\cO}_v \wh{\tensor}_{\ov{k}} R)[\f{1}{t_v}]$ and the corresponding restriction of the Kostant--Hitchin lift $\eps_{\sL}(a)$. Analogously to after \eqref{eqn:Hitchin-concrete}, the data being parametrized amount to a $G_{\wh{\cO}_v \wh{\tensor}_{\ov{k}} R}$-torsor $E_v$ equipped both with a $G$-equivariant $X$-morphism $E_v \ra \fg_{\sL^{\tensor 2}}$ that lifts $a|_{\wh{\cO}_v}$ and an isomorphism after localizing away from $v$ with the analogous data determined by $\eps_{\sL}(a)$. Since $t_v$ is a nonzerodivisor in $R\llb t_v\rrb$, this rigidification with respect to the Kostant--Hitchin section eliminates nontrivial automorphisms, so the functor $\cM_{\sL^{\tensor 2},\, a,\, v}$ is set-valued. 

By \cite{Ngo10}*{Proposition 3.2.1} (see also \cite[Section 2, Proposition 1]{KL88} and \cite{Yun17}*{Theorem~2.5.2}), if $a|_{\wh{\cO}_v[\f{1}{t_v}]}$ factors through $\ft^\rs_{\sL^{\tensor 2}}/W$, as happens if and only if $a \in \cA^\heart_{\sL^{\tensor 2}}(\ov{k})$, then $\cM_{\sL^{\tensor 2},\, a,\, v}$ is representable by an ind-scheme 
whose associated reduced $\cM^\red_{\sL^{\tensor 2},\, a,\, v}$ is a locally of finite type, finite-dimensional $\ov{k}$-scheme. 
\epp

\bpp[Symmetries of the Hitchin and  affine Springer fibers]\label{pp:symmetries}
Assume the setting of \S\S\ref{pp:Hitchin-fibration}--\ref{pp:affine-Springer-fibration} with $G$ quasi-split and $\#W$ invertible on $S$, and a point $a \in \cA_{\sL^{\tensor 2}}(\ov{k})$. 
We use the corresponding section $a\colon X_{\ov{k}} \ra (\ft_{\sL^{\tensor 2}}/W)_{X_{\ov{k}}}$ to pull back  the descent $J$ of the universal regular centralizer constructed in \Cref{rem:J-with-an-L}. Thanks to the $G$-equivariant homomorphism 
\[
J_{\fg_{\sL^{\tensor 2}}} \xra{\eqref{eqn:J-L-descent}} C_{\fg_{\sL^{\tensor 2}}}
\]
to the universal $G$-centralizer, the resulting $X_{\ov{k}}$-group $J_a$ acts on the objects of the $X_{\ov{k}}$-stack $[\fg_{\sL^{\tensor 2}}/G]_{X_{\ov{k}}}$. 
 Thus, descent allows us to twist these objects by $J_a$-torsors, so the $a$-fiber $\cM_{\sL^{\tensor 2},\, a}$ of the Hitchin fibration \eqref{eqn:HF} admits an action of the Picard $\ov{k}$-stack 
 \[
 \cP_a \ce (\pi_{\ov{k}})_*(\bbB(J_a)),
 \]
 where $\bbB(-)$ denotes the classifying stack. By the general criterion \cite{HR19}*{Corollary 9.2~(ii)} for representability and properties of restrictions of scalars, the stack $\cP_a$ is algebraic and its diagonal is  affine.


Under the assumption that $a \in \cA^\heart_\sL(\ov{k})$, for each closed point $v \in X_{\ov{k}}$ we also consider the functor $\cP_{a,\, v}$ that for a variable $\ov{k}$-algebra $R$ parametrizes $J_a$-torsors over $\wh{\cO}_v \wh{\tensor}_{\ov{k}} R$ equipped with a trivialization over $(\wh{\cO}_v \wh{\tensor}_{\ov{k}} R)[\f{1}{t_v}]$. The same considerations as for $\cP_a$ show that $\cP_{a,\, v}$ acts on the affine Springer fiber $\cM_{\sL^{\tensor 2},\, a,\, v}$. By \cite[Section 3.3, especially, Lemme 3.3.1]{Ngo10}, the functor $\cP_{a,\, v}$ is representable by an ind-scheme whose associated reduced $\cP^\red_{a,\, v}$ is a locally of finite type $\ov{k}$-scheme. Due to its group structure, $\cP^\red_{a,\, v}$ is even $\ov{k}$-smooth, and it inherits an action on the locally of finite type $\ov{k}$-scheme $\cM^\red_{\sL^{\tensor 2},\, a,\, v}$.
\epp

\bpp[The product formula morphism] \label{pp:product-morphism}
Assume the setting of \S\S\ref{pp:Hitchin-fibration}--\ref{pp:symmetries} with $G$ quasi-split and $\#W$ invertible on $S$, in particular, fix a point $a \in \cA^{\heart}_{\sL^{\tensor 2}}(\ov{k})$ valued in an algebraically closed field, so that $X_a^\rs \subset X_{\ov{k}}$ is a dense open, and let $U_a \subset X_a^\rs$ be a dense open. Due to the moduli interpretation of $\cM_{\sL^{\tensor 2},\, a}$ and $\cM_{\sL^{\tensor 2},\, a,\, v}$, 
Beauville--Laszlo glueing in the style of \Cref{BL-glue}
of the Kostant--Hitchin section $\eps_{\sL}(a)|_{U_a}$ to sections of affine Springer fibers at the points in $X_{\ov{k}} \setminus U_a$ gives a $\ov{k}$-stack morphism 
\be \label{eqn:first-glue}
\tst \prod_{v\in X_{\ov{k}}\setminus U_a}\cM_{\sL^{\tensor 2},\, a,\, v}\rightarrow \cM_{\sL^{\tensor 2},\, a}.
\ee
A similar glueing gives an action of $\prod_{v\in X_{\ov{k}} \setminus U_a}\cP_{a,\,v}$ on $\cP_{a}$. Thus, by twisting the glueing that gives the morphism \eqref{eqn:first-glue} by variable $J_a$-torsors as in \S\ref{pp:symmetries}, we obtain the $\ov{k}$-stack~morphism
\be \label{eqn:product-quotient}
\tst \prod_{v\in X_{\ov{k}} \setminus U_a}\cM_{\sL^{\tensor 2},\, a,\, v}\times^{\prod_{v\in X_{\ov{k}} \setminus U_a}\cP_{a,\,v}}\cP_{a}\rightarrow\cM_{\sL^{\tensor 2},\, a}
\ee
whose source is the stackification of the prestack quotient that is described in general\footnote{\label{foot:explicit-quotient}We recall that for a $1$-category $\sX$ and a group $G$ acting on $\sX$, the quotient $\sX/G$ is the $1$-category whose objects are those of $\sX$ and morphisms between objects $x$ and $x'$ are given by pairs $(g, \iota)$ with $g \in G$ and $\iota \in \Hom_\sX(gx, x')$. The source of \eqref{eqn:product-quotient} is the stackification of this construction performed on groupoids of sections.} in \cite{Ngo06}*{before Lemme 4.7}, see also \cite{Rom05}*{proof of Proposition 2.6} (since $\cM_{\sL^{\tensor 2},\, a,\, v}$ and $\cP_{a,\,v}$ are only ind-schemes, we do not claim any algebraicity for this source). The map \eqref{eqn:product-quotient} is fully faithful: since its target is already a stack and stackification sheafifies the morphism functors \cite{SP}*{Lemma~\href{https://stacks.math.columbia.edu/tag/02ZN}{02ZN}}, it suffices to see this on points valued in a $\ov{k}$-algebra $R$ before the stackification, and then we use the Beauville--Laszlo glueing as follows. An isomorphism in $\cM_{\sL^{\tensor 2},\, a}$ between the glueings of 
\[
((m_v), p) \qxq{and} ((m'_v), p')
\]
amounts to both a $J_a$-torsor isomorphism $p|_{U_a} \isomto p'|_{U_a}$ (see the gerbe aspect of \Cref{rem:J-with-an-L}), which, after twisting by uniquely determined $(p_v)$ in $\prod_{v\in X_{\ov{k}} \setminus U_a}\cP_{a,\,v}$ (the difference between the $J_a$-torsors $p'$ and $p$), extends to a $J_a$-torsor isomorphism $p \isomto p'$, and, granted this uniquely determined adjustment, isomorphisms 
\[
m_v \isomto m'_v \qxq{for} v \in X_{\ov{k}}\setminus U_a.
\] 
To counter the potentially nonalgebraic nature of the source of \eqref{eqn:product-quotient}, one considers the following variant. As reviewed in \S\S\ref{pp:affine-Springer-fibration}--\ref{pp:symmetries}, both $\cM_{\sL^{\tensor 2},\, a,\, v}^\red$ and $\cP_{a,\,v}^\red$ are locally of finite type $\ov{k}$-schemes, with $\cP_{a,\,v}^\red$ even a smooth $\ov{k}$-group that acts on $\cM_{\sL^{\tensor 2},\, a,\, v}^\red$. The smoothness 
ensures that the stack 
\[
\tst \prod_{v\in X_{\ov{k}} \setminus U_a}\cM_{\sL^{\tensor 2},\, a,\, v}^\red \times^{\prod_{v\in X_{\ov{k}} \setminus U_a}\cP_{a,\,v}^\red}\cP_{a}
\]
is algebraic and, by \cite{SP}*{Lemma \href{https://stacks.math.columbia.edu/tag/076V}{076V}} and \cite{Rom05}*{Theorem 4.1 and its proof}, may be formed in the \'{e}tale topology. Consequently, the explicit description of the quotients before the stackification reviewed in \cref{foot:explicit-quotient} and the agreement of $\cM_{\sL^{\tensor 2},\, a,\, v}^\red$ and $\cM_{\sL^{\tensor 2},\, a,\, v}$ (respectively, of $\cP_{a,\,v}^\red$ and $\cP_{a,\,v}$) on reduced rings implies that the morphism 
\be \label{eqn:now-reduced-quotient}
\tst \prod_{v\in X_{\ov{k}} \setminus U_a}\cM_{\sL^{\tensor 2},\, a,\, v}^\red \times^{\prod_{v\in X_{\ov{k}} \setminus U_a}\cP_{a,\,v}^\red}\cP_{a} \ra \prod_{v\in X_{\ov{k}} \setminus U_a}\cM_{\sL^{\tensor 2},\, a,\, v}\times^{\prod_{v\in X_{\ov{k}} \setminus U_a}\cP_{a,\,v}}\cP_{a}
\ee
is an equivalence on $R$-points for every \emph{reduced} $\ov{k}$-algebra $R$. The composition of \eqref{eqn:product-quotient} and \eqref{eqn:now-reduced-quotient} is the promised product formula morphism between locally of finite type algebraic $\ov{k}$-stacks:
\be \label{eqn:product-formula-morphism}
\tst \prod_{v\in X_{\ov{k}} \setminus U_a}\cM_{\sL^{\tensor 2},\, a,\, v}^\red \times^{\prod_{v\in X_{\ov{k}} \setminus U_a}\cP_{a,\,v}^\red}\cP_{a} \ra \cM_{\sL^{\tensor 2},\, a}.
\ee
By the above, this map is fully faithful on groupoids of $R$-points for every reduced $\ov{k}$-algebra $R$. 
\epp

The proof of the product formula in \Cref{thm:product-formula} will rely on the following general lemma, whose argument shows that the locally of finite type over a field assumption could be weakened significantly.

\begin{lem}\label{norm}
Let $f\colon \cX \ra \cY$ be a map  of algebraic stacks that are locally of finite type over a field $k$. If $f(R)$ is an equivalence for normal, strictly Henselian, local $k$-algebras $R$, then $f$ is a universal~homeomorphism.
\end{lem}

\begin{proof}
The assumption about $R$ is stable under base change along any $\cY' \ra \cY$. Thus, by passing to a smooth cover of $\cY$ and using the fppf local on the base nature of being a universal homeomorphism \cite{SP}*{Lemma \href{https://stacks.math.columbia.edu/tag/0DTQ}{0DTQ}}, we may assume that $\cY$ is an affine scheme. Moreover, by letting $R$ be a field, we see that the map $\cX \ra \cY$ induces a continuous bijection on the underlying topological spaces defined in \cite{SP}*{Lemma \href{https://stacks.math.columbia.edu/tag/04XL}{04XL}, Definition \href{https://stacks.math.columbia.edu/tag/04XG}{04XG}}, and continues to do so after any base change. 
Thus, we only need to argue that it is also universally closed. For this, we consider the normalization morphism $\wt{\cY} \ra \cY$, which, due to our setup over $k$ and \cite{SP}*{Lemmas \href{https://stacks.math.columbia.edu/tag/035Q}{035Q} and \href{https://stacks.math.columbia.edu/tag/035S}{035S}}, is finite and surjective. These properties are preserved after base change along $\cX \ra \cY$, so, for proving the remaining universal closedness of $\cX \ra \cY$, we may pass to $\cX \times_\cY \wt{\cY} \ra \wt{\cY}$ and assume that $\cY$ is normal. 
The assumption now implies that the morphisms $\Spec(\sO_{\cY,\, y}^\sh) \ra \cY$ lift to $\cX$, so, since the latter is locally of finite presentation over $k$, we conclude from a limit argument that the map $\cX \ra \cY$ has a section \'{e}tale locally on $\cY$. In particular, after base change to some \'{e}tale cover of $\cY$, by \cite{SP}*{Lemma~\href{https://stacks.math.columbia.edu/tag/0DTQ}{0DTQ}} again, we may assume that our universal continuous bijection $\cX \ra \cY$ has a section, a case in which it certainly is a universal homeomorphism.
\end{proof}

\begin{theorem} \label{thm:product-formula}
Let $S$ be a scheme, let $X$ be a proper, smooth $S$-scheme with geometrically connected fibers of dimension $1$, let $\sL$ be a $\bG_m$-torsor on $X$, let $G$ be a quasi-split reductive group over $X$, and let 
$\ov{k}$ be an algebraically closed $S$-field  in which the order of the Weyl group of $G$ is invertible.
\benum
\m \label{PF-a}
For an $a\in\cA^{\heartsuit}_{\sL^{\tensor 2}}(\overline{k})$ and a dense open $U_a \subset X_a^\rs$, the morphism of locally of finite type algebraic $\ov{k}$-stacks 
\be \label{eqn:comparison-map}
\tst \prod_{v\in X_{\ov{k}} \setminus U_a } \cM^{\red}_{\sL^{\tensor 2},\, a,\, v}\times^{\prod_{v\in X_{\ov{k}} \setminus U_a}\cP^{\red}_{a,v}}\cP_a \xra{\eqref{eqn:product-formula-morphism}} \cM_{\sL^{\tensor 2},\, a}
\ee
 is a universal homeomorphism that induces an equivalence on the groupoids of $R$-points for every seminormal, strictly Henselian, local $\ov{k}$-algebra $R$. 

\m \label{PF-b}
If $\ov{k}$ is the algebraic closure of a finite field\footnote{The only purpose of the finite field assumption is to be able to apply results from \cite{Ngo10}, especially, \cite{Ngo10}*{Proposition 4.15.1}.} and $a \in \cA^{\heartsuit}_{\sL^{\tensor 2}}(\overline{k})$ is anisotropic in the sense that the set of connected components $\pi_0(\cP_a)$ is finite, then the universal homeomorphism \eqref{eqn:comparison-map} is a finite morphism that is representable by schemes. 
\eenum
\end{theorem}

\begin{proof}
\Cref{norm} reduces the universal homeomorphism aspect to the claim about $R$-points. For the latter, since a seminormal $R$ is reduced (see \S\ref{conv}), we already know from \S\ref{pp:product-morphism} that the map \eqref{eqn:comparison-map} is fully faithful on $R$-points. For the essential surjectivity, fix an $m\in\cM_{\sL^{\tensor 2},\,a}(R)$. By the gerbe aspect of \Cref{rem:J-with-an-L} and the agreement of $\fg_{\sL^{\tensor 2}}$ and $\fg^\reg_{\sL^{\tensor 2}}$ over $\ft^\rs_\sL/W$ ensured by \Cref{prop:describe-Chevalley}, the restrictions $m|_{U_a}$ and $\eps_{\sL}(a)|_{U_a}$ differ by a uniquely determined $(J_a|_{U_a})$-torsor. Thus, it suffices to show that this $(J_a|_{U_a})$-torsor extends to a $J_a$-torsor over the entire $X_{\ov{k}}$---we would then be able to absorb it together with $m$ into the source of \eqref{eqn:comparison-map}. By Beauville--Laszlo glueing (see \Cref{BL-glue}~\ref{BLG-b}), restriction to the punctured formal neighborhoods $R\llp t_v \rrp$ of $(X_{\ov{k}})_R$ along the points $v \in X_{\ov{k}} \setminus U_a$ then reduces us to arguing the triviality of $J_a$-torsors over such $R\llp t_v \rrp$. However, the Galois-theoretic description supplied by \Cref{prop:galois-theoretic} implies that $J_a|_{R\llp t_v \rrp}$ is a torus that trivializes over some $W$-torsor. Since $\#W$ is invertible in $R$, \Cref{twisted} ensures that 
\[
H^1(R\llp t_v \rrp, J_a) = 0,
\]
and it follows that $m$ is in the essential image, as desired. 

Assume now the setup of \ref{PF-b}, so that, by \cite{Ngo10}*{Proposition 4.15.1}, the source of \eqref{eqn:comparison-map} is Deligne--Mumford with quasi-compact and separated diagonal. By \S\ref{pp:Hitchin-fibration}, the diagonal of the target of \eqref{eqn:comparison-map} is even affine. Thus, since the map \eqref{eqn:comparison-map} induces an equivalence on points valued in algebraically closed fields, by the criterion \cite{modular-description}*{Lemma 3.2.2~(b)} for which in the Deligne--Mumford case it suffices to consider field-valued points as indicated there, it is representable by algebraic spaces. By bootstrap \cite{SP}*{Lemma \href{https://stacks.math.columbia.edu/tag/050N}{050N}}, this map is also quasi-separated, so the valuative criterion \cite{SP}*{Lemma \href{https://stacks.math.columbia.edu/tag/03KV}{03KV}} and the first part of the claim for $R$ a valuation ring imply that \eqref{eqn:comparison-map} is separated. However, by \cite{Ryd10}*{Corollary 5.22}, a separated universal homeomorphism of algebraic spaces is representable by schemes. Thus, since it is also locally of finite type and, by \cite{SP}*{Lemma~\href{https://stacks.math.columbia.edu/tag/04DF}{04DF}}, integral, it is necessarily a finite morphism. 
\end{proof}


\end{document}